\documentclass[12pt]{article}
\usepackage[cp850]{inputenc}
\usepackage{latexsym}
\usepackage{amssymb}
\parskip=\medskipamount
\newtheorem{guia}{}[section]
\newtheorem{lemma}[guia]{Lemma} 
\newtheorem{example}[guia]{Example}
\newtheorem{examples}[guia]{Examples}
\newtheorem{definition}[guia]{Definition}
\newtheorem{proposition}[guia]{Proposition}
\newtheorem{corollary}[guia]{Corollary}
\newtheorem{theorem}[guia]{Theorem}
\newtheorem{remark}[guia]{Remark}

 1   
\font\ddpp=msbm10  scaled \magstep 1  

\newenvironment{proof}{\noindent{\bf Proof~:}}{\QED\medskip}
\def\QED{\hskip0.1em\hfill\null\ \null\nobreak\hfill
\kern3pt\lower1.8pt\vbox{\hrule\hbox
{\vrule\kern1pt\vbox{\kern1.7pt \hbox{$\scriptstyle
QED$}\kern0.2pt}\kern1pt\vrule}\hrule}}
\def\R{\hbox{\ddpp R}}               
\def\N{\hbox{\ddpp N}}    

\def\arco#1{\kern3pt \mathop{\vbox{\ialign{##\crcr\noalign{\kern1pt}
        $\braceld\leaders\vrule\hfill\leaders\vrule\hfill\bracerd$
    \crcr\noalign{\kern1pt\nointerlineskip}
        $\hfil\displaystyle{\kern-1pt#1\kern2pt}\hfil$\crcr}}}\limits}

\newcommand\map[3]{#1\ \colon\ #2\longrightarrow#3}
\def\lcf{\lbrack\! \lbrack}
\def\rcf{\rbrack\! \rbrack}
\newcommand{\br}[3][]{\lbrack\!\lbrack#2,#3\rbrack\!\rbrack^{#1}}
\newcommand{\an}[1][]{\rho^{#1}}
\renewcommand{\d}[1][]{d^{#1}}

\newcommand{\prd}{\tau^*}
\newcommand{\Prol}[2]{\mathcal{L}^{#1}#2}
\newcommand{\prol}[1][\tau]{\Prol{#1}{E}}
\newcommand{\prold}[1][\tau^*]{\Prol{#1}{E}}

\newcommand{\FL}{{Leg}_L}
\newcommand{\text}[1]{\mbox{\rm #1}}
\newcommand{\qquand}{\qquad\text{and}\qquad}

\newcommand{\cinfty}[1]{C^\infty(#1)}
\newcommand{\pd}[2]{\frac{\partial #1}{\partial #2}}
\renewcommand{\sec}[1]{\Gamma(#1)}
\begin{document}

\baselineskip=.55cm
\title{\bf Lagrangian submanifolds and dynamics on Lie algebroids
\footnotetext{$\mbox{}$ \kern-20pt{\it
Mathematics Subject Classification} (2000): 17B66, 53D12, 70G45, 70H03, 70H05, 70H20. \\
\noindent {\it Key words and phrases}: Lagrangian submanifolds,
Euler-Lagrange equations, Hamilton equations, Lie algebroids,
canonical involution, Tulczyjew's triple, symplectic Lie
algebroids, Atiyah algebroid, Lagrange-Poincar\'{e} equations,
Hamilton-Poincar\'{e} equations.} }
\author{\large Manuel  de Le\'on$^1$, Juan C.
Marrero$^{2}$, Eduardo Mart\'{\i}nez$^3$
\\[7pt]
{\small \it $^1$Instituto de Matem\'aticas y F{\'\i}sica
Fundamental,}\\[-8pt] {\small\it Consejo Superior de
In\-ves\-ti\-ga\-cio\-nes Ci\-en\-t{\'\i}\-fi\-cas,} \\[-8pt]
{\small\it Serrano 123, 28006 Madrid, SPAIN,}\\[-8pt] {\small\it
E-mail: mdeleon@imaff.cfmac.csic.es} \\ {\small\it
$^2$Departamento de Matem\'atica Fundamental, Facultad de
Matem\'aticas,}\\[-8pt] {\small\it Universidad de la Laguna, La
Laguna,} \\[-8pt] {\small\it Tenerife, Canary Islands,
SPAIN,}\\[-8pt] {\small\it E-mail: jcmarrer@ull.es } \\ {\small\it
$^3$Departamento de Matem\'atica Aplicada, Facultad de
Ciencias,}\\[-8pt] {\small \it Universidad de Zaragoza,}\\[-8pt]
{\small \it 50009 Zaragoza, SPAIN,}\\[-8pt] {\small \it E-mail:
emf@unizar.es}}
\date{\empty}

\maketitle

\begin{abstract}
{\footnotesize In some previous papers, a geometric description of
Lagrangian Mechanics on Lie algebroids has been developed. In the
present paper, we give a Hamiltonian description of Mechanics on
Lie algebroids. In addition, we introduce the notion of a
Lagrangian submanifold of a symplectic Lie algebroid and we prove
that the Lagrangian (Hamiltonian) dynamics on Lie algebroids  may
be described in terms of Lagrangian submanifolds of symplectic Lie
algebroids. The Lagrangian (Hamiltonian) formalism on Lie
algebroids permits to deal with Lagrangian (Hamiltonian) functions
not defined necessarily on tangent (cotangent) bundles. Thus, we
may apply our results to the projection of Lagrangian
(Hamiltonian) functions which are invariant under the action of a
symmetry Lie group. As a consequence, we obtain that
Lagrange-Poincar\'{e} (Hamilton-Poincar\'{e}) equations are the
Euler-Lagrange (Hamilton) equations associated with the
corresponding Atiyah algebroid. Moreover, we prove that
Lagrange-Poincar\'{e} (Hamilton-Poincar\'{e}) equations are the local
equations defining certain Lagrangian submanifolds of symplectic
Atiyah algebroids.}
\end{abstract}

\baselineskip=.521cm

\newpage

\tableofcontents

\section{Introduction}
\setcounter{equation}{0}

 Lie algebroids have deserved
a lot of interest in recent years. Since a Lie algebroid is a
concept which unifies tangent bundles and Lie algebras, one can
suspect their relation with Mechanics. In his paper
\cite{weinstein} A. Weinstein (see also the paper by P. Libermann
\cite{liber}), developed a generalized theory of Lagrangian
Mechanics on Lie algebroids  and obtained the equations of motion,
using the linear Poisson structure on the dual of the Lie
algebroid  and the Legendre transformation associated with the
Lagrangian $L$, when $L$ is regular. In that paper, he also asks
the question of whether it is possible to develop a formalism
similar on Lie algebroids to Klein's formalism \cite{Klein} in
ordinary Lagrangian Mechanics. This task was finally done by E.
Mart{\'\i}nez \cite{M} (see also \cite{CaM,CM,Ma2,Ma3}). The main
notion is that of prolongation of a Lie algebroid over a mapping,
introduced by P.J. Higgins and K. Mackenzie \cite{HM}. The purpose
of the present paper is to give a description of Hamiltonian and
Lagrangian dynamics on Lie algebroids in terms of Lagrangian
submanifolds of symplectic Lie algebroids.

One could ask about the interest to generalize Classical Mechanics
on tangent and cotangent bundles to Lie algebroids. However, it is
not a mere academic exercise. Indeed, if we apply our procedure to
Atiyah algebroids we recover in a very natural way the
Lagrange-Poincar\'e and Hamilton-Poincar\'e equations. In this
case, the Lagrangian and Hamiltonian functions are not defined on
tangent and cotangent bundles, but on the quotients by the
structure Lie group. This fact is a good motivation for our study.

The paper is organized as follows. In Section 2.1 we recall the
notion of prolongation ${\cal L}^{f}E$ of a Lie algebroid $\tau :
E \longrightarrow M$ over a mapping $f : M' \longrightarrow M$;
when $f$ is just the canonical projection $\tau$, then ${\cal
L}^{\tau}E$ will play the role of the tangent bundle. We also
consider action Lie algebroids, which permit to induce a Lie
algebroid structure on the pull-back of a Lie algebroid by a
mapping. The notion of quotient Lie algebroids is also discussed,
and in particular, we consider Atiyah algebroids. In Section 2.2
we develop the Lagrangian formalism on the prolongation ${\cal
L}^{\tau}E$ starting with a Lagrangian function $L : E
\longrightarrow \R$. Indeed, one can construct the
Poincar\'e-Cartan 1 and 2-sections (i.e. $\theta_L \in
\Gamma(({\cal L}^{\tau}E)^*)$ and $\omega_L \in
\Gamma(\Lambda^2({\cal L}^{\tau}E)^*)$, respectively) using the
geometry of ${\cal L}^{\tau}E$ provided by the Euler section
$\Delta$ and the vertical endomorphism $S$. The dynamics is given
by a SODE $\xi$ of ${\cal L}^{\tau}E$ (that is, a section $\xi$ of
${\cal L}^{\tau}E$ such that $S\xi = \Delta$) satisfying $i_{\xi}
\, \omega_L = d^{{\cal L}^{\tau}E}E_L$, where $E_L$ is the energy
associated with $L$ (along the paper $d^E$ denotes the
differential of the Lie algebroid $E$). As in Classical Mechanics,
$L$ is regular if and only if $\omega_L$ is a symplectic section,
and in this case $\xi=\xi_L$ is uniquely defined and a SODE. Its
solutions (curves in $E$) satisfy the Euler-Lagrange equations for
$L$.

Sections 3.1-3.4 are devoted to develop a Hamiltonian description
of Mechanics on Lie algebroids. Now, the role of the cotangent
bundle of the configuration manifold is played by the prolongation
${\cal L}^{\tau^*}E$ of $E$ along the projection $\tau^* : E^*
\longrightarrow M$, which is the dual bundle of $E$. We can
construct the canonical Liouville 1-section $\lambda_E$ and the
canonical symplectic 2-section $\Omega_E$ on ${\cal L}^{\tau^*}E$.
Theorem 3.4 and Corollary 3.6 are the Lie algebroid version of the
classical results concerning the universality of the standard
Liouville 1-form on cotangent bundles. Given a Hamiltonian
function $H : E^* \longrightarrow \R$, the dynamics are obtained
solving the equation $i_{\xi_H} \, \Omega_E = d^{{\cal
L}^{\tau^*}E}H$ with the usual notations. The solutions of $\xi_H$
(curves in $E^*$) are the ones of the Hamilton equations for $H$.
The Legendre transformation $Leg_L : E \longrightarrow E^*$
associated with a Lagrangian $L$ induces a Lie algebroid morphism
${\cal L}Leg_L : {\cal L}^{\tau}E \longrightarrow {\cal
L}^{\tau^*}E$, which permits in the regular case to connect
Lagrangian and Hamiltonian formalisms as in Classical Mechanics.
In Section 3.5 we develop the corresponding Hamilton-Jacobi
theory; we prove that the function $S : M \longrightarrow \R$
satisfying the Hamilton-Jacobi equation $d^E(H \circ d^ES) = 0$ is
just the action for $L$.

As is well-known, there is a canonical involution $\sigma_{TM} :
TTM \longrightarrow TTM$ defined by S. Kobayashi (\cite{koba}). In
Section 4 we prove that for an arbitrary Lie algebroid $\tau : E
\longrightarrow M$ there is a unique Lie algebroid isomorphism
$\sigma_E : {\cal L}^{\tau}E \longrightarrow \rho^*(TE)$ such that
$\sigma_E^2 = id$, where ${\cal L}^{\tau}E$ is the prolongation of
$E$ by $\tau$, and $\rho^*(TE)$ is the pull-back of the tangent
bundle prolongation $T\tau : TE \longrightarrow TM$ via the anchor
mapping $\rho : E \longrightarrow TM$ (Theorem 4.4). Note that
$\rho^*(TE) = {\cal L}^{\tau}E$. We also remark that $\rho^*(TE)$
carries a Lie algebroid structure over $E$ since the existence of
an action of the tangent Lie algebroid $T\tau : TE \longrightarrow
TM$ on $\tau$. When $E$ is the standard Lie algebroid $TM$ we
recover the standard canonical involution.

In \cite{Tul1,Tul2} W.M. Tulczyjew has interpreted the Lagrangian
and Hamiltonian dynamics as Lagrangian submanifolds of convenient
special symplectic manifolds. To do that, Tulczyjew introduced
canonical isomorphisms which commute the tangent and cotangent
functors. Section 5 is devoted to extend Tulczyjew' s
construction. First we define a canonical vector bundle
isomorphism $\flat_{E^*} : {\cal L}^{\tau^*}E \longrightarrow
({\cal L}^{\tau^*}E)^*$ which is given using the canonical
symplectic section of ${\cal L}^{\tau^*}E$. Next, using the
canonical involution $\sigma_E$ one defines a canonical vector
bundle isomorphism $A_E : {\cal L}^{\tau^*}E \longrightarrow
({\cal L}^{\tau}E)^*$. Both vector bundle isomorphims extend the
so-called Tulczyjew's triple for Classical Mechanics.

In Section 6 we introduce the notion of a symplectic Lie
algebroid. The definition is the obvious one: $\Omega$ is a
symplectic section on a Lie algebroid $\tau :E \longrightarrow M$
if it induces a nondegenerate bilinear form on each fibre of $E$
and, in addition, it is $d^E$-closed ($d^E \Omega = 0$). In this
case, the prolongation ${\cal L}^{\tau}E$ is symplectic too. The
latter result extends the well-known result which proves that the
tangent bundle of a symplectic manifold is also symplectic.

In Section 7 we consider Lagrangian Lie subalgebroids of
symplectic Lie algebroids; the definition is of course pointwise.
This definition permits to consider, in Section 8, the notion of a
Lagrangian submanifold of a symplectic Lie algebroid: a
submanifold $i : S \longrightarrow E$ is a Lagrangian submanifold
of the symplectic Lie algebroid $\tau : E \longrightarrow M$ with
anchor map $\rho : E \longrightarrow TM$ if the following
conditions hold:

\begin{itemize}
\item $\dim (\rho(E_{{\tau^S}(x)}) + (T_x \tau^S)(T_xS))$ does not depend on $x$, for all $x
\in S$;
\item the Lie subalgebroid ${\cal L}^{\tau^{S}}E$ of the symplectic Lie algebroid ${\cal L}^\tau E$, is
Lagrangian;
\end{itemize}
here $\tau^{S} = \tau \circ i : S \longrightarrow M$. The
classical results about Lagrangian submanifolds in symplectic
geometry are extended to the present context in the natural way.
Also, we generalize the interpretation of Tulzcyjew; for instance,
given a Hamiltonian $H : E^* \longrightarrow \R$ we prove that
$S_{H} = \xi_H(E^*)$ is a Lagrangian submanifold of the symplectic
extension ${\cal L}^{\tau^*}E$ and that there exists a bijective
correspondence between the admissible curves in $S_{H}$ and the
solutions of the Hamilton equations for $H$. For a Lagrangian $L :
E \longrightarrow \R$ we prove that $S_L = (A_E^{-1} \circ
d^{{\cal L}^{\tau}E}L)(E)$ is a Lagrangian submanifold of ${\cal
L}^{\tau^*}E$ and, furthermore, that there exists a bijective
correspondence between the admissible curves in $S_{L}$ and the
solutions of the Euler-Lagrange equations for $L$. In addition, we
deduce that for an hyperregular Lagrangian $L$, then $S_{\xi_L}=
\xi_L(E)$ is a Lagrangian submanifold of the symplectic extension
${\cal L}^\tau E$, and, moreover, we have that ${\cal
L}Leg_L(S_{\xi_L}) = S_L = S_H$, for $H = E_L \circ Leg_L^{-1}$.

Finally, Sections 9.1-9.4 are devoted to some applications.
Consider a (left) principal bundle $\pi : Q \longrightarrow M$
with structural group $G$. The Lie algebroid $\tau_{Q}|G : TQ/G
\longrightarrow M$ is called the Atiyah algebroid associated with
$\pi : Q \longrightarrow M$. One can prove that the prolongation
${\cal L}^{\tau_{Q}|G}(TQ/G)$ is isomorphic to the Atiyah
algebroid associated with the principal bundle $\pi_T : TQ
\longrightarrow TQ/G$, and, moreover, the dual vector bundle of
${\cal L}^{\tau_{Q}|G}(TQ/G)$ is isomorphic to the quotient vector
bundle of $\pi_{TQ}: T^*(TQ) \longrightarrow TQ$ by the canonical
lift action of $G$ on $T^*(TQ)$. Similar results are obtained for
cotangent bundles. In Section 9.2 (respectively, Section 9.3), we
prove that the solutions of the Hamilton-Poincar\'e equations for
a $G$-invariant Hamiltonian function $H : T^*Q \longrightarrow \R$
(resp., the Lagrange-Poincar\'e equations for a $G$-invariant
Lagrangian $L : TQ \longrightarrow \R$) are just the solutions of
the Hamilton equations (resp. the Euler-Lagrange equations) on
$T^*Q/G$ for the reduced Hamiltonian $h : T^*Q/G \longrightarrow
\R$ (resp. on $TQ/G$ for the reduced Lagrangian $l : TQ/G
\longrightarrow \R$). Moreover, in Sections 9.2 and 9.3 all these
equations are reinterpreted as those defining the corresponding
Lagrangian submanifolds. Finally, in Section 9.4 we show how our
formalism allows us to obtain in a direct way Wong' s equations.

Manifolds are real, paracompact and $C^{\infty}$. Maps are
$C^{\infty}$. Sum over crossed repeated indices is understood.

\section{Lie algebroids and Lagrangian Mechanics}
\subsection{Some algebraic constructions in the category of Lie
algebroids}

Let $E$ be a vector bundle of $rank$ $n$ over a manifold $M$ of
dimension $m$  and $\tau:E\to M$ be the vector bundle projection.
Denote by $\Gamma(E)$ the $C^\infty(M)$-module of sections of
$\tau:E\to M$. A {\it Lie algebroid structure }
$(\lcf\cdot,\cdot\rcf,\rho)$ on $E$ is a Lie bracket
$\lcf\cdot,\cdot\rcf$ on the space $\Gamma(E)$ and a bundle map
$\rho:E\to TM$, called {\it the anchor map}, such that if we also
denote by $\rho:\Gamma(E)\to {\frak X}(M)$ the homomorphism of
$C^\infty(M)$-modules induced by the anchor map then
\[
\lcf X,fY\rcf=f\lcf X,Y\rcf + \rho(X)(f)Y,
\]
for $X,Y\in \Gamma(E)$ and $f\in C^\infty(M)$. The triple
$(E,\lcf\cdot,\cdot\rcf,\rho)$ is called {\it a Lie algebroid
over} $M$ (see \cite{Ma}).

If $(E,\lcf\cdot,\cdot\rcf,\rho)$ is a Lie algebroid over $M,$
then the anchor map $\rho:\Gamma(E)\to {\frak X}(M)$ is a
homomorphism between the Lie algebras
$(\Gamma(E),\lcf\cdot,\cdot\rcf)$ and $({\frak
X}(M),[\cdot,\cdot])$.

Trivial examples of Lie algebroids are real Lie algebras of finite
dimension and the tangent bundle $TM$ of an arbitrary manifold
$M.$

Let $(E,\lcf\cdot,\cdot\rcf,\rho)$ be a Lie algebroid over $M$. We
consider the generalized distribution ${\cal F}^E$ on $M$ whose
characteristic space at a point $x\in M$ is given by
\[
{\cal F}^E(x)=\rho(E_x)
\]
where $E_x$ is the fiber of $E$ over $x.$ The distribution ${\cal
F}^E$ is finitely generated and involutive. Thus, ${\cal F}^E$
defines a generalized foliation on $M$ in the sense of Sussmann
\cite{Su}. ${\cal F}^E$ is {\it the Lie algebroid foliation } on
$M$ associated with $E$.

If $(E,\lcf\cdot,\cdot\rcf,\rho)$ is a Lie algebroid, one may
define {\it the differential of $E$}, $d^E:\Gamma(\wedge^k E^*)\to
\Gamma(\wedge^{k+1}E^*)$, as  follows
\begin{equation}\label{dif}
\begin{array}{rcl}
d^E \mu(X_0,\cdots, X_k)&=&\displaystyle \sum_{i=0}^{k}
(-1)^i\rho(X_i)(\mu(X_0,\cdots,
\widehat{X_i},\cdots, X_k)) \\[5pt] &&+ \displaystyle\sum_{i<
j}(-1)^{i+j}\mu(\lcf X_i,X_j\rcf,X_0,\cdots,
\widehat{X_i},\cdots,\widehat{X_j},\cdots ,X_k),
\end{array}
\end{equation}
for $\mu\in \Gamma(\wedge^k E^*)$ and $X_0,\dots ,X_k\in
\Gamma(E).$ It follows that $(d^E)^2=0$. Moreover, if $X$ is a
section of $E$, one may introduce, in a natural way, {\it the Lie
derivative with respect to $X$,} as the operator ${\cal
L}^E_X:\Gamma(\wedge^kE^*)\to \Gamma(\wedge^k E^*)$ given by
\[
{\cal L}^E_X=i_X\circ d^E + d^E\circ i_X.
\]

Note that if $E = TM$ and $X \in \Gamma(E) = {\frak X}(M)$ then
$d^{TM}$ and ${\cal L}_{X}^{TM}$ are the usual differential and
the usual Lie derivative with respect to $X$, respectively.

If we take local coordinates $(x^i)$ on $M$ and a local basis
$\{e_\alpha\}$ of sections of $E$, then we have the corresponding
local coordinates $(x^i,y^\alpha)$ on $E$, where $y^\alpha(a)$ is
the $\alpha$-th coordinate of $a\in E$ in the given basis. Such
coordinates determine local functions $\rho_\alpha^i$,
$C_{\alpha\beta}^{\gamma}$ on $M$ which contain the local
information of the Lie algebroid structure, and accordingly they
are called the {\it structure functions of the Lie algebroid.}
They are given by
\[
\rho(e_\alpha)=\rho_\alpha^i\frac{\partial }{\partial
x^i}\;\;\;\mbox{ and }\;\;\; \lcf
e_\alpha,e_\beta\rcf=C_{\alpha\beta}^\gamma e_\gamma.
\]
These functions should satisfy the relations
\begin{equation}\label{estruc1}
\rho_\alpha^j\frac{\partial \rho_\beta^i}{\partial x^j}
-\rho_\beta^j\frac{\partial \rho_\alpha^i}{\partial x^j}=
\rho_\gamma^iC_{\alpha\beta}^\gamma \end{equation}

and

\begin{equation}\label{estruc2}
\sum_{cyclic(\alpha,\beta,\gamma)}[\rho_{\alpha}^i\frac{\partial
C_{\beta\gamma}^\nu}{\partial x^i} + C_{\alpha\mu}^\nu
C_{\beta\gamma}^\mu]=0,
\end{equation}
which are usually called {\it the structure equations.}

If $f\in C^\infty(M)$, we have that
\begin{equation}\label{diff0}
d^E f=\frac{\partial f}{\partial x^i}\rho_\alpha^i e^\alpha,
\end{equation}
where $\{e^\alpha\}$ is the dual basis of $\{e_\alpha\}$. On the
other hand, if $\theta\in \Gamma(E^*)$ and $\theta=\theta_\gamma
e^\gamma$ it follows that
\begin{equation}\label{diff1}
d^E \theta=(\frac{\partial \theta_\gamma}{\partial
x^i}\rho^i_\beta-\frac{1}{2}\theta_\alpha
C^\alpha_{\beta\gamma})e^{\beta}\wedge e^\gamma.
\end{equation}
In particular, \[d^E x^i=\rho_\alpha^ie^\alpha,\;\;\; d^E
e^\alpha=-\frac{1}{2}C_{\beta\gamma}^{\alpha} e^\beta\wedge
e^\gamma.\]

On the other hand, if $(E,\lcf\cdot,\cdot \rcf,\rho)$ and
$(E',\lcf\cdot,\cdot\rcf', \rho')$ are Lie algebroids over $M$ and
$M'$, respectively, then a morphism of vector bundles $(F,f)$ of
$E$ on $E'$

\begin{picture}(375,90)(40,20)
\put(180,20){\makebox(0,0){$M$}}
\put(250,25){$f$}\put(210,20){\vector(1,0){80}}
\put(310,20){\makebox(0,0){$M'$}} \put(170,50){$\tau$}
\put(180,70){\vector(0,-1){40}} \put(320,50){$\tau'$}
\put(310,70){\vector(0,-1){40}} \put(180,80){\makebox(0,0){$E$}}
\put(250,85){$F$}\put(210,80){\vector(1,0){80}}
\put(310,80){\makebox(0,0){$E'$}} \end{picture}

 is a Lie algebroid
morphism if
\begin{equation}\label{Morph}
d^E ((F,f)^*\phi')= (F, f)^*(d^{E'}\phi'), \;\;\; \mbox{ for
}\phi'\in \Gamma(\wedge^k(E')^*).
\end{equation}
Note that $(F, f)^*\phi'$ is the section of the vector bundle
$\wedge^kE^*\to M$ defined by
\[
((F,f)^*\phi')_x(a_1,\dots ,a_k)=\phi'_{f(x)}(F(a_1),\dots
,F(a_k)),
\]
for $x\in M$ and $a_1,\dots ,a_k\in E_{x}$. We remark that
(\ref{Morph}) holds if and only if
\begin{equation}\label{Morph1}
\begin{array}{l}
d^E(g'\circ f)=(F, f)^{*}(d^{E'}g'), \;\;\; \mbox{for }g'\in
C^\infty(M'),\\ d^E((F,f)^*\alpha')=(F,
f)^*(d^{E'}\alpha'),\;\;\;\mbox{for } \alpha'\in \Gamma((E')^*).
\end{array} \end{equation}

If $M=M'$ and $f=id:M\to M$ then, it is easy to prove that the
pair $(F,id)$ is a Lie algebroid morphism if and only  if
\[F\lcf X,Y\rcf=\lcf FX,FY \rcf',\;\;\; \rho'(FX)=\rho(X)\]
for $X,Y\in \Gamma(E).$

 Other equivalent
definitions of a Lie algebroid morphism may be found in \cite{HM}.

Let $(E,\lcf\cdot,\cdot\rcf,\rho)$ be a Lie algebroid over $M$ and
$E^*$ be the dual bundle to $E.$ Then, $E^*$ admits a linear
Poisson structure $\Lambda_{E^*}$, that is, $\Lambda_{E^*}$ is a
$2$-vector on $E^*$ such that
\[
[\Lambda_{E^*},\Lambda_{E^*}]=0
\]
and if $y$ and $y'$ are linear functions on $E^*,$ we have that
$\Lambda_{E^*}(dy,dy')$ is also linear function on $E^*$. If
$(x^i)$ are local coordinates on $M$, $\{e_\alpha\}$ is a local
basis of $\Gamma(E)$ and $(x^i,y_\alpha)$ are the corresponding
coordinates on $E^*$ then the local expression of $\Lambda_{E^*}$
is
\begin{equation}\label{2.7'}
\Lambda_{E^*}=\frac{1}{2} C_{\alpha\beta}^\gamma
y_\gamma\frac{\partial}{\partial y_\alpha}\wedge \frac{\partial
}{\partial y_\beta} + \rho_\alpha^i\frac{\partial }{\partial
y_\alpha}\wedge \frac{\partial }{\partial x^i},
\end{equation}
where $\rho_\alpha^i$ and $C_{\alpha\beta}^\gamma$ are the
structure functions of $E$ with respect to the coordinates $(x^i)$
and to the basis $\{e_\alpha\}$. The Poisson structure
$\Lambda_{E^*}$ induces a linear Poisson bracket of functions on
$E^*$ which we will denote by $\{ \; , \; \}_{E^*}$. In fact, if
$F, G \in C^{\infty}(E^*)$ then
\[
\{F, G\}_{E^*} = \Lambda_{E^*}(d^{TE^*}F, d^{TE^*}G).
\]
On the other hand, if $f$ is a function on $M$ then the associated
basic function $f^v\in\cinfty{E^*}$ defined by $f$ is given by
\[
f^v = f \circ \tau^*.
\]
In addition, if $X$ is a section of $E$ then the linear function
$\hat{X}\in\cinfty{E^*}$ defined by $X$ is given by
\[
\hat{X}(a^*)=a^*(X(\tau^*(a^*))),\text{ for all $a^*\in E^*$.}
\]
The Poisson bracket $\{ \; , \; \}_{E^*}$ is characterized by the
following relations
\begin{equation}\label{e2.8'}
\{f^v, g^v \}_{E^*} = 0, \qquad \{\hat{X}, g^v \}_{E^*} =
(\rho(X)g)^v, \qquand \{\hat{X}, \hat{Y}\}_{E^*} =
\widehat{\br{X}{Y}},
\end{equation}
for $X, Y \in \Gamma(E)$ and $f, g \in C^{\infty}(M)$.

In local coordinates $(x^i,y_\alpha)$ on $E^*$ we have that
\[
\{x^i,x^j \}_{E^*} = 0 \qquad \{y_\alpha,x^j\}_{E^*} =
\rho^j_\alpha \qquand \{y_\alpha,y_\beta\}_{E^*} = y_\gamma
C^\gamma_{\alpha\beta},
\]
(for more details, see \cite{CDW,C}).

\subsubsection{ The prolongation of a Lie algebroid over a smooth
map}\label{seccion2.1.1}

In this section we will recall the definition of the Lie algebroid
structure on the prolongation of a Lie algebroid over a smooth
map. We will follow \cite{HM} (see Section $1$ in \cite{HM}).

Let $(E,\lcf\cdot,\cdot\rcf,\rho)$ be a Lie algebroid of rank $n$
over a manifold $M$ of dimension $m$ and $f:M'\to M$ be a smooth
map.

We consider the subset ${\cal L}^fE$ of $E\times TM'$ defined by
\[
{\cal L}^fE=\{(b,v')\in E\times TM'/\rho(b)=(Tf)(v')\}
\]
where $Tf:TM\to TM'$ is the tangent map to $f$.

Denote by $\tau^f:{\cal L}^fE\to M'$ the map given by
\[
\tau^f(b,v')=\tau_{M'}(v'),
\]
for $(b,v')\in {\cal L}^fE,$ $\tau_{M'}:TM'\to M'$ being the
canonical projection. If $x'$ is a point of $M'$, it follows that
\[
(\tau^f)^{-1}(x')=({\cal L}^fE)_{x'}=\{(b,v')\in E_{f(x')}\times
T_{x'}M'/\rho(b)=(T_{x'}f)(v')\}
\]
is a vector subspace of $E_{f(x')}\times T_{x'}M'$, where
$E_{f(x')}$ is the fiber of $E$ over the point $f(x')\in M$.
Moreover, if $m'$ is the dimension of $M'$, one may  prove that
\[
\dim({\cal L}^fE)_{x'}=n+m'-\dim (\rho(E_{f(x')}) +
(T_{x'}f)(T_{x'}M')).\] Thus, if we suppose that there exists
$c\in \N$ such that
\begin{equation}\label{(1)}
\dim (\rho(E_{f(x')}) + (T_{x'}f)(T_{x'}M'))=c,\;\;\mbox{ for all
}x'\in M',
\end{equation} then we conclude that ${\cal L}^fE$ is a vector
bundle over $M'$ with vector bundle projection $\tau^f:{\cal
L}^fE\to M'.$
\begin{remark}{\rm If $\rho$ and $T(f)$ are transversal, that is,
\begin{equation}\label{(2)}
\rho(E_{f(x')}) + (T_{x'}f)(T_{x'}M')=T_{f(x')}M,\mbox{ for all }
x'\in M',
\end{equation}
then it is clear that (\ref{(1)}) holds. Note that if $E$ is a
transitive Lie algebroid (that is, $\rho$ is an epimorphism of
vector bundles) or $f$ is a submersion, we deduce that (\ref{(2)})
holds.}
\end{remark}

Next, we will assume that condition (\ref{(1)}) holds and we will
describe the sections of the vector bundle $\tau^f:{\cal L}^fE\to
M'$.

Denote by $f^*E$ the pullback of $E$ over $f$, that is,
\[
f^*E=\{(x',b)\in M'\times E/f(x')=\tau(b)\}.
\]
$f^*E$ is a vector bundle over $M'$ with vector bundle projection
\[
pr_1:f^*E\to M',\;\;\; (x',b)\in f^*E\to x'\in M'.\] Furthermore,
if $\sigma$ is a section of $pr_1:f^*E\to M'$ then
\[
\sigma=h'_i(X_i\circ f),
\]
for suitable $h_i'\in C^\infty(M')$ and $X_i\in \Gamma(E)$.

On the other hand, if $X^{\wedge}$ is a section of the vector
bundle $\tau^f:{\cal L}^fE\to M'$, one may prove that there exists
a unique $\sigma\in \Gamma(f^*E)$ and a unique $X'\in {\frak
X}(M')$ such that
\begin{equation}\label{(3)}
(T_{x'}f)(X'(x'))=\rho(\sigma(x')),\;\;\;\mbox{ for all }x'\in M',
\end{equation}
and $X^{\wedge}(x')=(\sigma(x'),X'(x')).$ Thus,
\[
X^{\wedge}(x')=(h'_{i}(x')X_i(f(x')),X'(x'))
\]
for suitable $h_i'\in C^\infty(M'),X_i\in \Gamma(E)$ and, in
addition,
\[
(T_{x'}f)(X'(x'))=h'_i(x')\rho(X_i)(f(x')).
\]
Conversely, if $\sigma\in \Gamma(f^*E)$ and $X'\in {\frak X}(M')$
satisfy condition (\ref{(3)}) then the map $X^{\wedge}:M'\to {\cal
L}^fE$ given by
\[
X^{\wedge}(x')=(\sigma(x'),X'(x')),\;\;\mbox{ for all } x'\in M',
\]
is a section of the vector bundle $\tau^f:{\cal L}^fE\to M'$.

Now, we consider the homomorphism of $C^\infty(M')$-modules
$\rho^f:\Gamma({\cal L}^fE)\to {\frak X}(M')$ and the Lie bracket
$\lcf\cdot,\cdot\rcf^f: \Gamma({\cal L}^fE)\times \Gamma({\cal
L}^fE)\to \Gamma({\cal L}^fE)$  on the space $\Gamma({\cal L}^fE)$
defined as follows. If $X^{\wedge}\equiv (\sigma,X')\in
\Gamma(f^*E)\times {\frak X}(M')$ is a section of $\tau^f:{\cal
L}^fE\to M'$ then
\begin{equation}\label{2.10'}
\rho^f(X^{\wedge})=X'
\end{equation}

and if $(\displaystyle h_i'(X_i\circ f),X')$ and $(\displaystyle
s_j'(Y_j\circ f),Y')$ are two sections of $\tau^f:{\cal L}^fE\to
M'$, with $h_i',s_j'\in C^\infty(M'),$ $X_i,Y_j\in \Gamma(E)$ and
$X',Y'\in {\frak X}(M'),$ then
\begin{equation}\label{2.10''}
\begin{array}{rcl}
\kern-8pt \lcf ( h_i'(X_i\circ f),X'),( s_j'(Y_j\circ f),Y')\rcf^f
& \kern-20pt = \kern-20pt & ( h_i's_j'(\lcf X_i,Y_j\rcf\circ f) +
X'(s_j')(Y_j\circ f)\\&& - Y'(h_i')(X_i\circ
f),[X',Y']).\end{array}
\end{equation}
The pair $(\lcf\cdot,\cdot\rcf^f,\rho^f)$ defines a Lie algebroid
structure on the vector bundle $\tau^f:{\cal L}^fE\to M'$  (see
\cite{HM}).

$({\cal L}^fE,\lcf\cdot,\cdot\rcf^f,\rho^f)$ is {\it the
prolongation of the Lie algebroid $E$ over the map $f$} ({\it the
inverse-image Lie algebroid of $E$ over $f$} in the terminology of
\cite{HM}).

On the other hand, if $pr_1:{\cal L}^fE\to E$ is the canonical
projection on the first factor then the pair $(pr_1,f)$ is a
morphism between the Lie algebroids $({\cal
L}^fE,\lcf\cdot,\cdot\rcf^f,\rho^f)$ and $(E,\lcf\cdot,\cdot
\rcf,\rho)$ (for more details, see \cite{HM}).

\subsubsection{Action Lie algebroids}
In this section, we will recall the definition of the Lie
algebroid structure of an action Lie algebroid. We will follow
again \cite{HM}.

Let $(E,\lcf\cdot,\cdot\rcf,\rho)$ be a Lie algebroid over a
manifold $M$ and $f:M'\to M$ be a smooth map. Denote by $f^*E$ the
pull-back of $E$ over $f$. $f^*E$ is a vector bundle over $M'$
whose vector bundle projection is the restriction to $f^*E$ of the
first canonical projection $pr_1:M'\times E\to M'.$

However, $f^*E$ is not, in general, a Lie algebroid over $M'$.

Now, suppose that $\Psi:\Gamma(E)\to {\frak X}(M')$ is an action
of $E$ on $f$, that is, $\Psi$ is a $\R$-linear map which
satisfies the following conditions:
\begin{enumerate}
\item
$\Psi(hX)=(h\circ f)\Psi X,$
\item
$\Psi\lcf X,Y\rcf=[\Psi X,\Psi Y],$
\item $\Psi X(h\circ f)=\rho(X)(h)\circ f,$
\end{enumerate}
for $X,Y\in \Gamma(E)$ and $h\in C^\infty(M)$. The action $\Psi$
allows us to introduce a homomorphism of $C^\infty(M')$-modules
$\rho_\Psi:\Gamma(f^*E)\to {\frak X}(M')$ and a Lie bracket
$\lcf\cdot,\cdot\rcf_\Psi:\Gamma(f^*E)\times \Gamma(f^*E)\to
\Gamma(f^*E)$ on the space $\Gamma(f^*E)$ defined as follows. If
$\sigma=\displaystyle h_i'(X_i\circ f)$ and $\gamma=\displaystyle
s_j'(Y_j\circ f)$ are sections of $f^*E$, with $h_i',s_j'\in
C^\infty(M')$ and $X_i,Y_j\in \Gamma(E),$ then
\[
\rho_\Psi(\sigma)=h_i'\Psi(X_i),
\]
\[
\begin{array}{rcl}
\lcf \sigma,\gamma\rcf_\Psi&=&\displaystyle h_i's_j'(\lcf
X_i,Y_j\rcf\circ f)+ \displaystyle \Psi(X_i)(s_j')(Y_j\circ
f)\\&&-\displaystyle s_j'\Psi(Y_j)(h_i')(X_i\circ f).
\end{array}
\]

The pair $(\lcf\cdot, \cdot \rcf_\Psi,\rho_\Psi)$ defines a Lie
algebroid structure on $f^*E.$ The corresponding Lie algebroid is
denoted by $E\ltimes M'$ or $E\ltimes f$ and we call it {\it an
action Lie algebroid } (for more details, see \cite{HM}).

\begin{remark}\label{prolon-acc}
{\rm Let $(E,\lcf\cdot,\cdot\rcf,\rho)$ be a Lie algebroid over a
manifold $M$.
\begin{enumerate}
\item
{\it A Lie subalgebroid} is a morphism of Lie algebroids $j:F\to
E,i:N\to M$ such that the pair $(j,i)$ is a monomorphism of vector
bundles and $i$ is an injective inmersion (see \cite{HM}).

\item
Suppose that $f: M' \to M$ is a smooth map and that $\Psi: \Gamma
(E) \to {\frak X}(M')$ is an action of $E$ on $f$. The anchor map
$\rho_{\Psi}$ of $f^*E$ induces a morphism between the vector
bundles $f^*E$ and $TM'$ which we will also denote by
$\rho_{\Psi}$. Thus, if ${\cal L}^{f}E$ is the prolongation of $E$
over $f$, we may introduce the map
\[
(id_{E}, \rho_{\Psi}): f^*E \to {\cal L}^fE,
\]
given by
\[
(id_{E}, \rho_{\Psi})(x', a) = (a, \rho_{\Psi}(x', a)),
\]
for $(x', a) \in (f^*E)_{x'} \subseteq \{x'\} \times E_{f(x')}$,
with $x' \in M'$. Moreover, if $id_{M'}: M' \to M'$ is the
identity map then, it is easy to prove that the pair $((id_{E},
\rho_{\Psi}), id_{M'})$ is a Lie subalgebroid. In fact, the map
$(id_{E}, \rho_{\Psi}): f^*E \to {\cal L}^{f}E$ is a section of
the canonical projection
\[
(id_{E}, \tau_{M'}): {\cal L}^{f}E \to f^*E
\]
defined by
\[
(id_{E}, \tau_{M'})(a, X_{x'}) = (x', a)
\]
for $(a, X_{x'}) \in ({\cal L}^{f}E)_{x'} \subseteq E_{f(x')}
\times T_{x'}M'$, with $x' \in M'$.
\end{enumerate} }
\end{remark}

\subsubsection{Quotient Lie algebroids by the action of a Lie
group}\label{seccion2.1.3}

Let $\pi:Q\to M$ be a principal bundle with structural group $G$.
Denote by $\phi:G\times Q\to Q$ the free action of $G$ on $Q$.

Now, suppose that $\tilde{E}$ is a vector bundle over $Q$ of rank
$n$, with vector bundle projection $\tilde{\tau}:{\tilde{E}}\to Q$
and that $\tilde{\phi}:G\times \tilde{E}\to \tilde{E}$ is an
action of $G$ on $\tilde{E}$ such that:
\begin{enumerate}
\item[$i)$] For each $g\in G$, the pair $(\tilde{\phi}_g,\phi_g)$
induces an isomorphism of vector bundles. Thus, the following
diagram

\begin{picture}(375,85)(40,20)
\put(190,20){\makebox(0,0){$Q$}}
\put(240,25){$\phi_g$}\put(210,20){\vector(1,0){80}}
\put(310,20){\makebox(0,0){$Q$}} \put(180,50){$\tilde\tau$}
\put(190,70){\vector(0,-1){40}} \put(320,50){$\tilde\tau$}
\put(310,70){\vector(0,-1){40}}
\put(190,80){\makebox(0,0){$\tilde{E}$}}
\put(240,85){$\tilde\phi_g$}\put(210,80){\vector(1,0){80}}
\put(310,80){\makebox(0,0){$\tilde{E}$}} \end{picture}

\vspace{5pt}

 is commutative and for each $q\in Q,$ the map
\[
\tilde{\phi}_g:\tilde{E}_q\to \tilde{E}_{\phi_g(q)}
\]
is a linear isomorphism between the vector spaces $\tilde{E}_q$
and $\tilde{E}_{\phi_g(q)}.$

\item[$ii)$] $\tilde{E}$ is covered by the ranges of equivariant charts,
that is, around each $q_0\in Q$ there is a $\pi$-satured open set
$\tilde{U}=\pi^{-1}(U)$, where $U\subseteq M$ is an open subset
with $x_0=\pi(q_0)\in U$ and a vector bundle chart
$\tilde\varphi:\tilde{U}\times \R^n\to
\tilde{\tau}^{-1}(\tilde{U})$ for $\tilde{E}$ which is equivariant
in the sense that
\[
\tilde{\varphi}(\phi_g(q),p)=\tilde\phi_g({\tilde\varphi}(q,p))
\]
for all $g\in G$, $q\in \tilde{U}$ and $p\in \R^{n}.$
\end{enumerate}
 Under the conditions $i)$ and $ii)$, the orbit set
$E=\tilde{E}/G$ has a unique vector bundle structure over $M=Q/G$
of rank $n$ such that the pair $(\tilde{\pi}, \pi)$ is a morphism
of vector bundles and $\tilde{\pi}:\tilde{E}\to E=\tilde{E}/G$ is
a surjective submersion, where $\tilde{\pi}:\tilde{E}\to
E=\tilde{E}/G$ is the canonical projection. The vector bundle
projection $\tau=\tilde{\tau}|G:E\to M$ of $E$ is given by
\[
\tau[\tilde{u}]=[\tilde{\tau}(\tilde{u})],\;\;\; \mbox{for }
\tilde{u}\in \tilde{E}.\] Moreover, if $q\in Q$ and $\pi(q)=x$
then the map
\[
\tilde{\pi}_{|\tilde{E}_q}: \tilde{E}_q\to E_x,\;\;\;\;
\tilde{u}\to [\tilde{u}]
\]
is a linear isomorphism between the vector spaces $\tilde{E}_q$
and $E_x$.

We call $(E,\tau, M)$ {\it the quotient vector bundle of
$(\tilde{E},\tilde{\tau},Q)$ by the action of $G$} (see
\cite{Ma}).

On the other hand, a section $\tilde{X}:Q\to \tilde{E}$ of
$\tilde{\tau}:\tilde{E}\to Q$ is said to be {\it invariant} if the
map $\tilde{X}$ is equivariant, that is, the following diagram

\begin{picture}(375,85)(40,20)
\put(200,20){\makebox(0,0){$Q$}}
\put(250,25){$\tilde{X}$}\put(210,20){\vector(1,0){80}}
\put(300,20){\makebox(0,0){$\tilde{E}$}} \put(185,50){$\phi_g$}
\put(200,70){\vector(0,-1){40}} \put(310,50){$\tilde{\phi}_g$}
\put(300,70){\vector(0,-1){40}} \put(200,80){\makebox(0,0){$Q$}}
\put(250,85){$\tilde{X}$}\put(210,80){\vector(1,0){80}}
\put(300,80){\makebox(0,0){$\tilde{E}$}} \end{picture}

\vspace{5pt}

 is commutative, for all $g\in G.$

We will denote by $\Gamma(\tilde{E})^G$ the set of invariant
sections of the vector bundle $\tilde{\tau}:\tilde{E}\to Q.$
$\Gamma(\tilde{E})^G$ is a $C^\infty(M)$-module where
\[
f\tilde{X}=(f\circ \pi)\tilde{X},\;\;\;  \mbox{for } f\in
C^\infty(M) \mbox{ and } \tilde{X}\in \Gamma (\tilde{E})^G.\]
Furthermore, there exists an isomorphism between the
$C^\infty(M)$-modules $\Gamma(E)$ and $\Gamma(\tilde{E})^G$. In
fact, if $\tilde{X}\in \Gamma (\tilde{E})^G$ then the
corresponding section $X\in \Gamma(E)$ is given by
\[
X(x)=\tilde{\pi}(\tilde{X}(q)),\;\;\; \mbox{ for } x\in M,
\]
with $q\in Q$ and $\pi(q)=x$ (for more details, see \cite{Ma}).

\begin{examples}\label{ex2.1'}{\rm $(a)$ Suppose that $\tilde{E}=TQ$ and
that $\tilde{\phi}:G\times TQ\to TQ$ is the tangent lift $\phi^T$
of $\phi$ defined by
\[ \phi^T_g=T\phi_g,\;\;\; \mbox{ for all } g\in G.
\]
 Then, $\phi^T$ satisfies conditions $i)$ and $ii)$ and, thus, one may consider
the quotient vector bundle $(E=TQ/G, \tau_Q|G, M)$ of
$(TQ,\tau_Q,M)$ by the action of $G.$ The space $\Gamma(TQ/G)$ may
be identified with the set of the vector fields on $Q$ which are
$G$-invariant.

\medskip

$(b)$ Assume that $\tilde{E}=T^*Q$ and that $\tilde{\phi}:G\times
T^*Q\to T^*Q$ is the cotangent lift $\phi^{T^*}$ of $\phi$ defined
by
\[
\phi^{T^*}_g=T^*\phi_{g^{-1}}, \mbox{ for all } g\in G.
\]

Then, $\phi^{T^*}$ satisfies conditions $i)$ and $ii)$ and,
therefore, one may consider the quotient vector  bundle
$(T^*Q/G,\pi_{Q}|G,M)$ of $(T^*Q,\pi_Q,Q)$ by the action of $G.$
Moreover, if $(\tau_{Q}|G)^*:(TQ/G)^*\to M$  is the dual vector
bundle to $\tau_{Q}|G:TQ/G\to M$ it is easy to prove that the
vector bundles $(\tau_{Q}|G)^*:(TQ/G)^*\to M$ and $\pi_{Q}|G:
T^*Q/G\to M$ are isomorphic.

\medskip

$(c)$ Suppose that ${\frak g}$ is the Lie algebra of $G,$ that
$\tilde{E}$ is the trivial vector bundle $pr_1:Q\times {\frak
g}\to Q$ and that the action $\tilde{\phi}=(\phi, Ad)$ of $G$ on
$Q\times {\frak g}$ is given by
\begin{equation}\label{PhiAd}
(\phi,Ad)_g(q,\xi)=(\phi_g(q), Ad_g\xi), \;\;\; \mbox{for } g\in G
\mbox{ and } (q,\xi)\in Q\times {\frak g},
\end{equation}
where $Ad: G\times {\frak g}\to {\frak g}$ is the adjoint
representation of $G$ on ${\frak g}$. Note that the space
$\Gamma(Q\times {\frak g})$ may be identified with the set of
$\pi$-vertical vector fields on $Q$. In addition, $\tilde{\phi}$
satisfies conditions $i)$ and $ii)$ and the resultant quotient
vector bundle $pr_1|G:\tilde{\frak g}=(Q\times {\frak g})/G\to
M=Q/G$ is just {\it the adjoint bundle associated with the
principal bundle} $\pi:Q\to M.$ Furthermore, if for each $\xi\in
{\frak g},$ we denote by $\xi_Q$ the infinitesimal generator of
the action $\phi$ associated with $\xi$, then the map
\[
j:\tilde{\frak g}\to TQ/G,\;\;\;\; [(q,\xi)]\to [\xi_Q(q)]
\]
induces a monomorphism between the vector bundles $\tilde{\frak
g}$ and $TQ/G$. Thus, $\tilde{\frak g}$ may be considered as a
vector subbundle of $TQ/G$. In addition, the space
$\Gamma(\tilde{\frak g})$ may be identified with the set of vector
fields on $Q$ which are vertical and $G$-invariant (see
\cite{Ma}).
\begin{remark}\label{r2.1''}{\rm $a)$ The tangent map to $\pi$,
$T\pi:TQ\to TM$, induces an epimorphism $[T\pi]:TQ/G \to TM$,
between the vector bundles $TQ/G$ and $TM$ and, furthermore, $Im
j=\ker [T\pi].$ Therefore, we have an exact sequence of vector
bundles

\begin{picture}(375,35)(-150,0)
\put(0,20){\makebox(0,0){$\tilde{\frak g}$}}
\put(10,20){\vector(1,0){40}} \put(75,20){\makebox(0,0){$TQ/G$}}
\put(25,25){$j$}
\put(100,20){\vector(1,0){40}}\put(105,25){$[T\pi]$}
\put(160,20){\makebox(0,0){$TM$}}
\end{picture}

\vspace{-10pt}

which is just {\it the Atiyah sequence associated with the
principal bundle} $\pi:Q\to M$ (for more details, see
\cite{Ma}).}\end{remark}

$b)$ Recall that if $\pi: Q \to M$ is a principal bundle with
structural group $G$ then a principal connection $A$ on $Q$ is a
Lie algebra valued one form $A: TQ \to {\frak g}$ such that:

(i) For all $\xi \in {\frak g}$ and for all $q \in Q$,
$A(\xi_{Q}(q)) = \xi$, and

(ii) $A$ is equivariant with respect to the actions $\phi^{T}:G
\times TQ \to TQ$ and $Ad: G \times {\frak g} \to {\frak g}$.

 Any choice of a connection in the principal bundle
$\pi:Q\to M$ determines an isomorphism between the vector bundles
$TQ/G\to M$ and $TM\oplus {\tilde{\frak g}}\to M$. In fact, if
$A:TQ\to {\frak g}$ is a principal connection then the map
$I_A:TQ/G\to TM\oplus \tilde{\frak g}$ defined by
\begin{equation}\label{IA}
I_A[\tilde{X}_q]=(T_q\pi)(\tilde{X}_q) \oplus [(q,A(\tilde{X}_q))]
\end{equation}
for $\tilde{X}_q\in T_qQ,$ is a vector bundle isomorphism over the
identity $id:M\to M$ (see \cite{CMR,Ma}).

Next, using the principal connection $A$, we will obtain a local
basis of $\Gamma(TQ/G)\cong \Gamma(TM\oplus \tilde{\frak g})\cong
{\frak X}(M)\oplus \Gamma(\tilde{\frak g}).$ First of all, we
choose a local trivialization of the principal bundle $\pi:Q\to M$
to be $U\times G$, where $U$  is an open subset of $M$. Thus, we
consider  the trivial principal bundle $\pi:U\times G\to U$ with
structural group $G$ acting only on the second factor by left
multiplication. Let $e$ be the identity element of $G$ and assume
that there are local coordinates $(x^i)$ in $U$ and that
$\{\xi_a\}$ is a basis of ${\frak g}$. Denote by $\{\xi_a^L\}$ the
corresponding left-invariant vector fields on $G$, that is,
\[
\xi_a^L(g)=(T_eL_g)(\xi_a),\;\;\; \mbox{for } g\in G,
\]
where $L_g:G\to G$ is the left translation by $g$, and suppose
that
\[
A(\frac{\partial }{\partial x^i}_{|_{(x,e)}})=A_i^a(x)\xi_a,
\]
for $i\in  \{1,\dots ,m\}$ and $x\in U$. Then, the horizontal lift
of the vector field $\displaystyle\frac{\partial }{\partial x^i}$
on $U$ is the vector field $(\displaystyle\frac{\partial
}{\partial x^i})^h$ on $U\times G$ given by
\[
(\frac{\partial }{\partial x^i})^h=\frac{\partial }{\partial
x^i}-A_i^a\xi_a^L.
\]
Therefore, the vector fields on $U\times G$
\begin{equation}\label{basis}
\{e_i=\frac{\partial }{\partial x^i}-A_i^a\xi_a^L,e_b=\xi_b^L\}
\end{equation}
are $G$-invariant and they define a local basis $\{e'_{i},
e'_{b}\}$ of $\Gamma(TQ/G)\cong {\frak X}(M)\oplus
\Gamma(\tilde{\frak g}).$ We will denote by $(x^i, y^i,y^b)$ the
corresponding fibred coordinates on $TQ/G$. In the terminology of
\cite{CMR},
\[
y^i=\dot{x}^i,\;\;\; y^b=\bar{v}^b,\;\;\; \mbox{ for } i\in
\{1,\dots ,m\} \mbox{ and } b\in \{1,\dots ,n\}.
\]
}\end{examples}

 Now, we will return to the general case.

Assume that $\pi:Q\to M$ is a principal bundle with structural
group $G$, that $\tilde{E}$ is a vector bundle over $Q$ of rank
$n$ with vector bundle projection $\tilde{\tau}:\tilde{E}\to Q$
and that $\tilde{\phi}:G\times \tilde{E}\to \tilde{E}$ is an
action of $G$ on $\tilde{E}$ which satisfies conditions $i)$ and
$ii)$. Denote by $\phi:G\times Q\to Q$ the free action of $G$ on
$Q$.

We will also suppose that $(\lcf\cdot,\cdot\rcf^{\tilde{\;}},
\tilde{\rho})$ is a Lie algebroid structure on
$\tilde{\tau}:\tilde{E}\to Q$ such that the space
$\Gamma(\tilde{E})^G$ is a Lie subalgebra of the Lie algebra
$(\Gamma(\tilde{E}), \lcf\cdot,\cdot\rcf^{\tilde{\;}}).$ Thus, one
may define a Lie algebra structure
\[
\lcf\cdot,\cdot \rcf:\Gamma(E)\times \Gamma(E)\to \Gamma(E)
\]
on $\Gamma(E)$. The Lie bracket $\lcf\cdot,\cdot\rcf$ is the
restriction of $\lcf\cdot,\cdot \rcf^{\tilde{\;}}$ to
$\Gamma(\tilde{E})^G\cong \Gamma(E)$.

 On the other hand, the anchor map $\tilde{\rho}:\tilde{E}\to TQ$
is equivariant. In fact, if $\tilde{X}\in \Gamma(\tilde{E})^G,$
$f\in C^\infty(M)$ and $\tilde{Y}\in \Gamma(\tilde{E})^G$ then
\[
\lcf \tilde{X}, (f\circ \pi)\tilde{Y}\rcf^{\tilde{\;}}=(f\circ
\pi) \lcf\tilde{X},\tilde{Y}\rcf^{\tilde{\;}} +
\tilde{\rho}(\tilde{X})(f\circ \pi)\tilde{Y}
\]
is an invariant section. This implies that the function
$\tilde{\rho}(\tilde{X})(f\circ \pi)$ is projectable, that is,
there exists $\rho(\tilde{X})(f)\in C^\infty(M)$ such that
\[
\tilde{\rho}(\tilde{X})(f\circ \pi)=\rho(\tilde{X})(f)\circ \pi,
\;\;\; \forall f\in C^\infty(M). \] The map
$\rho(\tilde{X}):C^\infty(M)\to C^\infty(M)$ defines a vector
field $\rho(\tilde{X})$ on $M$ and $\tilde{\rho}(\tilde{X})$ is
$\pi$-projectable onto $\rho(\tilde{X})$.

This proves that $\tilde{\rho}:\tilde{E}\to TQ$ is equivariant
and, therefore, $\tilde{\rho}$ induces a bundle map
$\rho:E=\tilde{E}/G\to TM=T(Q/G)$ such that the following diagram
is commutative

\begin{picture}(375,90)(40,10)
\put(200,20){\makebox(0,0){$E$}}
\put(250,25){$\rho$}\put(210,20){\vector(1,0){80}}
\put(310,20){\makebox(0,0){$TM$}} \put(180,50){$\tilde\pi$}
\put(200,70){\vector(0,-1){40}} \put(320,50){$T\pi$}
\put(310,70){\vector(0,-1){40}}
\put(200,80){\makebox(0,0){$\tilde{E}$}}
\put(250,85){$\tilde\rho$}\put(210,80){\vector(1,0){80}}
\put(310,80){\makebox(0,0){$TQ$}} \end{picture}

\vspace{5pt}

 Moreover,  it follows that the pair $(\lcf\cdot,
\cdot\rcf, \rho)$ is a Lie algebroid structure on the quotient
vector bundle $\tau=\tilde{\tau}|G:E=\tilde{E}/G\to M=Q/G.$ In
addition, from the definition of $(\lcf\cdot,\cdot\rcf, \rho)$,
one deduces that the pair $(\tilde{\pi}, \pi)$ is a morphism
between the Lie algebroids
$(\tilde{E},\lcf\cdot,\cdot\rcf^{\tilde{\;}}, \tilde{\rho})$ and
$(E,\lcf\cdot,\cdot\rcf,\rho)$.

We call $(E,\lcf\cdot,\cdot\rcf,\rho)$ {\it the quotient Lie
algebroid of $(\tilde{E},
\lcf\cdot,\cdot\rcf^{\tilde{\;}},\tilde{\rho})$ by the action of
the Lie group} $G$ ( a more general definition of a quotient Lie
algebroid may be found in \cite{HM}).

\begin{examples}{\rm $(i)$Assume that $\tilde{E}=TQ$    and that
$\tilde{\phi}$ is the tangent action $\phi^T:G\times TQ\to TQ$.
Consider on the vector bundle $\tau_Q:TQ\to Q$ the standard Lie
algebroid structure $([\cdot,\cdot],id)$. Since the Lie bracket of
two $G$-invariant vector fields on $Q$ is also $G$-invariant, we
obtain a Lie algebroid structure $(\lcf\cdot,\cdot\rcf, \rho)$ on
the quotient vector bundle $\tau_Q|G:E=TQ/G\to M=Q/G.$ We call
$(E=TQ/G,\lcf\cdot,\cdot \rcf, \rho)$ the {\it Atiyah algebroid
associated with the principal bundle $\pi:Q\to M$} (see
\cite{Ma}).

\medskip

$(ii)$ Suppose that $\tilde{E}$ is the trivial vector bundle
$pr_1:Q\times {\frak g}\to Q$. The space $\Gamma(Q\times {\frak
g})$ is isomorphic to the space of $\pi$-vertical vector fields on
$Q$ and, thus, $\Gamma(Q\times {\frak g})$ is a Lie subalgebra of
$(\frak{X}(Q),[\cdot,\cdot])$. This implies that the vector bundle
$pr_1:Q\times {\frak g}\to {\frak g}$ admits a Lie algebroid
structure $(\lcf\cdot,\cdot\rcf^{\tilde{\;}},\tilde\rho).$ On the
other hand, denote by $(\phi,Ad)$ the action of $G$ on $Q\times
{\frak g}$ given by $(\ref{PhiAd})$. Then, the space
$\Gamma(Q\times{\frak g})^G$ is isomorphic to the space ${\frak
X}^v(Q)^G$ of $\pi$-vertical $G$-invariant vector fields on $Q$.
Since ${\frak X}^v(Q)^G$ is a Lie subalgebra of $({\frak
X}(Q),[\cdot,\cdot])$, one may define a Lie algebroid structure on
the adjoint bundle $pr_1|G:\tilde{\frak g}=(Q\times{\frak g})/G\to
M=Q/G$ with anchor map $\rho=0$, that is, the adjoint bundle is a
Lie algebra bundle (see \cite{Ma}).}\end{examples}

Now, let $A:TQ\to {\frak g}$ be a connection in the principal
bundle $\pi:Q\to M$ and $B:TQ\oplus TQ\to {\frak g}$ be the
curvature of $A$. Using the principal connection $A$ one may
identity the vector bundles $E=TQ/G\to M=Q/G$ and $TM\oplus
\tilde{\frak g}\to M$, via the isomorphism $I_A$ given by
(\ref{IA}). Under this identification, the Lie bracket
$\lcf\cdot,\cdot\rcf$ on $\Gamma(TQ/G)\cong \Gamma(TM\oplus
\tilde{\frak g})\cong {\frak X}(M)\oplus {\frak X}^v(Q)^G$ is
given by
\[
\lcf X\oplus \tilde{\xi},Y\oplus \tilde{\eta}\rcf=[X,Y] \oplus
([\tilde{\xi},\tilde{\eta}] +
[X^h,\tilde{\eta}]-[Y^h,\tilde{\xi}]-B(X^h,Y^h)),
\]
for $X,Y\in {\frak X}(M)$ and $\tilde{\xi}, \tilde{\eta}\in {\frak
X}^v(Q)^G,$ where $X^h\in {\frak X}(Q)$ (respectively, $Y^h\in
{\frak X}(Q)$) is the horizontal lift of $X$ (respectively, $Y$),
via the principal connection $A$ (see \cite{CMR}). The anchor map
$\rho:\Gamma (TQ/G)\cong {\frak X}(M)\oplus {\frak X}^v(Q)^G\to
{\frak X}(M)$ is given by
\[
\rho(X\oplus \tilde{\xi})=X.
\]
Next, using the connection $A$, we will obtain the (local)
structure functions of $(E,\lcf\cdot,\cdot\rcf,\rho)$ with respect
to a local trivialization of the vector bundle.

First of all, we choose a local trivialization $U\times G$ of the
principal bundle $\pi:Q\to M$, where $U$ is an open subset of $M$
such that there are local coordinates $(x^i)$ on $U$. We will also
suppose that $\{\xi_a\}$ is a basis of ${\frak g}$ and that
\begin{equation}\label{ABcom}
\begin{array}{rcl}
A(\displaystyle\frac{\partial }{\partial
x^i}_{|(x,e)})&=&A_i^a(x)\xi_a,\\ B(\displaystyle\frac{\partial
}{\partial x^i}_{|(x,e)}, \displaystyle\frac{\partial }{\partial
x^j}_{|(x,e)})&=&B_{ij}^a(x)\xi_a,
\end{array}
\end{equation}
for $i,j\in \{1,\dots , m\}$ and $x\in U$. If $c_{ab}^c$ are the
structure constants of ${\frak g}$ with respect to the basis
$\{\xi_a\}$ then

\begin{equation}\label{curvat}
B_{ij}^c = \frac{\partial A_i^c}{\partial x^j}-\frac{\partial
A_j^c}{\partial x^i}-c_{ab}^cA_i^aA_j^b.
\end{equation}

Moreover, if $\{e'_i,e'_b\}$ is the local basis of $\Gamma(TQ/G)$
considered in Remark \ref{r2.1''} (see (\ref{basis})) then, using
(\ref{curvat}), we deduce that
\[
\lcf e'_i,e'_j\rcf=-B_{ij}^ce'_c,\;\;\; \lcf
e'_i,e'_a\rcf=c_{ab}^cA_{i}^be'_c, \; \; \lcf
e'_a,e'_b\rcf=c_{ab}^ce'_c,\;\;\;\] \[ \rho(e'_i)=\frac{\partial
}{\partial x^i}, \;\;\; \rho(e'_a)=0,\] for $i,j\in \{1,\dots
,m\}$ and $a,b\in \{1,\dots ,n\}.$ Thus, the local structure
functions of the Atiyah algebroid $\tau_{Q}|G:E=TQ/G\to M=Q/G$
with respect to the local coordinates $(x^i)$ and to the local
basis $\{e'_i,e'_a\}$ of $\Gamma(TQ/G)$ are
\begin{equation}\label{Atstfu}
\begin{array}{c}
C_{ij}^k= C_{ia}^j=-C_{ai}^j=
C_{ab}^i=0,\;\;\;C_{ij}^a=-B_{ij}^a,\;\;\;
C_{ia}^c=-C_{ai}^c=c_{ab}^cA_i^b,\;\;\;
C_{ab}^c=c_{ab}^c,\\
\rho_i^j=\delta_{ij},\;\;\; \rho_i^a=\rho_a^i=\rho_a^b=0.
\end{array}
\end{equation}
On the other hand, as we know the dual vector bundle to the Atiyah
algebroid is the quotient vector bundle $\pi_{Q}|G:T^*Q/G\to
M=Q/G$ of the cotangent bundle $\pi_Q:T^*Q\to Q$ by the cotangent
action $\phi^{T^*}$ of $G$ on $T^*Q$ (see Example \ref{ex2.1'}).
Now,  let $\Omega_{TQ}$ be the canonical symplectic $2$-form of
$T^*Q$ and $\Lambda_{TQ}$ be the Poisson $2$-vector on $T^*Q$
associated with $\Omega_Q$. If $(q^\alpha)$ are local coordinates
on $Q$ and $(q^\alpha,p_\alpha)$ are the corresponding fibred
coordinates on $T^*Q$ then
\[
\Omega_{TQ}=dq^\alpha\wedge dp_\alpha, \;\;\;
\Lambda_{TQ}=\frac{\partial }{\partial p_\alpha}\wedge
\frac{\partial }{\partial q^\alpha}.
\]

Note that $\Lambda_{TQ}$ is the linear Poisson structure on $T^*Q$
associated with the standard Lie algebroid $\tau_Q:TQ\to Q$. In
addition, it is well-known that the cotangent action $\phi^{T^*}$
is symplectic, that is,
\[
(\phi^{T^*})_g:(T^*Q,\Omega_{TQ})\to (T^*Q,\Omega_{TQ})
\]
is a symplectomorphism, for all $g\in G.$ Thus, the $2$-vector
$\Lambda_{TQ}$ on $TQ$ is $G$-invariant and it induces a
$2$-vector $\widetilde{\Lambda_{TQ}}$ on the quotient manifold
$T^*Q/G$. Under the identification between the vector bundles
$(\tau_{Q}|G)^*:E^*=(TQ/G)^*\to M=Q/G$ and $\pi_{Q}|G:T^*Q/G\to
M=Q/G,$ $\widetilde{\Lambda_{TQ}}$ is just the linear Poisson
structure $\Lambda_E=\Lambda_{TQ/G}$ on $E^*=(TQ/G)^*$ associated
with the Atiyah algebroid $\tau_{Q}|G:TQ/G\to M= Q/G$.  If
$(x^i,\dot{x}^i,\bar{v}^a)$ are the local coordinates on the
vector bundle $TQ/G$ considered in Remark \ref{r2.1''} and $(x^i,
p_i,\bar{p}_a)$ are the corresponding coordinates on the dual
vector bundle $(TQ/G)^*\cong T^*Q/G$ then, using (\ref{2.7'}) and
(\ref{Atstfu}), we obtain that the local expression of
$\Lambda_{TQ/G}$ is
\[
\begin{array}{rcl}
\Lambda_{TQ/G}&=&\displaystyle\frac{\partial }{\partial p_i}\wedge
\displaystyle\frac{\partial }{\partial x^i} +
c_{ab}^cA_i^b\bar{p}_c\displaystyle\frac{\partial }{\partial
p_i}\wedge \displaystyle\frac{\partial }{\partial
\bar{p}_a}\\[8pt]&& +
\displaystyle\frac{1}{2}(c_{ab}^c\bar{p}_c\displaystyle\frac{\partial
}{\partial \bar{p}_a}\wedge \displaystyle\frac{\partial }{\partial
\bar{p}_b} -B_{ij}^c\bar{p}_c\displaystyle\frac{\partial
}{\partial p_i}\wedge \displaystyle\frac{\partial }{\partial
p_j}).
\end{array}
\]

\subsection{Lagrangian mechanics on Lie algebroids}
In this section, we will recall some results about a geometric
description of Lagrangian Mechanics on Lie algebroids which has
been  developed by Mart\'{\i}nez in \cite{M}.

\subsubsection{The prolongation of a Lie algebroid over the vector
bundle projection}\label{seccion2.2.1} Let
$(E,\lcf\cdot,\cdot\rcf,\rho)$ be a Lie algebroid of rank $n$ over
a manifold $M$ of dimension $m$ and $\tau:E\to M$ be the vector
bundle projection.

If $f\in C^\infty(M)$ we will denote by $f^c$ and $f^v$ {\it the
complete and vertical lift} to $E$ of $f$. $f^c$ and $f^v$ are the
real functions on $E$ defined by
\begin{equation}\label{fcv}
f^c(a)=\rho(a)(f),\;\;\;\; f^v(a)=f(\tau(a)),
\end{equation}
for all $a\in E$.

Now, let $X$ be a section of $E$. Then, we can consider {\it the
vertical lift of $X$} as the vector field on $E$ given by
\[
X^v(a)=X(\tau(a))_a^v,\;\;\; \mbox{ for } a\in E,
\]
where $\;^v_{a}:E_{\tau(a)}\to T_a(E_{\tau(a)})$ is the canonical
isomorphism between the vector spaces $E_{\tau(a)}$ and
$T_a(E_{\tau(a)})$.

On the other hand, there exists a unique vector field $X^c$ on
$E$, {\it the complete lift of $X$,} satisfying the two following
conditions:
\begin{enumerate}
\item $X^c$ is $\tau$-projectable on $\rho(X)$ and \item
$X^c(\hat{\alpha})=\widehat{{\cal L}^E_X\alpha},$
\end{enumerate}
for all $\alpha\in \Gamma(E^*)$ (see \cite{GU,GU1}). Here, if
$\beta\in\Gamma(E^*)$ then $\hat\beta$ is the linear function on
$E$ defined by
\[
\hat\beta(b)=\beta(\tau(b))(b),\;\;\; \mbox{ for all  $b\in E$.}
\]
We have that (see \cite{GU,GU1})
\begin{equation}\label{cvstru}
[X^c,Y^c]=\lcf X,Y\rcf^c,\;\;\; [X^c,Y^v]=\lcf X,Y\rcf^v, \;\;\;
[X^v,Y^v]=0.\end{equation}

Next, we consider the prolongation ${\cal L}^\tau E$ of $E$ over
the projection $\tau$ (see Section \ref{seccion2.1.1}). ${\cal
L}^\tau E$ is a vector bundle over $E$ of rank $2n$. Moreover, we
may introduce {\it the vertical lift} $X^{\bf v}$ and {\it the
complete
 lift } $X^{\bf c}$ of a section $X\in \Gamma(E)$ as the sections of ${\cal L}^\tau E\to E$
 given by
 \begin{equation}\label{Defvc}
 X^{\bf v}(a)=(0,X^{v}(a)),\;\;\; X^{\bf c}(a)=(X(\tau(a)),X^{c}(a))
 \end{equation}
 for all $a\in E$. If $\{X_i\}$ is a local basis of
 $\Gamma(E)$ then $\{X_i^{\bf v},X_i^{\bf c}\}$ is a local
 basis of $\Gamma({\cal L}^\tau E)$ (see \cite{M}).

 Now, denote by $(\lcf\cdot,\cdot\rcf^\tau,\rho^\tau)$ the Lie
 algebroid structure on ${\cal L}^\tau E$ (see Section
 \ref{seccion2.1.1}). It follows that
 \begin{equation}\label{cvstru1}
 \begin{array}{ccc}
 \lcf X^{\bf c},Y^{\bf c}\rcf^\tau=\lcf X,Y\rcf^{\bf c},&\lcf X^{\bf c},Y^{\bf v}\rcf^\tau=\lcf
 X,Y\rcf^{\bf v},&\lcf X^{\bf v},Y^{\bf v}\rcf^\tau=0,\\
 \rho^\tau(X^{\bf c})(f^c)=(\rho(X)(f))^c,&\rho^\tau(X^{\bf c})(f^v)
 =(\rho(X)(f))^v,&\\
 \rho^\tau(X^{\bf v})(f^{c})=(\rho(X)(f))^{v},&\rho^\tau(X^{\bf v})(f^{v})=0,&
 \end{array}
 \end{equation}
 for $X,Y\in \Gamma(E)$ (see \cite{M}).

 Other two canonical objects on ${\cal L}^\tau E$ are {\it the
 Euler section } $\Delta$ and {\it the vertical endomorphism $S$}.
 $\Delta$ is the section of ${\cal L}^\tau E\to E$ defined by
 \[
 \Delta(a)=(0,a_a^v),\;\;\;\;\mbox{for all } a\in E\]
 and $S$ is the section of the vector bundle $({\cal L}^\tau
 E)\oplus ({\cal L}^\tau E)^*\to E$ characterized by the following
 conditions
 \[
 S(X^{\bf v})=0,\;\;\;\; S(X^{\bf c})=X^{\bf v},
 \]
for $X\in \Gamma(E).$

Finally, a section $\xi$ of ${\cal L}^\tau E\to E$ is said to be a
{\it second order differential equation} (SODE) on $E$ if
$S(\xi)=\Delta$ (for more details, see \cite{M}).

\begin{remark}\label{r2.2}
{\rm If $E$ is the standard Lie algebroid $TM$, then ${\cal
L}^\tau E=T(TM)$ and the Lie algebroid structure
$(\lcf\cdot,\cdot\rcf^\tau,\rho^\tau)$ is the usual one on the
vector bundle $T(TM)\to TM.$ Moreover, $\Delta$ is the Euler
vector field on $TM$ and $S$ is the vertical endomorphism on $TM.$
}\end{remark}
\begin{remark}\label{r2.3}{\rm Suppose that $(x^i)$ are
coordinates on an open subset $U$ of $M$ and that $\{e_\alpha\}$
is a basis of sections of $\tau^{-1}(U)\to U.$ Denote by
$(x^i,y^\alpha)$ the corresponding coordinates on $\tau^{-1}(U)$
and by $\rho_\alpha^i$ and $C_{\alpha\beta}^\gamma$ the
corresponding structure functions of $E$. If $X$ is a section of
$E$ and on $U$
\[
X=X^\alpha e_\alpha
\]
then $X^v$ and $X^c$ are the vector fields on $E$ given by
\begin{equation}\label{comverl}\begin{array}{rcl} X^v&=&\displaystyle
X^\alpha\displaystyle\frac{\partial}{\partial y^\alpha},\\
X^c&=&X^\alpha \rho_\alpha^i\displaystyle\frac{\partial }{\partial
x^i} + (\rho_\beta^i\displaystyle\frac{\partial X^\alpha}{\partial
x^i}-X^\gamma
C_{\gamma\beta}^\alpha)y^\beta\displaystyle\frac{\partial
}{\partial y^\alpha}.
\end{array}\end{equation}

In particular,
\begin{equation}\label{e.comve}
e_\alpha^v=\frac{\partial }{\partial y^\alpha},\;\;\;
e_\alpha^c=\rho_\alpha^i\frac{\partial }{\partial
x^i}-C_{\alpha\beta}^\gamma y^\beta\frac{\partial }{\partial
y^\gamma}.
\end{equation}

On the other hand, if $V_\alpha$, $T_\alpha$ are the sections of
${\cal L}^\tau E$ defined by $V_\alpha=e_\alpha^{\bf v}$,
$T_\alpha=e_\alpha^{\bf c}$ then
\begin{equation}\label{Xvc}
\begin{array}{rcl}
X^{\bf v}&=&X^\alpha V_\alpha,\\ X^{\bf c}&=&(\displaystyle
\rho^i_\beta\frac{\partial X^\alpha}{\partial x^i}y^\beta)V_\alpha
+ X^\alpha T_\alpha.
\end{array}
\end{equation}

Thus,
\begin{equation}\label{cvstrul}
\begin{array}{rcl}
\lcf T_\alpha, T_\beta\rcf^\tau &=&
(\rho^i_\delta\displaystyle\frac{\partial
C_{\alpha\beta}^\gamma}{\partial x^i}y^\delta)V_\gamma +
C_{\alpha\beta}^\gamma T_\gamma,\\ \lcf T_\alpha,V_\beta\rcf^\tau
&=& C_{\alpha\beta}^\gamma V_\gamma,\\ \lcf
V_\alpha,V_\beta\rcf^\tau&=&0.
\end{array}
\end{equation}
The local expression of $\Delta$ and $S$ are the following ones
\begin{equation}\label{2.20'}
\Delta=y^\alpha V_\alpha,\;\;\; S=T^\alpha\otimes V_\alpha,
\end{equation} where  $\{T^\alpha,V^\alpha\}$ is the dual basis of
$\{T_\alpha,V_{\alpha}\}$. Thus, a section $\xi$ of ${\cal L}^\tau
E$ is a SODE if and only if the local expression of $\xi$ is of
the form
\[
\xi=y^\alpha T_\alpha + \xi^\alpha V_\alpha,
\]
where $\xi^\alpha$ are arbitrary local functions on $E$.

Note that
\begin{equation}\label{difpro}
\begin{array}{rcl}
d^{{\cal L}^{\tau}E}f&=& (\rho_\alpha^i\displaystyle\frac{\partial
f}{\partial x^i}-C_{\alpha\beta}^\gamma
y^\beta\displaystyle\frac{\partial f}{\partial y^\gamma})T^\alpha
+ \frac{\partial f}{\partial y^\alpha}V^\alpha, \\ d^{{\cal
L}^{\tau}E}T^\gamma&=& \displaystyle -\frac{1}{2}
C_{\alpha\beta}^\gamma T^\alpha\wedge T^\beta, \\ d^{{\cal
L}^{\tau}E}V^\gamma&=&\displaystyle
-\frac{1}{2}(\rho_\delta^i\displaystyle\frac{\partial
C_{\alpha\beta}^\gamma}{\partial x^i}y^\delta)T^\alpha\wedge
T^\beta + C_{\alpha\beta}^\gamma T^\alpha \wedge V^\beta,
\end{array}
\end{equation}
for $f\in C^\infty(E)$ and $\gamma\in \{1,\dots ,n\}$.

We also remark that there exists another local basis of sections
on ${\cal L}^\tau E.$ In fact, we may define the local section
$\tilde{T}_\alpha$ as follows \begin{equation}\label{2.21'}
\tilde{T}_\alpha(a)=(e_\alpha(\tau(a)),\rho_\alpha^i\frac{\partial
}{\partial x^i}_{|a}),\;\;\;\mbox{ for all } a\in \tau^{-1}(U).
\end{equation}
Then, $\{\tilde{T}_\alpha,\tilde{V}_\alpha=V_\alpha\}$ is a local
basis of sections of ${\cal L}^\tau E$.

Note that
\begin{equation}\label{TtilT}
\tilde{T}_\alpha=T_\alpha + C_{\alpha\beta}^\gamma y^\beta
V_\gamma
\end{equation}
and thus,
\begin{equation}\label{tilde}
\begin{array}{lcl}
\lcf \tilde{T}_\alpha,\tilde{T}_\beta\rcf^\tau   =
C_{\alpha\beta}^\gamma\tilde{T}_\gamma,\;\;\;\; \lcf
\tilde{T}_\alpha,\tilde{V}_\beta\rcf^\tau  =  0,\;\;\;\; \lcf
\tilde{V}_\alpha,\tilde{V}_\beta\rcf^\tau  =  0, \\
\rho^{\tau}(\tilde{T}_{\alpha})  =  \displaystyle \rho_{\alpha}^i
\frac{\partial}{\partial x^i}, \;\;\;\;
\rho^{\tau}(\tilde{V}_{\alpha}) = \displaystyle
\frac{\partial}{\partial y^\alpha},
\end{array}
\end{equation}
for all $\alpha$ and $\beta$. Using the local basis
$\{\tilde{T}_\alpha,\tilde{V}_\alpha\}$ one may introduce, in a
natural way, local coordinates $(x^i,y^\alpha;z^\alpha,v^\alpha)$
on ${\cal L}^\tau E.$ If $\omega$ is a point of
$\tau^\tau(\tau^{-1}(U))$ $(\tau^\tau:{\cal L}^\tau E\to E$ being
the vector bundle projection) then $(x^i,y^\alpha)$ are the
coordinates of the point $\tau^\tau(\omega)\in \tau^{-1}(U)$ and
\[
\omega=z^\alpha\tilde{T}_\alpha(\tau^\tau(\omega)) + v^\alpha
\tilde{V}_\alpha(\tau^\tau(\omega)).
\]
In addition, the anchor map $\rho^\tau$ is given by
\begin{equation}\label{ancprol}
\rho^\tau(x^i,y^\alpha;z^\alpha,v^\alpha)=(x^i,y^\alpha;\rho_\alpha^iz^\alpha,v^\alpha)
\end{equation}
 and if $\{\tilde{T}^\gamma,\tilde{V}^\gamma\}$ is the dual basis
 of $\{\tilde{T}_\gamma,\tilde{V}_\gamma\}$ then
 \begin{equation}\label{2.24'}
 S=\tilde{T}^\alpha \otimes \tilde{V}_\alpha.
 \end{equation}
 and
\begin{equation}\label{Difftil}
\begin{array}{rcl}
d^{{\cal L}^{\tau}E}f&=& \rho_\gamma^i\displaystyle\frac{\partial
f}{\partial x^i}\tilde{T}^\gamma + \displaystyle\frac{\partial
f}{\partial y^\gamma}\tilde{V}^\gamma, \\ d^{{\cal
L}^{\tau}E}\tilde{T}^\gamma&=&-\displaystyle\frac{1}{2}
C_{\alpha\beta}^\gamma \tilde{T}^\alpha\wedge \tilde{T}^\beta,\;\;
d^{{\cal L}^{\tau}E}\tilde{V}^\gamma=0.
\end{array}
\end{equation}
} \end{remark}

\subsubsection{The Lagrangian formalism on Lie
algebroids}\label{seccion2.2.2}

Let $(E,\lcf\cdot,\cdot\rcf,\rho)$ be a Lie algebroid of rank $n$
over a manifold $M$ of dimension $m$ and  $L:E\to \R$ be a
Lagrangian function.

In this section, we will develop a geometric framework, which
allows to write the Euler-Lagrange equations associated with the
Lagrangian function $L$ in an intrinsic way (see \cite{M}).

First of all, we introduce {\it the Poincar\'{e}-Cartan $1$-section }
$\theta_L\in \Gamma(({\cal L}^\tau E)^*)$ associated with $L$
defined by
\begin{equation}\label{cartan1}
\theta_L(a)(\hat{X}_a)=(d^{{\cal L}^{\tau}E}
L(a))(S_a(\hat{X}_a))=\rho^\tau(S_a(\hat{X}_a))(L)
\end{equation}
for $a\in E$ and $\hat{X}_a\in ({\cal L}^\tau E)_a, ({\cal L}^\tau
E)_a$ being the fiber of ${\cal L}^\tau E\to E$ over the point
$a$. Then, {\it the Poincar\'{e}-Cartan $2$-section $\omega_L$ }
associated with $L$ is, up to the sign, the differential of
$\theta_L$, that is,
\begin{equation}\label{cartan2}
\omega_L=-d^{{\cal L}^{\tau}E}\theta_L
\end{equation}
and {\it the  energy function} $E_L$ is
\begin{equation}\label{EnerL}
E_L=\rho^\tau(\Delta)(L)-L.
\end{equation}

Now, let $\gamma:I=(-\varepsilon,\varepsilon)\subseteq \R\to E$ be
a curve in $E$. Then, $\gamma$ is {\it a solution of the
Euler-Lagrange equations } associated with $L$ if and only if:

\begin{enumerate}
\item $\gamma$ is admissible, that is,
$(\gamma(t),\dot{\gamma}(t))\in ({\cal L}^\tau E)_{\gamma(t)},$
for all $t$.
\item
$i_{(\gamma(t),\dot{\gamma}(t))}\omega_L(\gamma(t))=(d^{{\cal
L}^{\tau} E}E_L)(\gamma(t)),$ for all $t$.
\end{enumerate}

If $(x^i)$ are coordinates on $M$, $\{e_\alpha\}$ is a local basis
of $\Gamma(E)$, $(x^i,y^\alpha)$ are the corresponding coordinates
on $E$ and
\[ \gamma(t)=(x^i(t),y^\alpha(t)),
\]
then,  $\gamma$ is a solution of  the Euler-Lagrange equations if
and only if
\begin{equation}\label{2.38'}
\frac{d x^i}{dt}=\rho_\alpha^i y^\alpha,\;\;\;
\frac{d}{dt}(\frac{\partial L}{\partial
y^\alpha})=\rho_\alpha^i\frac{\partial L}{\partial
x^i}-C_{\alpha\beta}^\gamma y^\beta\frac{\partial L}{\partial
y^\gamma},
\end{equation}
for $i\in \{1,\dots ,m\}$ and $\alpha\in \{1,\dots ,n\},$ where
$\rho_\alpha^i$ and $C_{\alpha\beta}^\gamma$ are the structure
functions of the Lie algebroid $E$ with respect the coordinates
$(x^i)$ and the local basis $\{e_\alpha\}$ (see Section
\ref{seccion2.1.1}).

In particular, if $\xi\in \Gamma({\cal L}^\tau E)$ is a SODE and
\begin{equation}\label{*} i_\xi\omega_L=d^{{\cal L}^{\tau}E}E_L
\end{equation}
then the integral curves of $\xi$, that is, the integral curves of
the vector field $\rho^\tau(\xi)$ are solutions of the
Euler-Lagrange equations associated with $L$.

If the Lagrangian $L$ is regular, that is, $\omega_L$ is a
nondegenerate section then there exists a unique solution $\xi_L$
of the equation (\ref{*}) and $\xi_L$ is a SODE. In such a case,
$\xi_L$ is called {\it the Euler-Lagrange section associated with
$L$ } (for more details, see \cite{M}).

If $E$ is the standard Lie algebroid $TM$ then $\theta_L$
(respectively, $\omega_L$ and $E_L$) is the usual Poincar\'{e}-Cartan
$1$-form (respectively, the usual Poincar\'{e}-Cartan $2$-form and the
Lagrangian energy) associated with the Lagrangian function
$L:TM\to \R$. In this case, if $L:TM\to \R$ is regular, $\xi_L$ is
the Euler-Lagrange  vector field.

\begin{remark}
{\rm Suppose that $(x^i)$ are coordinates on $M$ and that
$\{e_\alpha\}$ is a local basis of sections of $E$. Denote by
$(x^i,y^\alpha)$ the corresponding coordinates on $E$, by
$\rho_\alpha^i$ and $C_{\alpha\beta}^\gamma$ the corresponding
structure functions of $E$ and by
$\{\tilde{T}_\alpha,\tilde{V}_\alpha\}$ the local basis of
$\Gamma({\cal L}^\tau E)$ considered in Remark \ref{r2.3}.

If $\{\tilde{T}^\alpha,\tilde{V}^\alpha\}$ is the dual basis of
$\{\tilde{T}_\alpha,\tilde{V}_\alpha\}$ then
\[
\begin{array}{rcl}
\theta_L&=&\displaystyle\frac{\partial L}{\partial
y^\alpha}\tilde{T}^\alpha,\\[8pt]
\omega_L&=&\displaystyle\frac{\partial^2L}{\partial
y^\alpha\partial y^\beta}\tilde{T}^\alpha\wedge \tilde{V}^\beta +
(\displaystyle\frac{1}{2}\displaystyle\frac{\partial L}{\partial
y^\gamma}C_{\alpha\beta}^\gamma-\rho_\alpha^i\displaystyle\frac{\partial^2
L}{\partial x^i\partial y^\beta})\tilde{T}^\alpha\wedge
\tilde{T}^\beta,\\[8pt] E_L&=&\displaystyle\frac{\partial
L}{\partial y^\alpha}y^\alpha-L.
\end{array}
\]
Thus, the Lagrangian $L$ is regular if and only if the matrix
$(W_{\alpha\beta})=( \displaystyle \frac{\partial^2 L}{\partial
y^\alpha\partial y^\beta})$ is regular. Moreover, if $L$ is
regular then the Euler-Lagrange section $\xi_L$ associated with
$L$ is given by

\begin{equation}\label{locxiL}
\xi_L=y^\alpha \tilde{T}_\alpha +
W^{\alpha\beta}(\rho^i_\beta\frac{\partial L}{\partial
x^i}-\rho_\gamma^iy^\gamma\frac{\partial^2 L}{\partial x^i\partial
y^\beta} + y^\gamma C_{\gamma\beta}^\nu\frac{\partial L}{\partial
y^\nu})\tilde{V}_\alpha,
\end{equation}

where $(W^{\alpha\beta})$ is the inverse matrix of
$(W_{\alpha\beta})$}.\end{remark}

\setcounter{equation}{0}
\section{Lie algebroids and Hamiltonian mechanics}
\subsection[The prolongation of a Lie algebroid over the dual vector bundle projection]
{The prolongation of a Lie algebroid over the vector
bundle projection of the dual bundle}\label{seccion3.1} Let
$(E,\lcf\cdot,\cdot \rcf,\rho)$ be a Lie algebroid of rank $n$
over a manifold $M$ of dimension $m$ and $\tau^*:E^*\to M$ be the
vector bundle projection of  the dual bundle $E^*$ to $E$.

We consider the prolongation ${\cal L}^{\tau^*}E$ of $E$ over
$\tau^*,$
\[
{\cal L}^{\tau^*}E=\{(b,v)\in E\times TE^*/
\rho(b)=(T\tau^*)(v)\}.
\]
${\cal L}^{\tau^*}E$ is a Lie algebroid over $E^*$ of rank $2n$
with Lie algebroid structure $(\lcf\cdot,\cdot\rcf^{\tau^*},
\rho^{\tau^*})$ defined as follows (see Section
\ref{seccion2.1.1}). If $(\displaystyle f_i'(X_i\circ \tau^*),X')$
and $(\displaystyle s_j'(Y_j\circ \tau^*),Y')$ are two sections of
${\cal L}^{\tau^*}E\to E^*$, with $f_i'$, $s_j'\in C^\infty(E^*),
X_i,Y_j\in \Gamma(E)$ and $X',Y'\in {\frak X}(E^*),$ then
\[
\begin{array}{l}
\lcf (\displaystyle f_i'(X_i\circ \tau^*),X'),(\displaystyle
s_j'(Y_j\circ \tau^*), Y')\rcf^{\tau^*}=
(\displaystyle f_i's_j'(\lcf X_i,Y_j\rcf\circ \tau^*) \\
+ \displaystyle X'(s_j')(Y_j\circ \tau^*)
-\displaystyle Y'(f_i')(X_i\circ \tau^*), [X',Y']),\\
\rho^{\tau^*}(\displaystyle f_i'(X_i\circ \tau^*),X')=X'.
\end{array}
\]
Now, if $(x^i)$ are local coordinates on an open subset $U$ of
$M$, $\{e_\alpha\}$ is a basis of sections of the vector bundle
$\tau^{-1}(U)\to U$ and $\{e^\alpha\}$ is the dual basis of
$\{e_\alpha\}$, then $\{\tilde{e}_\alpha,\bar{e}_\alpha\}$ is a
basis of sections of the vector bundle
$\tau^{\tau^*}((\tau^*)^{-1}(U))\to (\tau^*)^{-1}(U)$, where
$\tau^{\tau^*}:{\cal L}^{\tau^*}E\to E^*$ is the vector bundle
projection and
\begin{equation}\label{tilbar}
\begin{array}{rcl}
\tilde{e}_\alpha(a^*)&=&(e_\alpha(\tau^*(a^*)),\rho_\alpha^i\displaystyle\frac{\partial
}{\partial x^i}_{|a^*}),\\
\bar{e}_\alpha(a^*)&=&(0,\displaystyle\frac{\partial }{\partial
y_{\alpha}}_{|a^*})
\end{array}
\end{equation}
for $a^*\in (\tau^*)^{-1}(U).$ Here, $\rho_\alpha^i$ are the
components of the anchor map with respect to the basis
$\{e_\alpha\}$ and $(x^i,y_\alpha)$ are the local coordinates on
$E^*$ induced by the local coordinates $(x^i)$ and the basis
$\{e^\alpha\}$. In general, if $X=X^\gamma e_\gamma$ is a section
of the vector bundle $\tau^{-1}(U)\to U$ then one may consider the
sections $\tilde{X}$ and $\bar{X}$ of
$\tau^{\tau^*}((\tau^*)^{-1}(U))\to (\tau^*)^{-1}(U)$ defined by
\[
\begin{array}{rcl}
\tilde{X}(a^*)&=&(X(\tau^*(a^*)),(\rho_\gamma^iX^\gamma)\displaystyle\frac{\partial
}{\partial x^i}_{|a^*}),\\
\bar{X}(a^*)&=&(0,X^\gamma(\tau^*(a^*))\displaystyle\frac{\partial}{\partial
y_\gamma}_{|a^*}),
\end{array}
\]
for $a^*\in (\tau^*)^{-1}(U).$ Using the local basis
$\{\tilde{e}_\alpha,\bar{e}_\alpha\}$ one may introduce, in a
natural way, local coordinates $(x^i,y_\alpha;z^\alpha,v_\alpha)$
on ${\cal L}^{\tau^*}E.$ If $\omega^*$ is a point of
$\tau^{\tau^*}((\tau^*)^{-1}(U)$ then $(x^i,y_\alpha)$ are the
coordinates of the point $\tau^{\tau^*}(\omega^*)\in
(\tau^*)^{-1}(U)$ and
\[
\omega^*=z^\alpha\tilde{e_\alpha}(\tau^{\tau^*}(\omega^*)) +
v_\alpha \bar{e}_\alpha(\tau^{\tau^*}(\omega^*)).
\]
On the other hand, a direct computation proves that

\begin{equation}\label{dualstru}
\begin{array}{rcl}
\rho^{\tau^*}(\tilde{e}_\alpha)&=&\rho_\alpha^i\displaystyle\frac{\partial
}{\partial x^i},\;\;\;
\rho^{\tau^*}(\bar{e}_\alpha)=\displaystyle\frac{\partial
}{\partial y_\alpha},\\
\lcf\tilde{e}_\alpha,\tilde{e}_\beta\rcf^{\tau^*}&=&\widetilde{\lcf
e_\alpha, e_\beta\rcf}=C_{\alpha\beta}^\gamma\tilde{e}_\gamma,\\
\lcf\tilde{e}_\alpha,\bar{e}_\beta\rcf^{\tau^*}&=&\lcf
\bar{e}_\alpha, \bar{e}_\beta\rcf^{\tau^*}=0,
\end{array}
\end{equation}
for all $\alpha$ and $\beta$, $C_{\alpha\beta}^\gamma$ being the
structure functions of the Lie bracket $\lcf\cdot,\cdot\rcf$ with
respect to the basis $\{e_\alpha\}$. Thus, if
$\{\tilde{e}^\alpha,\bar{e}^\alpha\}$ is the dual basis of
$\{\tilde{e}_\alpha,\bar{e}_\alpha\}$ then
\begin{equation}\label{dif*}
\begin{array}{rcl}
d^{{\cal L}^{\tau^*}E}f&=&
\rho_\alpha^i\displaystyle\frac{\partial f}{\partial
x^i}\tilde{e}^\alpha + \displaystyle\frac{\partial f}{\partial
y_\alpha} \bar{e}^\alpha,\\ d^{{\cal
L}^{\tau^*}E}\tilde{e}^\gamma&=& \displaystyle -\frac{1}{2}
C_{\alpha\beta}^{\gamma} \tilde{e}^\alpha \wedge \tilde{e}^\beta,
\makebox[1cm]{} d^{{\cal L}^{\tau^*}E}\bar{e}^\gamma=0,
\end{array}
\end{equation}
for $f\in C^\infty(E^*).$

\begin{remark}{\rm  If $E$ is the tangent Lie algebroid $TM$, then the
Lie algebroid $({\cal
L}^{\tau^*}E,\lcf\cdot,\cdot\rcf^{\tau^*},\rho^{\tau^*})$ is the
standard Lie algebroid $(T(T^*M),[\cdot,\cdot],Id).$}
\end{remark}

\subsection{The canonical symplectic section of ${\cal
L}^{\tau^*}E$}\label{seccion3.2} Let $(E,\lcf\;,\;\rcf,\rho)$ be a
Lie algebroid of rank $n$ over a manifold $M$ of dimension $m$ and
${\cal L}^{\tau^*}E$ be the prolongation of $E$ over the vector
bundle projection $\tau^*:E^*\to M$. We may introduce a canonical
section $\lambda_E$ of the vector bundle $({\cal L}^{\tau^*}E)^*$
as follows. If $a^*\in E^*$ and $(b,v)$ is a point of the fiber of
${\cal L}^{\tau^*}E$ over $a^*$ then
\begin{equation}\label{Lio}
\lambda_E(a^*)(b,v)=a^*(b). \end{equation} $\lambda_E$  is called
{\it the Liouville section of $({\cal L}^{\tau^*}E)^*$.}

Now, {\it the canonical symplectic section } $\Omega_E$ is defined
by
\begin{equation}\label{sym}
\Omega_E=-d^{{\cal L}^{\tau^*}E} \lambda_E.\end{equation}

We have:
\begin{theorem}
$\Omega_E$ is a symplectic section of the Lie algebroid $({\cal
L}^{\tau^*}E, \lcf\cdot,\cdot\rcf^{\tau^*},\rho^{\tau^*})$, that
is,
\begin{enumerate}
\item $\Omega_E$ is a nondegenerate $2$-section and \item
$d^{{\cal L}^{\tau^*}E}\Omega_E=0.$
\end{enumerate}
\end{theorem}
\begin{proof}
It is clear that $d^{{\cal L}^{\tau^*}E}\Omega_E=0$.

On the other hand, if $(x^i)$ are local coordinates on an open
subset $U$ of $M$, $\{e_\alpha\}$ is a basis of sections of the
vector bundle $\tau^{-1}(U)\to U,$ $(x^i,y_\alpha)$ are the
corresponding local coordinates of $E^*$ on $(\tau^*)^{-1}(U)$ and
$\{\tilde{e}_{\alpha},\bar{e}_{\alpha}\}$ is the basis of the
vector bundle $\tau^{\tau^*}((\tau^*)^{-1}(U))\to
(\tau^*)^{-1}(U)$ then, using (\ref{Lio}), it follows that
\begin{equation}\label{locLio}
\lambda_E(x^i,y_\alpha)=y_\alpha\tilde{e}^\alpha
\end{equation}
where $\{\tilde{e}^\alpha,\bar{e}^\alpha\}$ is the dual basis to
$\{\tilde{e}_\alpha,\bar{e}_\alpha\}.$ Thus, from (\ref{dif*}),
(\ref{sym}) and (\ref{locLio}), we obtain that
\begin{equation}\label{locsym}
\Omega_E=\tilde{e}^{\alpha}\wedge \bar{e}^\alpha + \frac{1}{2}
C_{\alpha\beta}^\gamma y_\gamma \tilde{e}^\alpha\wedge
\tilde{e}^\beta
\end{equation}
$C_{\alpha\beta}^\gamma$ being the structure functions of the Lie
bracket $\lcf\cdot,\cdot\rcf$ with respect to the basis
$\{e_\alpha\}.$

Therefore, using (\ref{locsym}), we deduce that $\Omega_E$ is
nondegenerate $2$-section.
\end{proof}

\begin{remark}{\rm If $E$ is the standard Lie algebroid $TM$ then
$\lambda_E=\lambda_{TM}$ (respectively, $\Omega_E=\Omega_{TM}$) is
the usual Liouville $1$-form (respectively, the canonical
symplectic $2$-form) on $T^*M.$}
\end{remark}

Let $(E,\lcf\cdot,\cdot\rcf, \rho)$ be a Lie algebroid over a
manifold $M$ and $\gamma$ be a section of the dual bundle $E^*$ to
$E$.

Consider the morphism $((Id,T\gamma\circ \rho),\gamma)$  between
the vector bundles $E$ and ${\cal L}^{\tau^*}E$

\begin{picture}(375,90)(40,20)
\put(190,20){\makebox(0,0){$M$}}
\put(250,25){$\gamma$}\put(210,20){\vector(1,0){80}}
\put(310,20){\makebox(0,0){$E^*$}} \put(170,50){$$}
\put(190,70){\vector(0,-1){40}} \put(320,50){$$}
\put(310,70){\vector(0,-1){40}} \put(190,80){\makebox(0,0){$E$}}
\put(220,85){$(Id,T\gamma\circ
\rho)$}\put(210,80){\vector(1,0){80}}
\put(310,80){\makebox(0,0){${\cal L}^{\tau^*}E$}}
\end{picture}

defined by $(Id, T\gamma\circ \rho)(a)=(a,(T_x\gamma)(\rho(a))),$
for $a\in E_x$ and $x\in M$.
\begin{theorem}\label{t3.1''}If $\gamma$ is a section of $E^*$ then the pair
$((Id,T\gamma\circ \rho),\gamma)$ is a morphism between the Lie
algebroids $(E,\lcf\cdot,\cdot\rcf,\rho)$ and $({\cal
L}^{\tau^*}E, \lcf\cdot,\cdot\rcf^{\tau^*},\rho^{\tau^*}).$
Moreover,
\begin{equation}\label{3.35'} ((Id,T\gamma\circ
\rho),\gamma)^*\lambda_E=\gamma,\;\;\; ((Id,T\gamma\circ
\rho),\gamma)^*(\Omega_E)=-d^E \gamma.
\end{equation}
\end{theorem}
\begin{proof}
Suppose that $(x^i)$ are local coordinates on $M$, that
$\{e_\alpha\}$ is a local basis of $\Gamma(E)$ and that
\[
\gamma=\gamma_\alpha e^\alpha,
\]
with $\gamma_\alpha$ local real functions on $M$ and
$\{e^\alpha\}$ the dual basis to $\{e_\alpha\}$. Denote by
$\{\tilde{e}_\alpha,\bar{e}_\alpha\}$ the corresponding local
basis of $\Gamma({\cal L}^{\tau^*}E)$. Then, using (\ref{tilbar}),
it follows that
\begin{equation}\label{local}
((Id,T\gamma\circ \rho),\gamma)\circ e_\alpha=(\tilde{e}_\alpha +
\rho_\alpha^i\frac{\partial \gamma_\nu}{\partial
x^i}\bar{e}_\nu)\circ \gamma ,
\end{equation}
for $\alpha\in \{1,\dots ,n\}$, $\rho_\alpha^i$ being the
components of the anchor map with respect  to the local
coordinates $(x^i)$ and to the basis $\{e_\alpha\}.$ Thus, from
(\ref{diff0}),
\[
\begin{array}{rcl}
((Id,T\gamma\circ \rho),\gamma)^*(\tilde{e}_\alpha)&=&e^\alpha,\\
((Id,T\gamma\circ
\rho),\gamma)^*(\bar{e^\alpha})&=&\rho_\beta^i\displaystyle\frac{\partial
\gamma_\alpha}{\partial x^i}e^\beta = d^E \gamma_{\alpha},
\end{array}
\]
where $\{\tilde{e}^\alpha,\bar{e}^\alpha\}$ is the dual basis to
$\{\tilde{e}_\alpha,\bar{e}_\alpha\}.$

Therefore, from (\ref{diff0}), (\ref{diff1}) and (\ref{dif*}), we
obtain that the pair $((Id,T\gamma\circ \rho),\gamma)$ is a
morphism between the Lie algebroids $E\to M$ and ${\cal
L}^{\tau^*}E\to E^*.$

On the other hand, if $x$ is a point of $M$ and $a\in E_x$ then,
using (\ref{Lio}), we have that
\[
\{\{((Id,T\gamma\circ \rho),\gamma)^*(\lambda_E)\}(x)\}(a)=
\lambda_E(\gamma(x))(a,(T_x\gamma)(\rho(a)))=\gamma(x)(a)
\]
that is,
\[
((Id,T\gamma\circ \rho),\gamma)^*(\lambda_E)=\gamma.
\]
Consequently, from (\ref{sym}) and since the pair
$((Id,T\gamma\circ \rho),\gamma)$ is a morphism between the Lie
algebroids $E$ and ${\cal L}^{\tau^*}E$, we deduce that
\[
((Id,T\gamma\circ \rho),\gamma)^*(\Omega_E)=-d^E \gamma.
\]
\end{proof}

\begin{remark}{\rm Let $\gamma: M \longrightarrow T^*M$ be a
$1$-form on a manifold $M$ and $\lambda_{TM}$ (respectively,
$\Omega_{TM}$) be the Liouville $1$-form (respectively, the
canonical symplectic $2$-form) on $T^*M$. Then, using Theorem
\ref{t3.1''} (with $E = TM$), we deduce a well-known result (see,
for instance, \cite{AM}):
\[
\gamma^*(\lambda_{TM}) = \gamma, \makebox[.4cm]{}
\gamma^*(\Omega_{TM}) = -d^{TM} \gamma.
\]
}
\end{remark}

 From Theorem \ref{t3.1''}, we also obtain
\begin{corollary}\label{c3.4''}
If $\gamma \in \Gamma(E^*)$ is a $1$-cocycle  of $E$ then
\[
((Id,T\gamma\circ \rho),\gamma)^*(\Omega_E)=0.
\]
\end{corollary}

\subsection{The Hamilton equations}\label{seccion3.3}

Let $H:E^*\to \R$ be a Hamiltonian function. Then, since
$\Omega_E$ is a symplectic section of $({\cal
L}^{\tau^*}E,\lcf\cdot,\cdot\rcf^{\tau^*},\rho^{\tau^*})$ and
$dH\in \Gamma(({\cal L}^{\tau^*}E)^*),$ there exists a unique
section $\xi_H\in \Gamma({\cal L}^{\tau^*}E)$ satisfying
\begin{equation}\label{defxiH}
i_{\xi_H}\Omega_E=d^{{\cal L}^{\tau^*}E} H.
\end{equation}
With respect to the local basis
$\{\tilde{e}_\alpha,\bar{e}_\alpha\}$ of $\Gamma({\cal
L}^{\tau^*}E),$ $\xi_H$ is locally expressed as follows:
\begin{equation}\label{xiH0}
\xi_H=\frac{\partial H}{\partial
y_\alpha}\tilde{e}_{\alpha}-(C_{\alpha\beta}^\gamma
y_\gamma\frac{\partial H}{\partial y_\beta} +
\rho_\alpha^i\frac{\partial H}{\partial x^i})\bar{e}_\alpha.
\end{equation}

 Thus, the vector field $\rho^{\tau^*}(\xi_H)$ on $E^*$ is
given by
\[
\rho^{\tau^*}(\xi_H)=\rho_\alpha^i\frac{\partial H}{\partial
y_\alpha}\frac{\partial }{\partial x^i}-(C_{\alpha\beta}^\gamma
y_\gamma \frac{\partial H}{\partial y_\beta} +
\rho_\alpha^i\frac{\partial H}{\partial x^i})\frac{\partial
}{\partial y_\alpha}.
\]
Therefore, $\rho^{\tau^*}(\xi_H)$ is just {\it the Hamiltonian
vector field of $H$} with respect to the linear Poisson structure
$\Lambda_{E^*}$ on $E^*$ induced by Lie algebroid structure
$(\lcf\cdot,\cdot\rcf,\rho),$ that is,
\[
\rho^{\tau^*}(\xi_H)=i_{d^{TE^*}H}\Lambda_{E^*}.
\]
Consequently, the integral curves of $\xi_H$ (i.e., the integral
curves of the vector field $\rho^{\tau^*}(\xi_H)$) satisfy the
Hamilton equations for $H$,

\begin{equation}\label{Hameq}
\frac{d x^i}{dt}=\rho_\alpha^i\frac{\partial H}{\partial
y_\alpha},\;\;\; \frac{dy_\alpha}{dt}=-(C_{\alpha\beta}^\gamma
y_\gamma \frac{\partial H}{\partial y_\beta} +
\rho_\alpha^i\frac{\partial H}{\partial x^i}),
\end{equation}
for $i\in \{1,\dots ,m\}$ and $\alpha\in \{1,\dots ,n\}.$

\begin{remark}{\rm If $E$ is the standard Lie algebroid $TM$ then
$\xi_H$ is the Hamiltonian vector field of $H:T^*M\to \R$ with
respect to the  canonical symplectic structure of $T^*M$ and
equations (\ref{Hameq}) are the usual Hamilton equations
associated with $H.$}
\end{remark}

\subsection{Complete and vertical lifts}

On $E^*$ we have similar concepts of complete and vertical lifts
to those in $E$ (see \cite{Ma2}). The existence of a vertical lift
is but a consequence of $E^*$ being a vector bundle. Explicitly,
given a section $\alpha\in\sec{E^*}$ we define the vector field
$\alpha^v$ on $E^*$ by
$$
\alpha^v(a^*)=\alpha(\tau^*(a^*))^v_{a^*},\quad\text{for $a^*\in
E^*$,}
$$
where
$\map{{}^v_{a^*}}{E^*_{\tau^*(a^*)}}{T_{a^*}(E^*_{\tau^*(a^*)})}$
is the canonical isomorphism between the vector spaces
$E^*_{\tau^*(a^*)}$ and $T_{a^*}(E^*_{\tau^*(a^*)})$.

On the other hand, if $X$ is a section of $\tau: E \to M$, there
exists a unique vector field $X^{*c}$ on $E^*$, {\it the complete
lift of $X$} to $E^*$ satisfying the two following conditions:
\begin{enumerate}
\item $X^{*c}$ is $\tau^*$-projectable on $\rho(X)$ and
\item $X^{*c}(\hat{Y})=\widehat{[X,Y]},$
\end{enumerate}
for all $Y\in \Gamma(E)$ (see \cite{GU}). Here, if $X$ is a
section of $E$ then $\hat{X}$ is the linear function
$\hat{X}\in\cinfty{E^*}$ defined by
$$
\hat{X}(a^*)=a^*(X(\tau^*(a^*))),\text{ for all $a^*\in E^*$.}
$$
We have that (see \cite{GU})
\begin{equation}\label{cvstru-dual}
[X^{*c},Y^{*c}]=\lcf X,Y\rcf^{*c},\;\;\;
[X^{*c},\beta^v]=(\mathcal{L}^E_X\beta)^v, \;\;\;
[\alpha^v,\beta^v]=0.
\end{equation}

Next, we consider the prolongation $\prold$ of $E$ over the
projection $\tau^*$ (see Section \ref{seccion2.1.1}). $\prold$ is
a vector bundle over $E^*$ of rank $2n$. Moreover, we may
introduce {\it the vertical lift} $\alpha^{\bf v}$ and {\it the
complete
 lift } $X^{*\bf c}$ of a section $\alpha\in\sec{E^*}$ and a section  $X\in \sec{E}$ as the sections of ${\cal L}^{\tau^*}
 E \to E^*$ given by
 \begin{equation}\label{Defvc-dual}
 \alpha^{\bf v}(a^*)=(0,\alpha^{v}(a^*)),\;\;\;
 X^{* \bf c}(a^*)=(X(\tau^*(a^*)),X^{*c}(a^*))
 \end{equation}
 for all $a^*\in E^*$. If $\{X_i\}$ is a local basis of
 $\sec{E}$ and $\{\alpha_i\}$ is the dual basis of $\sec{E^*}$ then $\{\alpha_i^{\bf v},X_i^{* \bf c}\}$ is a local
 basis of $\sec{\prold}$.

 Now, denote by $(\lcf\cdot,\cdot\rcf^{\tau^*},\rho^{\tau^*})$ the Lie
 algebroid structure on $\prold$ (see Section
 \ref{seccion3.1}). It follows that
 \begin{equation}\label{cvstru1-dual-1}
 \lcf X^{* \bf c},Y^{* \bf c}\rcf^{\tau^*}=\lcf X,Y\rcf^{* \bf c},
 \qquad
 \lcf X^{* \bf c},\beta^{\bf v}\rcf^{\tau^*}=(\mathcal{L}^E_X\beta
)^{\bf v}, \qquad \lcf \alpha^{\bf v},\beta^{\bf
v}\rcf^{\tau^*}=0,
\end{equation}
and
\begin{equation}\label{cvstru1-dual-2}
\begin{array}{ll}
 \rho^{\tau^*}(X^{* \bf c})(f^v)
 =(\rho(X)(f))^v,\qquad
 & \rho^{\tau^*}(\alpha^{\bf v})(f^{v})=0,\\
 \rho^{\tau^*}(X^{* \bf c})(\hat{Y})
 =\widehat{\br{X}{Y}},
 & \rho^{\tau^*}(\alpha^{\bf v})(\hat{Y})=
  \alpha(Y)^v
 \end{array}
  \end{equation}
 for $X,Y\in\sec{E}$, $\alpha,\beta\in\sec{E^*}$ and $f \in
 C^{\infty}(M)$. Here, $f^v \in C^{\infty}(E^*)$ is the basic
 function on $E^*$ defined by
 \[
 f^v = f \circ \tau^*.
 \]


\bigskip

Suppose that $(x^i)$ are coordinates on an open subset $U$ of $M$
and that $\{e_\alpha\}$ is a basis of sections of $\tau^{-1}(U)\to
U$, and $\{e^\alpha\}$ is the dual basis of sections of
$\tau^*{}^{-1}(U)\to U$. Denote by $(x^i,y_\alpha)$ the
corresponding coordinates on $\tau^*{}^{-1}(U)$ and by
$\rho_\alpha^i$ and $C_{\alpha\beta}^\gamma$ the corresponding
structure functions of $E$. If $\theta$ is a section of $E^*$ and
on $U$
\[
\theta=\theta_\alpha e^\alpha
\]
and $X$ is a section of $E$ and on $U$
\[
X=X^\alpha e_\alpha,
\]
then $\theta^v$ and $X^{*c}$ are the vector fields on $E^*$ given
by
\begin{equation}\label{comverl-dual}\begin{array}{rcl} \theta^v&=&\displaystyle
\theta_\alpha\displaystyle\frac{\partial}{\partial y_\alpha}\\
X^{*c}&=&X^\alpha \rho_\alpha^i\displaystyle\pd{}{x^i}
-(\rho^i_\alpha \pd{X^\beta}{x^i}y_\beta +
C^\gamma_{\alpha\beta}y_\gamma X^\beta)\pd{}{y_\alpha}.
\end{array}\end{equation}

In particular,
\begin{equation}\label{e.comve-dual}
e_\alpha^v=\frac{\partial }{\partial y^\alpha},\;\;\;
e_\alpha^{*c}=\rho_\alpha^i\frac{\partial }{\partial x^i} -
C_{\alpha\beta}^\gamma y_\gamma\frac{\partial }{\partial y_\beta}.
\end{equation}

In terms of the basis $\{ \tilde{e}_\alpha,\bar{e}_\alpha \}$ of
sections of $\prold \to E^*$ we have the local expressions
$$
\theta^{\bf v}=\theta_\alpha\bar{e}_\alpha
$$
and
$$
X^{*{\bf c}}=X^\alpha\tilde{e}_\alpha-(\rho^i_\alpha
\pd{X^\beta}{x^i}y_\beta + C^\gamma_{\alpha\beta}y_\gamma
X^\beta)\bar{e}_{\alpha}.
$$

\begin{remark}
{\rm If $E$ is the standard Lie algebroid $TM$, then
$\prold=T(T^*M)$ and the vertical and complete lifts of sections
are the usual vertical and complete lifts.}
\end{remark}

The following properties relate vertical and complete lifts with
the Liouville 1-section and the canonical symplectic 2-section.
\begin{proposition}
If $X,Y$ are sections of $E$, $\beta,\theta$ are sections of $E^*$
and $\lambda_E$ is the Liouville 1-section, then
\begin{equation}\label{lifts-liouville}
\lambda_E(\theta^{\bf v})=0\qquand\lambda_E(X^{*\bf c})=\hat{X}.
\end{equation}
If $\Omega_E$ is the canonical symplectic 2-section then
\begin{equation}\label{lifts-symplectic}
\Omega_E(\beta^{\bf v},\theta^{\bf v})=0,\qquad \Omega_E(X^{* \bf
c},\theta^{\bf v})=\theta(X)^v \qquand \Omega_E(X^{* \bf c},Y^{*
\bf c})=- \widehat{\br{X}{Y}}.
\end{equation}
\end{proposition}
\begin{proof}
Indeed, for every $a^*\in E^*$,
$$
\lambda_E(\theta^{\bf v})(a^*)=a^*(pr_1(\theta^{\bf v}(a^*)))
=a^*(0_{\tau^*(a^*)})=0,
$$
and
$$
\lambda_E(X^{* \bf c})(a^*)=a^*(pr_1(X^{* \bf c}(a^*)))
=a^*(X(\tau^*(a^*)))=\hat{X}(a^*),
$$
which proves~(\ref{lifts-liouville}).

For~(\ref{lifts-symplectic}) we take into
account~(\ref{lifts-liouville}) and the definition of the
differential $\d$, so that
\begin{eqnarray*}
&&\Omega_E(\theta^{\bf v},\beta^{\bf v}) =\rho^{\tau^*}(\beta^{\bf
v})(\lambda_E(\theta^{\bf v})) -\rho^{\tau^*}(\theta^{\bf
v})(\lambda_E(\beta^{\bf v})) +\lambda_E(\lcf \theta^{\bf v},
\beta^{\bf v} \rcf^{\tau^*})
= 0, \\
&&\Omega_E(X^{* \bf c},\beta^{\bf v}) =\rho^{\tau^*}(\beta^{\bf
v})(\lambda_E(X^{* \bf c})) -\rho^{\tau^*}(X^{* \bf
c})(\lambda_E(\beta^{\bf v})) +\lambda_E(\lcf X^{* \bf c},
\beta^{\bf v} \rcf ^{\tau^*})\\
&& \makebox[2cm]{} =  \rho^{\tau^*}(\beta^{\bf v})\hat{X}=\beta(X)^v, \\
&&\Omega_E(X^{* \bf c},Y^{* \bf c}) =\rho^{\tau^*}(Y^{* \bf
c})(\lambda_E(X^{* \bf c})) -\rho^{\tau^*}(X^{* \bf
c})(\lambda_E(Y^{* \bf c}))
+\lambda_E(\lcf X^{* \bf c}, Y^{* \bf c} \rcf ^{\tau^*}) \\
&&\qquad\qquad\quad{}=\rho^{\tau^*}(Y^{* \bf
c})\hat{X}-\rho^{\tau^*}(X^{* \bf c})\hat{Y}+\widehat{\br{X}{Y}}=-
\widehat{\br{X}{Y}},
\end{eqnarray*}
where we have used~(\ref{cvstru1-dual-1})
and~(\ref{cvstru1-dual-2}).
\end{proof}

Noether's theorem has a direct generalization to the theory of
dynamical Hamiltonian systems on Lie algebroids. By an
infinitesimal symmetry of a section $X$ of a Lie algebroid we mean
another section $Y$ which commutes with $X$, that is,
$\br{X}{Y}=0$.
\begin{theorem}
Let $H\in\cinfty{E^*}$ be a Hamiltonian function and $\xi_H$ be
the corresponding Hamiltonian section. If $X\in\sec{E}$ is a
section of $E$ such that $\rho^{\tau^*}(X^{* \bf c})H=0$ then
$X^{* \bf c}$ is a symmetry of $\xi_H$ and the function $\hat{X}$
is a constant of the motion, that is
$\rho^{\tau^*}(\xi_H)\hat{X}=0$.
\end{theorem}
\begin{proof}
Using (\ref{defxiH}) and since $\rho^{\tau^*}(X^{* \bf{c}})(H) =
0$, it follows that
\begin{equation}\label{N1}
{\cal L}^{{\cal L}^{\tau^*}E}_{X^{* \bf{c}}}
(i_{\xi_{H}}\Omega_{E}) = d^{{\cal L}^{\tau^*}E}({\cal L}^{{\cal
L}^{\tau^*}E}_{X^{* \bf{c}}}H) = 0.
\end{equation}
On the other hand, from (\ref{cvstru1-dual-2}) and
(\ref{lifts-symplectic}), we obtain that
\begin{equation}\label{N2}
i_{X^{* \bf{c}}}\Omega_{E} = - d^{{\cal L}^{\tau^*}E} \hat{X}
\end{equation}
and thus
\[
{\cal L}^{{\cal L}^{\tau^*}E}_{X^{* \bf{c}}} \Omega_{E} = 0.
\]
Therefore, using (\ref{N1}), we deduce that
\[
i_{\lcf \xi_{H}, X^{* \bf{c}} \rcf^{\tau^*}}\Omega_{E} = 0,
\]
which implies that $\lcf \xi_{H}, X^{* \bf{c}}\rcf^{\tau^*} = 0$,
that is, $X^{* \bf{c}}$ is a symmetry of $\xi_{H}$.

In addition, from (\ref{defxiH}) and (\ref{N2}), we conclude that
\[
0 = \Omega_{E}(\xi_{H}, X^{* \bf{c}}) = {\cal L}^{{\cal
L}^{\tau^*}E}_{\xi_{H}}\hat{X} = \rho^{\tau^*}(\xi_{H})(\hat{X}).
\]
\end{proof}

\subsection{Poisson bracket}

Let $E$ be a Lie algebroid over $M$. Then, as we know, $E^*$
admits a linear Poisson structure $\Lambda_{E^*}$ with linear
Poisson bracket $\{ \; , \; \}_{E^*}$ (see Section 2.1).

Next, we will prove that the Poisson bracket $\{ \; , \; \}_{E^*}$
can also be defined in terms of the canonical symplectic 2-section
$\Omega_E$.
\begin{proposition}
Let $F,G\in\cinfty{E^*}$ be functions on $E^*$ and $\xi_F,\xi_G$
be the corresponding Hamiltonian sections. Then, the Poisson
bracket of $F$ and $G$ is
$$
\{F,G\}_{E^*} = -\Omega_E(\xi_F,\xi_G).
$$
\end{proposition}
\begin{proof}
We will see that the Hamiltonian section defined by a basic
function $f^v$ is $-(\d[E]f)^{\bf v}$, the vertical lift of
$-\d[E]f\in\sec{E^*}$, and that the Hamiltonian section defined by
a linear function $\hat{X}$ is $X^{* \bf{c}}$, the complete lift
of $X\in\sec{E}$.

Indeed, using (\ref{cvstru1-dual-2}) and (\ref{lifts-symplectic}),
we deduce that
\begin{eqnarray*}
&&\Omega_E(\xi_{f^v},\theta^{\bf v})=(\d[\prold]f^v)(\theta^{\bf
v})=0
=\Omega_E((-\d[E]f)^{\bf v},\theta^{\bf v}),\\
&&\Omega_E(\xi_{f^v},Y^{* \bf c})=(\d[\prold]f^v)(Y^{* \bf
c})=(\d[E]f(Y))^v =\Omega_E((-\d[E]f)^{\bf v},Y^{* \bf c}),
\end{eqnarray*}
and similarly
\begin{eqnarray*}
&&\Omega_E(\xi_{\hat{X}},\theta^{\bf
v})=(\d[\prold]\hat{X})(\theta^{\bf v})=\rho^{\tau^*}(\theta^{\bf
v})\hat{X}=(\theta(X))^v
=\Omega_E(X^{* \bf c},\theta^{\bf v}),\\
&&\Omega_E(\xi_{\hat{X}},Y^{* \bf c})=(\d[\prold]\hat{X})(Y^{* \bf
c})=\rho^{\tau^*}(Y^{* \bf c})\hat{X}=-
\widehat{\br{X}{Y}}=\Omega_E(X^{* \bf c},Y^{* \bf c}),
\end{eqnarray*}
for every $Y\in\sec{E}$ and every $\theta\in\sec{E^*}$. The proof
follows using (\ref{e2.8'}), (\ref{lifts-symplectic}) and by
taking into account that complete and vertical lifts generate
$\sec{\prold}$.
\end{proof}

\subsection[The Legendre transformation]{The Legendre transformation and the equivalence
between the Lagrangian and Hamiltonian
formalisms}\label{seccion3.4}

Let $L:E\to \R$ be a Lagrangian function and $\theta_L\in
\Gamma(({\cal L}^\tau E)^*)$ be the Poincar\'{e}-Cartan $1$-section
associated with $L$.

We introduce {\it the Legendre transformation associated with $L$}
as the smooth map $Leg_L:E\to E^*$ defined by
\begin{equation}\label{LegL}
Leg_L(a)(b)=\theta_L(a)(z),
\end{equation}
for $a,b\in E_x,$ where $E_x$ is the fiber of $E$ over the point
$x\in M$ and $z$ is a point in the fiber of  ${\cal L}^\tau E$
over the point $a$ such that
\[
pr_1(z)=b,
\]
$pr_1:{\cal L}^\tau E\to E$ being the restriction to ${\cal
L}^\tau E$ of the first canonical projection $pr_1:E\times TE\to
E$.

The map $Leg_L$ is well-defined and its local expression in fibred
coordinates on $E$ and  $E^*$ is
\begin{equation}\label{locLegL}
Leg_L(x^i,y^\alpha)=(x^i,\frac{\partial L}{\partial y^\alpha}).
\end{equation}
The Legendre transformation induces a map ${\cal L}{Leg_L}:{\cal
L}^\tau E\to {\cal L}^{\tau^*}E$ defined by
\begin{equation}\label{LLegL}
({\cal L}{Leg_L})(b,X_a)=(b,(T_aLeg_L)(X_a)),
\end{equation}
for $a,b\in E$ and $(b,X_a)\in ({\cal L}^\tau E)_a\subseteq
E_{\tau(a)}\times T_aE,$ where $TLeg_L:TE\to TE^*$ is the tangent
map of $Leg_L.$ Note that $\tau^*\circ Leg_L=\tau$ and thus ${\cal
L}{Leg_L}$ is well-defined.

Using (\ref{locLegL}), we deduce that the local expression of
${\cal L}{Leg_L}$ in the coordinates of ${\cal L}^\tau E$ and
${\cal L}^{\tau^*} E$ (see Sections \ref{seccion2.2.1} and
\ref{seccion3.1}) is
\begin{equation}\label{LLegloc}
{\cal L}{Leg_L}(x^i,y^\alpha;
z^\alpha,v^\alpha)=(x^i,\frac{\partial L}{\partial y^\alpha};
z^\alpha,\rho_\beta^iz^\beta\frac{\partial^2 L}{\partial
x^i\partial y^\alpha} + v^\beta\frac{\partial^2 L}{\partial
y^\alpha\partial y^\beta}). \end{equation}

\begin{theorem}\label{t3.2}
The pair $({\cal L}{Leg_L},Leg_L)$ is a morphism between the Lie
algebroids $({\cal L}^\tau E, \linebreak
\lcf\cdot,\cdot\rcf^{\tau},\rho^\tau)$ and $({\cal
L}^{\tau^*}E,\lcf\cdot,\cdot\rcf^{\tau^*},\rho^{\tau^*}).$
Moreover, if $\theta_L$ and $\omega_L$ (respectively, $\lambda_E$
and $\Omega_E)$ are the Poincar\'{e}-Cartan $1$-section and
$2$-section associated with $L$ (respectively, the Liouville
$1$-section and the canonical symplectic section on ${\cal
L}^{\tau^*}E$) then
\begin{equation}\label{pullback}
({\cal L}{Leg_L},Leg_L)^*(\lambda_E)=\theta_L,\;\;\; ({\cal
L}{Leg_L},Leg_L)^*(\Omega_E)=\omega_L.
\end{equation}
\end{theorem}
\begin{proof}
Suppose that $(x^i)$ are local coordinates on $M,$ that
$\{e_\alpha\}$ is a local basis of $\Gamma(E)$ and denote by
$\{\tilde{T}_\alpha,\tilde{V}_\alpha\}$ (respectively,
$\{\tilde{e}_\alpha,\bar{e}_\alpha\}$) the corresponding local
basis of $\Gamma({\cal L}^\tau E)$ (respectively, $\Gamma({\cal
L}^{\tau^*}E$)). Then, using (\ref{Difftil}) and (\ref{LLegloc}),
we deduce that
\[
({\cal
L}{Leg_L},Leg_L)^*(\tilde{e}^\gamma)=\tilde{T}^\gamma,\;\;\;
({\cal L}{Leg_L},Leg_L)^*(\bar{e}^\gamma)=d^{{\cal
L}^{\tau}E}(\frac{\partial L}{\partial y^\gamma}),\;\;\mbox{ for
all } \gamma.
\]
Thus, from (\ref{Difftil}) and (\ref{dif*}), we conclude that
\begin{equation}\label{e3.42'}
({\cal L}{Leg_L},Leg_L)^*(d^{{\cal L}^{\tau^*}E} f') = d^{{\cal
L}^\tau E}(f' \circ Leg_{L}),
\end{equation}
\begin{equation}\label{e3.42''}
 ({\cal L}{Leg_L},Leg_L)^*(d^{{\cal L}^{\tau^*}E}
\tilde{e}^\gamma)= d^{{\cal L}^\tau E}(({\cal
L}Leg_L,Leg_L)^*\tilde{e}^\gamma),
\end{equation}
\begin{equation}\label{e3.42'''}
({\cal L}{Leg_L},Leg_L)^*(d^{{\cal L}^{\tau^*}E} \bar{e}^\gamma)=
d^{{\cal L}^\tau E}(({\cal L}Leg_L,Leg_L)^*\bar{e}^\gamma),
\end{equation}
 for all $f' \in C^{\infty}(E^*)$
and for all $\gamma$. Consequently, the pair $({\cal
L}Leg_L,Leg_L)$ is a Lie algebroid morphism. This result also
follows using Proposition 1.8 in \cite{HM}.

Now, from (\ref{Lio}) and (\ref{LegL}), we obtain that
\[
({\cal L}{Leg_L},Leg_L)^*(\lambda_E)=\theta_L.
\]
Thus, using (\ref{cartan2}), (\ref{sym}) and the first part of the
theorem, we deduce that
\[
({\cal L}{Leg_L},Leg_L)^*(\Omega_E)=\omega_L.
\]
\end{proof}

 We also may prove the following result.

\begin{proposition}
The Lagrangian $L$ is regular if and only if the Legendre
transformation $Leg_L:E\to E^*$ is a local diffeomorphism.
\end{proposition}
\begin{proof}
$L$ is regular if and only if the matrix
$(\displaystyle\frac{\partial^2 L}{\partial y^\alpha\partial
y^\beta})$ is regular. Therefore, using (\ref{locLegL}), we deduce
the result.
\end{proof}

Next, we will assume that $L$ is {\it hyperregular}, that is,
$Leg_L$ is a global diffeomorphism. Then, from (\ref{LLegL}) and
Theorem \ref{t3.2}, we conclude that the pair $({\cal
L}Leg_L,Leg_L)$ is a Lie algebroid isomorphism. Moreover, we way
consider the Hamiltonian function $H:E^*\to \R$ defined by
\[
H=E_L\circ Leg_L^{-1},
\]
where $E_L: E\to \R$ is the Lagrangian energy associated with $L$
given by (\ref{EnerL}). The Hamiltonian section $\xi_H\in
\Gamma({\cal L}^{\tau^*}E)$ is characterized by the condition
\begin{equation}\label{xiH}
i_{\xi_H}\Omega_E=d^{{\cal L}^{\tau^*}E} H
\end{equation}
and we have following.

\begin{theorem}\label{t3.4}
If the Lagrangian $L$ is hyperregular then the Euler-Lagrange
section $\xi_L$ associated with $L$ and the Hamiltonian section
$\xi_H$ are $({\cal L}Leg_L,Leg_L)$-related, that is,
\begin{equation}\label{Related}
\xi_H\circ Leg_L={\cal L}Leg_L\circ \xi_L.
\end{equation}

Moreover, if $\gamma:I\to E$ is a solution of the Euler-Lagrange
equations associated with $L$, then $\mu=Leg_L\circ \gamma:I\to
E^*$ is a solution of the Hamilton equations associated with $H$
and, conversely, if $\mu:I\to E^*$ is a solution of the Hamilton
equations for $H$ then $\gamma=Leg_L^{-1}\circ \mu$ is a solution
of the Euler-Lagrange equations for $L$.
\end{theorem}
\begin{proof}
From (\ref{pullback}), (\ref{e3.42'}) and (\ref{xiH}), we obtain
that (\ref{Related}) holds. Now, using (\ref{Related}) and Theorem
\ref{t3.2}, we deduce the second part of the theorem.
\end{proof}

\begin{remark}{\rm If $E$ is the standard Lie algebroid $TM$ then
$Leg_L:TM\to T^*M$ is the usual Legendre transformation associated
with $L:TM\to \R$ and Theorem \ref{t3.4} gives the equivalence
between the Lagrangian and Hamiltonian formalisms in Classical
Mechanics.}
\end{remark}




\subsection{The Hamilton-Jacobi equation}
The aim of this section is to prove the following result.

\begin{theorem}\label{t3.10}
Let $(E,\br{\cdot}{\cdot},\an)$ be a Lie algebroid over a manifold $M$ and
$(\br[\prd]{\cdot}{\cdot},\an[\prd])$ be the Lie algebroid structure on
$\prold$. Let $\map{H}{E^*}{\R}$ be a Hamiltonian function and
$\xi_H\in\sec{\prold}$ be the corresponding Hamiltonian section.
Let $\alpha\in\sec{E^*}$ be a 1-cocycle, $\d[E]\alpha=0$, and
denote by $\sigma\in\sec{E}$ the section
$\sigma=pr_1\circ\xi_H\circ\alpha$. Then, the following two
conditions are equivalent:
\begin{enumerate}
\item For every curve $t\to c(t)$ in $M$ satisfying
\begin{equation}\label{proyxi}
\an(\sigma)(c(t))=\dot{c}(t),\text{ for all $t$,}
\end{equation}
the curve $t\to \alpha(c(t))$ on $E^*$ satisfies the Hamilton
equations for $H$.

\item $\alpha$ satisfies the Hamilton-Jacobi equation $\d[E](H\circ\alpha)=0$,
that is, the function $\map{H\circ\alpha}{M}{\R}$ is constant on
the leaves of the Lie algebroid foliation associated with $E$.
\end{enumerate}
\end{theorem}
\begin{proof}
For a curve $\map{c}{I = (-\epsilon, \epsilon) \subset \R}{M}$ on
the base we define the curves $\map{\mu}{I}{E^*}$ and
$\map{\gamma}{I}{E}$ by
$$
\mu(t)=\alpha(c(t))\qquand\gamma(t)=\sigma(c(t)).
$$
Since $\alpha$ and $\sigma$ are sections, it follows that both
curves project to $c$.

We consider the curve $v=(\gamma,\dot{\mu})$ in $E\times TE^*$ and
notice the following important facts about $v$:
\begin{itemize}
\item $v(t)$ is in $\prold$, for every $t\in I$, if and
only if $c$ satisfies (\ref{proyxi}). Indeed
$\an\circ\gamma=\an\circ\sigma\circ c$ while
$T\prd\circ\dot{\mu}=\dot{c}$.
\item In such a case, $\mu$ is a solution of the Hamilton equations if and only if $v(t)=\xi_H(\mu(t))$,
for every $t\in I$. Indeed, the first components coincide
$pr_1(v(t))=\gamma(t)$ and
$pr_1(\xi_H(\mu(t)))=pr_1(\xi_H(\alpha(c(t))))
=\sigma(c(t))=\gamma(t)$, and the equality of the second
components is just $\dot{\mu}(t)=\an[\prd](\xi_H(\mu(t)))$.
\end{itemize}
We also consider the map $\map{\Phi_{\alpha}}{E}{\prold}$ given by
$\Phi_{\alpha}=(Id,T\alpha\circ\rho)$, and we recall that
$\Omega_E(\alpha(x))(\Phi_{\alpha}(a),\Phi_{\alpha}(b))=0$, for
all $a,b\in E_{x}$, because of Corollary \ref{c3.4''}.

\smallskip\noindent[$(ii)$ $\Rightarrow$ $(i)$]
Assume that $c$ satisfies (\ref{proyxi}), so that $v(t)$ is a
curve in $\prold$. We have to prove that $v(t)$ equals to
$\xi_H(\mu(t))$, for every $t\in I$.

The difference $d(t)=v(t)-\xi_H(\mu(t))$ is vertical, that is,
$pr_{1}(d(t)) = 0$, for all $t$ (note that $pr_{1}(v(t)) =
pr_{1}(\xi_{H}(\mu(t))) = \gamma (t)$, for all $t$). Therefore, we
have that $\Omega_E (\mu(t))(d(t), \eta(t))= 0$, for every
vertical curve $t \to \eta(t)$ (see (\ref{locsym})).

Let $\map{a}{I}{E}$ be any curve on $E$ over $c$ (that is, $\tau
\circ a = c$) and consider its image under $\Phi_{\alpha}$, that
is $\zeta(t)=\Phi_{\alpha}(a(t))=(a(t),T\alpha(\an(a(t))))$. Since
$v(t)=\Phi_{\alpha}(\gamma(t))$ is also in the image of
$\Phi_{\alpha}$ we have that $\Omega_E(\mu(t))(v(t),\zeta(t))=0$.
Thus
\begin{eqnarray*}
\Omega_E (\mu(t))(d(t),\zeta(t)) &=&-\Omega_E
(\mu(t))(\xi_H(\mu(t)),\zeta(t))
=-\langle\d[\prold]H,\zeta(t)\rangle \\
& =&-T\alpha(\an(a(t)))H =-\langle \d[E](H\circ\alpha),a(t)\rangle
\end{eqnarray*}
which vanishes because $\d[E](H\circ\alpha)=0$.

Since any element in $(\prold)_{\alpha(x)}$, with $x \in M$, can
be obtained as a sum of an element in the image of $\Phi$ and a
vertical, we conclude that $\Omega_E
(\mu(t))(d(t),\eta(t))=\Omega_E (\alpha(c(t)))(d(t),\eta(t)) = 0$
for every curve $t \to \eta(t)$, which amounts $d(t)=0$.

\smallskip\noindent[$(i)$ $\Rightarrow$ $(ii)$]
Suppose that $x$ is a point of $M$ and that $b\in E_x$. We will
show that
\begin{equation}\label{HamJa}
\langle\d[E](H\circ\alpha),b\rangle=0.
\end{equation}
Let $c: I = (-\epsilon, \epsilon) \to M$ be the integral curve of
$\an(\sigma)$ such that $c(0)=x$. It follows that $c$ satisfies
(\ref{proyxi}). Let $\mu=\alpha\circ\ c$, $\gamma=\sigma\circ c$
and $v=(\gamma,\dot{\mu})$ as above. Since $\mu$ satisfies the
Hamilton equations, we have that $v(t)=\xi_H(\mu(t))$ for all $t$.
As above we take any curve $t \to a(t)$ in $E$ over $c$ such that
$a(0)=b$. Since $v(t)=\Phi_{\alpha}(\gamma(t))$ we have that
\begin{eqnarray*}
0&&=\Omega_E(\mu(t))(\Phi_{\alpha}(\gamma(t)),\Phi_{\alpha}(a(t)))
=\Omega_E(\mu(t))(v(t),\Phi_{\alpha}(a(t)))\\
&&=\Omega_E(\mu(t))(\xi_H(\mu(t)),\Phi_{\alpha}(a(t))) =\langle
\d[\prold]H,\Phi_{\alpha}(a(t))\rangle =\langle
\d[E](H\circ\alpha),a(t)\rangle.
\end{eqnarray*}
In particular, at $t=0$ we have that $\langle
\d[E](H\circ\alpha),b\rangle=0$.
\end{proof}

\begin{remark}
{\rm Obviously, we can consider as a cocycle $\alpha$ a
$1$-coboundary $\alpha=\d[E]S$, for some function $S$ on $M$.
Nevertheless, it should be noticed that on a Lie algebroid there
exist, in general, 1-cocycles that are not locally
$1$-coboundaries.}
\end{remark}

\begin{remark}\label{r3.10'}
{\rm If we apply Theorem \ref{t3.10} to the particular case when
$E$ is the standard Lie algebroid $TM$ then we directly deduce a
well-known result (see Theorem 5.2.4 in \cite{AM}).}
\end{remark}

Let $\map{L}{E}{\R}$ be an hyperregular Lagrangian and
$\map{\FL}{E}{E^*}$ be the Legendre transformation associated with
$L$. Denote by $\map{H}{E^*}{\R}$ the corresponding Hamiltonian
function, that is,
\[
H=E_L\circ\FL^{-1},
\]
$E_L$ being the  energy for $L$.

Now, suppose that $\alpha=\d[E]S$, for $\map{S}{M}{\R}$ a function
on $M$, is a solution of the Hamilton-Jacobi equation
$\d[E](H\circ \d[E]S)=0$ and that
$\map{\mu}{I=(-\varepsilon,\varepsilon)}{E^*}$ is a solution of
the Hamilton equations for $H$ such that $\mu(0)=d^E S(x),$ $x$
being a point of $M$. If $\map{c}{I}{M}$ is the projection of the
curve $\mu$ to $M$ (that is, $c=\prd\circ\mu$) then, from Theorem
\ref{t3.10}, we deduce that
\[
\mu = d^E S \circ c.
\]
On the other hand, the curve $\map{\gamma}{I}{E}$ given by
\[
\gamma=\FL^{-1}\circ\mu
\]
is a solution of the Euler-Lagrange equations associated with $L$
(see Theorem \ref{t3.4}). Moreover, since $c = \tau \circ \gamma$
and $\gamma$ is admissible, it follows that
\[
\dot{c}(t)=\rho(\gamma(t)),\text{ for all $t$}.
\]
Thus, if ${\cal F}^E$ is the Lie algebroid foliation associated
with $E$, we have that
\[
c(I)\subseteq {\cal F}_x^E,
\]
where ${\cal F}_x^E$ is the leaf of ${\cal F}^E$ over the point
$x$.

In addition, $H\circ\d[E]S$ is constant on ${\cal F}_x^E$. We will
show next that in the case that the constant is 0 the function $S$
is but the action.

\begin{proposition}\label{S.is.the.action}
Let $\alpha=\d[E]S$ be a solution of the Hamilton-Jacobi equation
such that
\[
H\circ d^E S=0\quad\text{on ${\cal F}_x^E$.}
\]
Then,
\[
\frac{d(S\circ c)}{dt}=L\circ \gamma.
\]
\end{proposition}
\begin{proof}
Since $\FL\circ\gamma=\mu$, we have that $(\Delta
L)\circ\gamma=\langle\mu,\gamma\rangle$, and from $E_L=\Delta L-L$
we get $L\circ\gamma=\langle\mu,\gamma\rangle-H\circ\mu$.
Moreover,
$$
\langle\mu(t),\gamma(t)\rangle=
\langle\d[E]S(c(t)),\gamma(t)\rangle
=\an(\gamma(t))S=\dot{c}(t)S=\frac{d}{dt}(S\circ c).
$$

In our case, $c$ is a curve on the leaf $\mathcal{F}_x^E$ and
$\mu=\d[E]S\circ c$, so that  $H\circ\mu=0$. Therefore, we
immediately get $L\circ\gamma= \displaystyle \frac{d}{dt}(S\circ
c)$.
\end{proof}

In particular, if we consider a curve $\map{c}{[0,T]}{M}$ such
that $c(T)=y$ then
$$
S(y)=\int_{0}^{T}L(\gamma(t))\,dt,
$$
where we have put $S(x)=0$, since $S$ is undetermined by a
constant.

\setcounter{equation}{0}
\section{The canonical involution for Lie algebroids}
\label{seccion4} Let $M$ be a smooth manifold, $TM$ be its tangent
bundle and $\tau_{M}: TM \to M$ be the canonical projection. Then,
the tangent bundle to $TM$, $T(T(M))$, admits two vector bundle
structures. The vector bundle projection of the first structure
(respectively, the second structure) is the canonical projection
$\tau_{TM}: T(T(M)) \to TM$ (respectively, the tangent map
$T(\tau_{M}): T(T(M)) \to TM$ to $\tau_{M}: TM \to M$). Moreover,
the canonical involution $\sigma_{TM}: T(T(M)) \to T(T(M))$ is an
isomorphism between the vector bundles $\tau_{TM}: T(T(M)) \to TM$
and $T(\tau_{M}): T(T(M)) \to TM$. We recall that if $(x^{i})$ are
local coordinates on $M$ and $(x^{i}, y^{i})$ (respectively,
$(x^{i}, y^{i}; \dot{x}^{i}, \dot{y}^{i})$) are the corresponding
fibred coordinates on $TM$ (respectively, $T(T(M))$) then the
local expression of $\sigma_{TM}$ is
\[
\sigma_{TM}(x^{i}, y^{i}; \dot{x}^{i}, \dot{y}^{i}) = (x^{i},
\dot{x}^{i}; y^{i}, \dot{y}^{i})
\]
(for more details, see \cite{Go,koba}).

Now, suppose that $(E,\lcf\cdot,\cdot,\rcf,\rho)$ is a Lie
algebroid of rank $n$ over a manifold $M$ of dimension $m$ and
that $\tau:E\to M$ is the vector bundle projection.

\begin{lemma}\label{vector.involucion}
Let $(b,v)$ be a point in $(\prol)_{a}$, with $a \in E_{x}$, that
is, $\tau^{\tau}(b,v)=a$. Then, there exists one and only one
tangent vector $\bar{v}\in T_bE$ such that:
\begin{enumerate}
\item $\bar{v}(f^v)=(\d[E]f)(x)(a)$, and
\item $\bar{v}(\hat{\theta})=v(\hat{\theta})+(\d[E]\theta)(x)(a,b)$
\end{enumerate}
for all $f\in\cinfty{M}$ and $\theta\in\sec{E^*}$.
\end{lemma}
\begin{proof}
Conditions $(i)$ and $(ii)$ determine a vector on $E$ provided
that they are compatible. To prove compatibility, we take $f \in
C^{\infty}(M)$ and $\theta \in \Gamma (E^*)$. Then,
$$
\bar{v}(\widehat{f\theta})=\bar{v}(f^v\hat{\theta})
=(\bar{v}(f^v))\hat{\theta}(b)+f^v(b)(\bar{v}(\hat{\theta})) =(d^E
f)(x)(a)\,\theta(x)(b)+f(x)\bar{v}(\hat{\theta}),
$$
where $x=\tau(a)=\tau(b)$, and on the other hand,
\begin{eqnarray*}
v(\widehat{f\theta})+(d^E(f\theta))(x)(a,b) &&=
(v(f^v))\hat{\theta}(a)+f^v(a)v(\hat{\theta})
+ (d^E f\wedge\theta)(x)(a,b) \\
&& +f(x)(d^E \theta)(x)(a,b)\\
&&= f(x)[v(\hat{\theta})+(d^E\theta)(x)(a,b)]+(d^E
f)(x)(a)\,\theta(x)(b),
\end{eqnarray*}
which are equal.
\end{proof}

The tangent vector $\bar{v}$ in the above
lemma~\ref{vector.involucion} projects to $\rho(a)$ since
$$
((T_b\tau)(\bar{v}))f=\bar{v}(f\circ\tau)=\bar{v}(f^v)=(d^E
f)(x)(a)=\rho(a)f,
$$
for all functions $f\in\cinfty{M}$, and thus
$(T_b\tau)(\bar{v})=\rho(a)$. It follows that $(a,\bar{v})$ is an
element of $(\prol)_{b}$, and we have defined a map from $\prol$
to $\prol$.

\begin{theorem}
\label{involucion} The map $\map{\sigma_E}{\prol}{\prol}$ given by
$$
\sigma_E(b,v)=(a,\bar{v}),
$$
where $\bar{v}$ is the tangent vector whose existence is ensured
by lemma~\ref{vector.involucion} is a smooth involution
interwining the projections $\tau^{\tau}$ and $pr_1$, that is
\begin{enumerate}
\item $\sigma_E^2=Id$, and
\item $pr_1 \circ \sigma_E= \tau^{\tau}$.
\end{enumerate}
\end{theorem}
\begin{proof}
We will find its coordinate expression. Suppose that $(x^i)$ are
local coordinates on an open subset $U$ of $M$, that
$\{e_{\alpha}\}$ is a basis of sections of the vector bundle
$\tau^{-1}(U) \to U$ and that $(x^i, y^{\alpha}; z^{\alpha},
v^{\alpha})$ are the coordinates of $(b, v)$, so that
$(x^i,y^\alpha)$ are the coordinates of $a$, $(x^i,z^\alpha)$ are
the coordinates of $b$ and $v=\displaystyle \rho^i_\alpha
z^\alpha\pd{}{x^i}+v^\alpha\pd{}{y^\alpha}$. Then
$\bar{v}=\displaystyle \rho^i_\alpha
y^\alpha\pd{}{x^i}+\bar{v}^\alpha\pd{}{y^\alpha}$ where we have to
determine $\bar{v}^\alpha$. Denote by $\{e^{\alpha}\}$ the dual
basis of $\{e_{\alpha}\}$. Applying $\bar{v}$ to
$y^\alpha=\widehat{e^\alpha}$ and taking into account the
definition of $\bar{v}$ we get
$$
\bar{v}^\alpha=\bar{v}(y^\alpha)=v(y^\alpha)+(\d[E]e^\alpha)(x)(a,b)
=v^\alpha-C^\alpha_{\beta\gamma}y^\beta z^\gamma.
$$
Therefore, in coordinates
\begin{equation}
\label{sigmaE.coordenadas} \sigma_E(x^i,y^\alpha;
z^\alpha,v^\alpha) =(x^i, z^\alpha;
y^\alpha,v^\alpha+C^\alpha_{\beta\gamma}z^\beta y^\gamma)
\end{equation}
which proves that $\sigma_E$ is smooth.

Moreover $\tau^{\tau}(b,v)=a$ and
$pr_1(\sigma_E(b,v))=pr_1(a,\bar{v})=a$. Thus, $(i)$ holds.
Finally, $\sigma_E$ is an involution. In fact, if $(b, v) \in
({\cal L}^{\tau}E)_{a}$ then
$$
\sigma_E^2(b,v)=\sigma_E(a,\bar{v})=(b,\bar{\bar{v}})
$$
and $\bar{\bar{v}}=v$ because both projects to $\rho(b)$ and over
linear functions
$$
\bar{\bar{v}}(\hat{\theta})=\bar{v}(\hat{\theta})+(d^E
\theta)(x)(b,a) =v(\hat{\theta})+(d^E\theta)(x)(a,b)+(d^E
\theta)(x)(b,a)=v(\hat{\theta}),
$$
which concludes the proof.
\end{proof}

\begin{definition}
The map $\sigma_E$ will be called {\it the canonical involution
associated with the Lie algebroid $E$.}
\end{definition}

If $E$ is the standard Lie algebroid $TM$ then $\sigma_E\equiv
\sigma_{TM}$ is the usual canonical involution
$\sigma_{TM}:T(TM)\to T(TM).$ In this case the canonical
involution has the following interpretation~(see~\cite{Tul1}). Let
$\map{\chi}{\R^2}{M}$ be a map locally defined in a neighborhood
of the origin in $\R^2$. We can consider $\chi$ as a 1-parameter
family of curves in two alternative ways. If $(s,t)$ are the
coordinates in $\R^2$, then we can consider the curve
$\chi_t:s\mapsto\chi(s,t)$, for fixed~$t$. If we take the tangent
vector at $s=0$ for every $t$ we get a curve $A(t)= \displaystyle
\frac{d\chi_t}{ds}\Big|_{s=0}=\pd{\chi}{s}(0,t)$ in $TM$ whose
tangent vector at $t=0$ is
$$
v=\dot{A}(0)=\frac{d}{dt}\frac{d\chi_t}{ds}\Big|_{s=0}\Big|_{t=0}
=\pd{}{t}\pd{\chi}{s}(0,0).
$$
$v$ is a tangent vector to $TM$ at $a=A(0)$. On the other hand, we
can consider the curve $\chi^s:t\mapsto\chi(s,t)$, for fixed~$s$.
If we take the tangent vector at $t=0$ for every $s$ we get a
curve $B(s)=\displaystyle
\frac{d\chi^s}{dt}\Big|_{t=0}=\pd{\chi}{t}(s,0)$ in $TM$ whose
tangent vector at $s=0$ is
$$
\bar{v}=\dot{B}(0)
=\frac{d}{ds}\frac{d\chi^s}{dt}\Big|_{t=0}\Big|_{s=0}
=\pd{}{s}\pd{\chi}{t}(0,0).
$$
$\bar{v}$ is a tangent vector to $TM$ at $b=B(0)$. We have that
$v\in T_a(TM)$ projects to $b$ and $\bar{v}\in T_b(TM)$ projects
to $a$. The canonical involution on $TM$ maps one of these vectors
into the other one, that is, $\sigma_{TM}(v)=\bar{v}$. Notice that
in terms of the tangent map $T\chi$ the curves $A$ and $B$ are
given by
$$
A(t)=T\chi\left(\pd{}{s}\Big|_{(0,t)}\right) \qquand
B(s)=T\chi\left(\pd{}{t}\Big|_{(s,0)}\right).
$$

\medskip

We look for a similar description in the case of an arbitrary Lie
algebroid. For that we consider a morphism $\map{\Phi}{T\R^2}{E}$,
locally defined in  $\tau_{\R^2}^{-1}(\mathcal{O})$ for some open
neighborhood $\mathcal{O}$ of the origin, from the standard Lie
algebroid $\map{\tau_{\R^2}}{T\R^2}{\R^2}$ to our Lie algebroid
$\map{\tau}{E}{M}$ and denote by $\chi$ the base map
$\map{\chi}{\R^2}{M}$, locally defined in $\mathcal{O}$. (The map
$\Phi$ plays the role of $T\chi$ in the standard case.)

The vector fields $\{\displaystyle \pd{}{s},\displaystyle
\pd{}{t}\}$ are a basis of $\sec{T\R^2}$ with dual basis
$\{ds,dt\}$. Thus we have the curves $\map{A}{\R}{E}$ and
$\map{B}{\R}{E}$ given by
$$
A(t)=\Phi\left(\pd{}{s}\Big|_{(0,t)}\right) \qquand
B(s)=\Phi\left(\pd{}{t}\Big|_{(s,0)}\right).
$$
Then $(B(0),\dot{A}(0))$ is an element of $(\prol)_{A(0)}$ and
$(A(0),\dot{B}(0))$ is an element of $(\prol)_{B(0)}$ and the
canonical involution maps one into the other. Indeed, let us write
$a=A(0)$, $b=B(0)$, $v=\dot{A}(0)$ and $\bar{v}=\dot{B}(0)$. Then,
applying the equality $T\chi=\rho\circ\Phi$ to $\displaystyle
\pd{}{t}\Big|_{(0,0)}$ we get that $\rho(b)=\displaystyle
\pd{\chi}{t}(0,0)$ and thus
$$
(T_a\tau)(v)=(T_a\tau)(\dot{A}(0))=\frac{d(\tau\circ
A)}{dt}\Big|_{t=0}=\pd{\chi}{t}(0,0)=\rho(b),
$$
which proves that $(b,v)\in\prol$ is an element of $(\prol)_{a}$.
Similarly, applying $T\chi=\rho\circ\Phi$ to $\displaystyle
\pd{}{s}\Big|_{(0,0)}$ we get that $\rho(a)=\displaystyle
\pd{\chi}{s}(0,0)$ and thus
$$
(T_b\tau)(\bar{v})=(T_b\tau)(\dot{B}(0))=\frac{d(\tau\circ
B)}{ds}\Big|_{s=0}=\pd{\chi}{s}(0,0)=\rho(a),
$$
which proves that $(a,\bar{v})\in\prol$ is an element of
$(\prol)_{b}$. Finally, to prove that $\sigma_E(b,v)=(a,\bar{v})$
we just need to prove that
$\bar{v}(\hat\theta)=v(\hat{\theta})+(\d[E]\theta)(x)(a,b)$, for
every $\theta\in\sec{E^*}$, with $x= \tau(a) = \tau(b)$. Applying
the condition $(\Phi,\chi)^*(\d[E]\theta)
=\d[T\R^2]((\Phi,\chi)^*\theta)$ to $\displaystyle
\pd{}{s}\Big|_{(0,0)}$ and $\displaystyle \pd{}{t}\Big|_{(0,0)}$,
we have
$$
(\d[E]\theta)(x)(a,b)=
\pd{}{s}\theta(B(s))\Big|_{s=0}-\pd{}{t}\theta(A(t))\Big|_{t=0}
=\dot{B}(0)(\hat{\theta})-\dot{A}(0)(\hat{\theta})
=\bar{v}(\hat{\theta})-v(\hat{\theta}),
$$
which proves $\sigma_E(b,v)=(a,\bar{v})$.

\bigskip

We will see that $\sigma_E$ is a morphism of Lie algebroids where
on $\map{pr_1}{\prol}{E}$ we have to consider the action Lie
algebroid structure that we are going to explain.

It is well-known (see, for instance,~\cite{Go}) that the tangent
bundle to $E$, $TE$, is a vector bundle over $TM$ with vector
bundle projection the tangent map to $\tau$, $T\tau:TE\to TM$.
Moreover, if $X$ is a section of $\tau: E \to M$ then the tangent
map to $X$
\[
TX:TM\to TE\] is a section of the vector bundle $T\tau:TE\to TM$.
We  may also consider the section $\hat{X}^0:TM\to TE$ of
$T\tau:TE\to TM$ given by
\begin{equation}\label{Xhat}
\hat{X}^0(u)=(T_x0)(u) +  X(x){\;^{v}_{0(x)}}
\end{equation}
for $u\in T_xM,$ where $0:M\to E$ is the zero section of $E$ and
$\;^{v}_{0(x)}:E_x\to T_{0(x)}(E_x)$ is the canonical isomorphism
between $E_x$ and $T_{0(x)}(E_x)$.

If $\{e_{\alpha}\}$ is a local basis of $\Gamma(E)$ then
$\{Te_{\alpha},\hat{e}_{\alpha}^0\}$ is a local basis of
$\Gamma(TE).$

Now, suppose that $(x^i)$ are local coordinates on an open subset
$U$ of $M$ and that $\{e_\alpha\}$ is a basis of the vector bundle
$\tau^{-1}(U)\to U.$ Denote by $(x^i,y^\alpha)$ (respectively,
$(x_i,\dot{x}_i)$ and $(x^i,y^\alpha;\dot{x}^i, \dot{y}^\alpha))$
the corresponding coordinates on the open subset $\tau^{-1}(U)$ of
$E$ (respectively, $\tau_M^{-1}(U)$ of $TM$ and
$\tau_E(\tau^{-1}(U))$ of $TE$, $\tau_M:TM\to M$ and $\tau_E:TE\to
E$ being the canonical projections). Then, the local expression of
$T\tau:TE\to TM$ is
\[
(T\tau)(x^i,y^\alpha;\dot{x}^i,\dot{y}^\alpha)=(x^i,\dot{x}^i)
\]
and the sum and product by real numbers in the vector bundle
$T\tau:TE\to TM$ are locally given by
\[
\begin{array}{l}
(x^i,y^\alpha;\dot{x}^i,\dot{y}^\alpha)\oplus (x^i,\bar{y}^\alpha;
\dot{x}^i,\dot{\bar{y}}^\alpha)=(x^i,y^\alpha+\bar{y}^\alpha;\dot{x}^i,\dot{y}^\alpha
+ \dot{\bar{y}}^\alpha)\\ \lambda\cdot
(x^i,y^\alpha;\dot{x}^i,\dot{y}^\alpha)=(x^i,\lambda
y^\alpha;\dot{x}^i,\lambda\dot{y}^\alpha)
\end{array}
\]

Moreover, if $X$ is a section of the vector bundle $\tau:E\to M$
and
\[
X=X^\alpha e_\alpha,
\]
then
\begin{equation}\label{TXhat}
\begin{array}{rcl}
TX(x^i,\dot{x}^i)&=&(x^i,X^\alpha;\dot{x}^i,\dot{x}^i\displaystyle\frac{\partial
X^\alpha}{\partial x^i}),\\ \hat{X}^0(x^i,\dot{x}^i)&=&(x^i,0;
\dot{x}^i,X^\alpha).
\end{array}
\end{equation}

Next, following \cite{MX} we define a Lie algebroid structure
$(\lcf\cdot,\cdot\rcf^T,\rho^T)$ on the vector bundle $T\tau:TE\to
TM$. The anchor map $\rho^T$ is given by $\rho^T=\sigma_{TM}\circ
T(\rho):TE\to T(TM)$, $\sigma_{TM}:T(TM)\to T(TM)$ being the
canonical involution and $T(\rho):TE\to T(TM)$ the tangent map to
$\rho:E\to TM$. The Lie bracket $\lcf\cdot,\cdot\rcf^T$ on the
space $\Gamma(TE)$ is characterized by the following equalities
\begin{equation}\label{Tstru}
\begin{array}{rcl}
\lcf TX,TY\rcf^T&=&T\lcf X,Y\rcf,\;\;\;\; \lcf TX,
\hat{Y}^0\rcf^T=\widehat{\lcf X,Y\rcf}^0 \\
\lcf\hat{X}^0,\hat{Y}^0\rcf^T&=&0,
\end{array}
\end{equation}

for $X,Y\in \Gamma(E).$

Now, we consider the pull-back vector bundle $\rho^*(TE)$ of
$T\tau:TE\to TM$ over the anchor map $\rho:E\to TM.$ Note that
\[
\rho^*(TE)=\{(b,v)\in E\times TE/\rho(b)=(T\tau)(v)\},
\]
that is, the total space of the vector bundle is just ${\cal
L}^\tau E$, the prolongation of $E$ over $\tau$. Thus,
$\rho^*(TE)$ is a vector bundle over $E$ with vector bundle
projection $pr_1:\rho^*(TE)\to E$ given by $pr_1(a,u)=a,$ for
$(a,u)\in \rho^*(TE).$

If $(x^i)$ are local coordinates on $M$, $\{e_\alpha\}$ is a local
basis of sections of the vector bundle $\tau:E\to M$ and
$(x^i,y^\alpha)$ (respectively,
$(x^i,y^\alpha;z^\alpha,v^\alpha)$) are the corresponding
coordinates on $E$ (respectively, ${\cal L}^\tau
E\equiv\rho^*(TE)$) then the local expressions of the vector
bundle projection, the sum and the product by real numbers in the
vector bundle $\rho^*(TE)\to E$ are
\[
\begin{array}{l}
pr_1(x^i,y^\alpha;z^\alpha,v^\alpha)=(x^i,z^\alpha)\\
(x^i,y^\alpha;z^\alpha,v^\alpha)\oplus
(x^i,\bar{y}^\alpha;z^\alpha,\bar{v}^\alpha)=(x^i,y^\alpha+\bar{y}^\alpha;z^\alpha,
v^\alpha + \bar{v}^\alpha)\\
\lambda(x^i,y^\alpha;z^\alpha,v^\alpha)=(x^i,\lambda
y^\alpha;z^\alpha,\lambda v^\alpha).
\end{array}
\]
Moreover, we have the following result.

\begin{theorem}
There exists a unique action $\Psi:\Gamma(TE)\to {\frak X}(E)$ of
the Lie algebroid $(TE,\lcf\;,\;\rcf^T,\rho^T)$ over the anchor
map $\rho:E\to TM$ such that
\begin{equation}\label{4.52'}
 \Psi(TX)=X^c,\;\;\;
\Psi(\hat{X}^0)=X^v,
\end{equation}
 for $X\in \Gamma(E)$, where
$X^c$ (respectively, $X^v$) is the complete lift (respectively,
the vertical lift) of $X$. In fact, if $\{e_\alpha\}$ is a local
basis of $\Gamma(E)$ and $\bar{X}$ is a section of $T\tau:TE\to
TM$ such that
\[
\bar{X}=\bar{f}^\alpha T(e_\alpha) + \bar{g}^\alpha
\hat{e}^0_\alpha
\]
with $\bar{f}^\alpha,\bar{g}^\alpha$ local real functions on $TM$,
then
\[
\Psi(\bar{X})=(\bar{f}^\alpha\circ \rho)e_\alpha^c +
(\bar{g}^\alpha\circ \rho) {e}^{v}_\alpha.
\]
\end{theorem}

\begin{proof}
A direct computation, using (\ref{estruc1}), (\ref{comverl}) and
(\ref{TXhat}) proves that
\begin{equation}\label{projec}
T\rho\circ X^c=\rho^T\circ TX\circ \rho,\;\;\; T\rho\circ
X^v=\rho^T\circ \hat{X}^0 \circ \rho,
\end{equation}
for $X\in \Gamma(E)$, that is, the vector field $X^c$
(respectively, $X^v$) is $\rho$-projectable to the vector field
$\rho^T(TX)$ (respectively, $\rho^T(\hat{X}^0)$) on $TM.$

Now, from (\ref{cvstru}), (\ref{Tstru}) and (\ref{projec}), we
deduce the result.
\end{proof}

\begin{corollary}\label{c4.2}
There exists a Lie algebroid structure
$(\lcf\cdot,\cdot\rcf_\Psi^T,\rho_\Psi^T)$ on the vector bundle
$\rho^*(TE) \linebreak \to E$ such that if $h^i(\bar{X}_i\circ
\rho)$ and $s^j(\bar{Y}_j\circ \rho)$ are two sections of
$\rho^*(TE)\to E$ with $h^i,s^j\in C^\infty(E)$ and
$\bar{X}_i,\bar{Y}_j$ sections of $T\tau:TE\to TM$ then
\[
\begin{array}{rcl}
\rho_\Psi^T( h^i(\bar{X}_i\circ
\rho))&=&h^i\Psi(\bar{X}_i),\\
\lcf h^i(\bar{X}_i\circ \rho), s^j(\bar{Y}_j\circ
\rho)\rcf_\Psi^T&=&h^is^j(\lcf\bar{X}_i,\bar{Y}_j\rcf^T\circ
\rho)\\ &&+ h^i\Psi(\bar{X}_i)(s^j)(\bar{Y}_j\circ \rho)- s^j\Psi
(\bar{Y}_j)(h^i)(\bar{X_i}\circ \rho).
\end{array}
\]
\end{corollary}
Next, we will give a local description of the Lie algebroid
structure $(\lcf\cdot,\cdot\rcf_\Psi^T,\rho_\Psi^T)$. For this
purpose, we consider local coordinates $(x^i)$ on $M$ and a local
basis $\{e_\alpha\}$ of $\Gamma(E)$. Denote by $\rho^i_\alpha$ and
$C_{\alpha\beta}^\gamma$ the structure functions of $E$ with
respect to $(x^i)$ and $\{e_\alpha\}$, by $(x^i,y^\alpha;
z^\alpha,v^\alpha)$ the corresponding coordinates on
$\rho^*(TE)\equiv {\cal L}^\tau E$ and by
$(x^i,y^\alpha;\dot{x}^i,\dot{y}^\alpha)$ the corresponding
coordinates on $TE.$

If $X$ is a section of $E$ and $X=X^\alpha e_\alpha$, we may
introduce the sections $T^\rho X$ and $\hat{X}^\rho$ of
$\rho^*(TE)$ given by
\[
T^\rho X=TX\circ \rho,\;\;\; \hat{X}^\rho=\hat{X}^0 \circ \rho.
\]
We have that (see (\ref{TXhat}))

\begin{equation}\label{TroXXro}
\begin{array}{rcl}
T^\rho X(x^i,z^\alpha)&=&(x^i,X^\alpha;
\rho_\alpha^iz^\alpha,\rho_\beta^iz^\beta\displaystyle\frac{\partial
X^\alpha}{\partial x^i}),\\
\hat{X}^\rho(x^i,z^\alpha)&=&(x^i,0;z^\alpha,X^\alpha).
\end{array}
\end{equation}
Thus, $\{T^\rho e_\alpha,\hat{e}_\alpha^\rho\}$ is a local basis
of $\Gamma(\rho^*(TE))$. Moreover, from (\ref{Tstru}),
(\ref{4.52'}) and Corollary \ref{c4.2}, it follows that
\begin{equation}\label{e4.54'}
\begin{array}{lcl}
\lcf T^{\rho}X, T^{\rho}Y \rcf^T_{\Psi} = T^{\rho}\lcf X, Y\rcf ,
\makebox[.2cm]{} \lcf T^{\rho}X, \hat{Y}^{\rho} \rcf^T_{\Psi} =
\widehat{\lcf X, Y\rcf}^{\rho}, \makebox[.2cm]{} \lcf
\hat{X}^{\rho}, \hat{Y}^{\rho} \rcf^{T}_{\Psi} = 0, \\
\rho_\Psi^T(T^\rho X) = X^c, \makebox[.2cm]{}
\rho^T_\Psi(\hat{X}^\rho) = X^v,
\end{array}
\end{equation}
for $X, Y \in \Gamma(E)$. In particular,
\begin{equation}\label{Tlhstru}
\begin{array}{lcl}
\lcf T^\rho e_\alpha, T^\rho e_\beta\rcf_\Psi^T = T^\rho\lcf
e_\alpha,e_\beta\rcf, \;\;\; \lcf T^\rho
e_\alpha,\hat{e}_\beta^\rho\rcf_\Psi^T = \widehat{\lcf
e_\alpha,e_\beta\rcf^\rho}, \;\;\;
\lcf\hat{e}_\alpha^\rho,\hat{e}_\beta^\rho\rcf_\Psi^T = 0,
\\ \rho_\Psi^T(T^\rho e_\alpha) = e_\alpha^c, \;\;\;
\rho^T_\Psi(\hat{e_\alpha}^\rho)=e_\alpha^v.
\end{array}
\end{equation}

Therefore,
\begin{equation}\label{roTpsi}
\rho^T_\Psi(x^i,y^\alpha;z^\alpha,v^\alpha)=(x^i,z^\alpha;
\rho^i_\alpha y^\alpha, v^\gamma + C_{\alpha\beta}^\gamma y^\beta
z^\alpha),
\end{equation}
\begin{equation}\label{AnLieT}
\begin{array}{rcl} \lcf T^\rho e_\alpha,T^\rho
e_\beta\rcf_\Psi^T (x^{i}, z^{\alpha}) &=&
C_{\alpha\beta}^\gamma(T^\rho e_\gamma) + (\rho_\nu^i
\displaystyle\frac{\partial C_{\alpha\beta}^\gamma}{\partial
x^i}z^\nu)\hat{e}_\gamma^\rho ,\\ \lcf T^\rho
e_\alpha,\hat{e}_\beta^\rho\rcf_\Psi^T (x^{i}, z^{\alpha})
&=&C_{\alpha\beta}^\gamma \hat{e}_\gamma^\rho,\;\;\; \lcf
\hat{e}_\alpha^\rho,\hat{e}_\beta^\rho\rcf_\Psi^T (x^{i},
z^{\alpha}) = 0.
\end{array}
\end{equation}

\begin{remark}
{\rm  If $E$ is the standard Lie algebroid $TM$ then $\rho^*(TE)$
is the
 tangent bundle to $TM$ and the vector bundle $\rho^*(TE)\to TM$
 is $T(TM)$ with vector bundle projection $T\tau_M:T(TM)\to TM$,
 where $\tau_M:TM\to M$ is the canonical projection. Moreover,
 following \cite{MX}, one may consider the tangent Lie algebroid
 structure $([\cdot,\cdot]^T,\sigma_{TM})$ on the vector bundle
 $T\tau_M:T(TM)\to TM$ and it is easy to prove that the Lie
 algebroid structures $([\cdot,\cdot]^T,\sigma_{TM})$ and
 $([\cdot,\cdot]_\Psi^T, (Id)_\Psi^T)$ coincide. }
 \end{remark}

Next, we will prove that the canonical involution is a morphism of
Lie algebroids.

\begin{theorem}\label{t4.4}
The canonical involution $\sigma_E$ is the unique Lie algebroid
morphism $\sigma_E:{\cal L}^\tau E\to \rho^*(TE)$ between the Lie
algebroids $({\cal L}^\tau E,\lcf\cdot,\cdot\rcf^\tau,\rho^\tau)$
and $(\rho^*(TE), \lcf\cdot,\cdot\rcf_\Psi^T,\rho_\psi^T)$ such
that $\sigma_E$ is an involution, that is, $\sigma_E^2=Id.$
\end{theorem}
\begin{proof}
Using (\ref{Xvc}), (\ref{TtilT}), (\ref{sigmaE.coordenadas}) and
(\ref{TroXXro}), we deduce that
\begin{equation}\label{e4.58'}
 \sigma_E\circ X^{\bf c}=T^\rho X,\;\;\;
\sigma_E\circ X^{\bf v}=\hat{X}^\rho, \end{equation}
 for $X\in
\Gamma(E)$, where $X^{\bf c}$ and $X^{\bf v}$ are the complete and
vertical lift, respectively, of $X$ to ${\cal L}^\tau E.$ Thus,
from (\ref{Defvc}), (\ref{cvstru1}) and (\ref{e4.54'}), we obtain
that $\sigma_E$ is a morphism between the Lie algebroids $({\cal
L}^\tau E,\lcf\cdot,\cdot\rcf^\tau,\rho^\tau)$ and
$(\rho^*(TE),\lcf\cdot,\cdot\rcf_\Psi^T,\rho^T_\Psi)$.

Moreover, if $\sigma'_E$ is another morphism which satisfies the
above conditions then, using that $\sigma_E'$ is a vector bundle
morphism  and the fact that $\sigma'_E$ is an involution, we
deduce that the local expression of $\sigma'_E$ is
\[
\sigma'_E(x^i, y^\gamma;
z^\gamma,v^\gamma)=(x^i,z^\gamma;y^\gamma,g^\gamma),
\]
where $g^\gamma$ is a linear function in the coordinates
$(z^\gamma,v^\gamma)$. Now, since
\[
\rho^T_\Psi\circ \sigma'_E=\rho^\tau
\]
we conclude that (see (\ref{ancprol}) and (\ref{roTpsi}))
\[
g^\gamma=v^\gamma + C_{\alpha\beta}^\gamma z^\alpha y^\beta,\;\;\;
\mbox{ for all } \gamma.
\]
Therefore, $\sigma_E=\sigma'_E.$
\end{proof}

\setcounter{equation}{0}
\section{Tulczyjew's triple on Lie algebroids}\label{seccion5}
Let $(E,\lcf\;,\;\rcf,\rho)$ be a Lie algebroid of rank $n$ over a
manifold $M$ of dimension $m$. Denote by $\tau:E\to M$ the vector
bundle projection of $E$ and by $\tau^*:E^*\to M$ the vector
bundle projection  of the dual bundle $E^*$ to $E$. Then, $TE^*$
is a vector bundle over $TM$ with vector bundle projection
$T\tau^*:TE^*\to TM$ and we may consider the pullback vector
bundle $\rho^ *(TE^*)$ of $T\tau^*:TE^*\to TM$ over $\rho$, that
is, \[ \rho^*(TE^*)=\{(b,v)\in E\times
TE^*/\rho(b)=(T\tau^*)(v)\}.
\]
It is clear that $\rho^*(TE^*)$ is the prolongation ${\cal
L}^{\tau^*}E$ of $E$ over $\tau^*$ (see Section \ref{seccion3.1}).
Thus, ${\cal L}^{\tau^*}E=\rho^*(TE^*)$ is a vector bundle over
$E$ with vector bundle projection
\[
pr_1:\rho^*(TE^*)\to E,\;\;\; (b,v)\to pr_1(b,v)=b.
\]
Note that if $b_0\in E$ then
\[
(pr_1)^{-1}(b_0)\cong (T\tau^*)^{-1}(\rho(b_0)).
\]
Moreover, if $(x^i)$ are coordinates on $M$, $\{e_\alpha\}$ is a
local basis of $\Gamma(E)$, $(x^i,y^\alpha)$ are the corresponding
coordinates on $E$ and $(x^i,y_\alpha;z^\alpha,v_\alpha)$ are the
corresponding coordinates on $\rho^*(TE^*)\equiv {\cal
L}^{\tau^*}E,$ then the local expression of the vector bundle
projection $pr_1:\rho^*(TE^*)\equiv {\cal L}^{\tau^*}E\to E$, the
sum and the product by real numbers in the vector bundle
$pr_1:\rho^*(TE^*)\equiv {\cal L}^{\tau^*}E\to E$ are:
\[
\begin{array}{rcl}
pr_1(x^i,y_\alpha;z^\alpha,v_\alpha)&=&(x^i,z^\alpha),\\(x^i,y_\alpha;z^\alpha,v_\alpha)\oplus
(x^i,\bar{y}_\alpha;z^\alpha,\bar{v}_\alpha)&=&(x^i,y_\alpha+\bar{y}_\alpha;z^\alpha,v_\alpha+
\bar{v}_\alpha),\\
\lambda(x^i,y_\alpha;z^\alpha,v_\alpha)&=&(x^i,\lambda
y_\alpha;z^\alpha,\lambda v_\alpha),
\end{array}
\]
for
$(x^i,y_\alpha;z^\alpha,v_\alpha),(x^i,\bar{y}_\alpha;z^\alpha,\bar{v}_\alpha)\in
\rho^*(TE^*)$ and $\lambda\in \R.$

Next, we consider the following vector bundles:

\begin{itemize}
\item $(\tau^\tau)^*:({\cal L}^\tau E)^*\to E$ (the dual vector bundle to the
prolongation of $E$ over $\tau$).
\item $pr_1:\rho^*(TE^*)\to E$ (the pull-back vector bundle of
$T\tau^*:TE^*\to TM$ over $\rho).$
\item
$\tau^{\tau^*}:{\cal L}^{\tau^*}E \to E^*$ (the prolongation of
$E$ over $\tau^*:E^*\to M)$.

\item $(\tau^{\tau^*})^{*}:({\cal L}^{\tau^*}E)^*\to E^*$ (the dual
vector bundle to $\tau^{\tau^*}:{\cal L}^{\tau^*}E\to E^*).$

\end{itemize}

Then, the aim of this section is to introduce two vector bundle
isomorphisms
\[
A_E:\rho^*(TE^*)\to ({\cal L}^\tau E)^*,\;\;\; \flat_{E^*}:{\cal
L}^{\tau^*}E\to ({\cal L}^{{\tau^*}} E)^*
\]
is such a way that the following diagram to be commutative
\begin{equation}\label{Diag}
\begin{picture}(480,90)(150,20)
\put(250,20){\makebox(0,0){$E$}} \put(310,70){\vector(-3,-2){60}}
 \put(480,70){\vector(-3,-2){60}}
\put(420,20){\makebox(0,0){$E^*$}} \put(170,40){$(\tau^\tau)^*$}
\put(360,40){$\tau^{\tau^*}$}
\put(460,40){$(\tau^{\tau^*})^*$}\put(180,70){\vector(3,-2){60}}
\put(290,40){$pr_1$} \put(350,70){\vector(3,-2){60}}
\put(180,80){\makebox(0,0){$({\cal L}^\tau E)^* $}}
\put(250,85){$A_E$}\put(400,85){$\flat_{E^*}$}\put(280,80){\vector(-1,0){80}}
\put(310,80){\makebox(0,0){$\rho^*(TE^*)\equiv$}}
\put(350,80){\makebox(0,0){${\cal L}^{\tau^*}E$}}
\put(370,80){\vector(1,0){80}} \put(480,80){\makebox(0,0){$({\cal
L}^{\tau^*}E)^*$}}\end{picture}\end{equation}

First, we will define the isomorphism $\flat_{E^*}$.

If $\Omega_E$ is the canonical symplectic section of ${\cal
L}^{\tau^*}E$ then $\flat_{E^*}:{\cal L}^{\tau^*}E\to ({\cal
L}^{\tau^*}E)^*$ is the vector bundle isomorphism induced by
$\Omega_E$, that is,
\[
\flat_{E^*}(b,v)=i(b,v)(\Omega_E(\tau^{\tau^*}(b,v))),
\]
for $(b,v)\in {\cal L}^{\tau^*}E.$

If $(x^i)$ are local coordinates on $M$, $\{e_\alpha\}$ is a local
basis of sections of $\Gamma(E)$ and we consider the local
coordinates $(x^i,y_\alpha;z^\alpha,v_\alpha)$ on ${\cal
L}^{\tau^*}E$ (see Section \ref{seccion3.1}) and the corresponding
coordinates on the dual vector bundle $({\cal L}^{\tau^*}E)^*$
then, using (\ref{locsym}), we deduce that the local expression of
$\flat_{E^*}$ in these coordinates is
\[
\flat_{E^*}(x^i,y_\alpha;z^\alpha,v_\alpha)=(x^i,y_\alpha;-v_\alpha-C_{\alpha\beta}^\gamma
y_\gamma z^\beta,z^\alpha )
\]
where $C_{\alpha\beta}^\gamma$ are the structure functions of the
Lie bracket $\lcf\cdot,\cdot\rcf$ with respect to the basis
$\{e_\alpha\}$.

Now, we define the vector bundle isomorphism
\[
A_E:\rho^*(TE^*)\to ({\cal L}^\tau E)^*
\]
as follows.

Let $<\cdot,\cdot>:E\times_ME^*\to \R$ be the natural pairing
given by
\[
<a,a^*>=a^*(a),
\]
for $a\in E_x$ and $a^*\in E^*_x$, with $x\in M.$ If $b\in E$ and
\[
(b,u_a)\in \rho^*(TE)_b,\;\; (b,v_{a^*})\in \rho^*(TE^*)_b
\]
then
\[
(u_a,v_{a^*})\in T_{(a,a^*)}(E\times_M E^*)=\{(u'_{a},v_{a^*}')\in
T_aE\times T_{a^*}E^*/(T_a\tau)(u_a')=(T_{a^*}\tau^*)(v_{a^*}')\}
\]
and we may consider the map
\[
\widetilde{T{<\cdot,\cdot>}}:\rho^*(TE)\times_E\rho^*(TE^*)\to \R
\]
defined by
\[
\widetilde{T{<\cdot,\cdot>}}((b,u_a),(b,v_{a^*}))=dt_{<a,a^*>}
((T_{(a,a^*)}<\cdot,\cdot>)(u_a,v_{a^*})),
\]
where $t$ is the usual coordinate on $\R$ and
$T{<\;,\;>}:T(E\times_M E^*)\to T\R$ is the tangent map to
$<\;,\;>:E\times_M E^*\to \R$.

If $(x^i,y^\alpha;z^\alpha,v^\alpha)$ are local coordinates on
${\cal L}^\tau E\equiv \rho^*(TE)$ as in Section
\ref{seccion2.2.1} then the local expression of
$\widetilde{T<\cdot,\cdot>}$ is
\begin{equation}\label{Tilpai1}
\widetilde{T{<\cdot,\cdot>}}((x^i,y^\alpha;z^\alpha,v^\alpha),
(x^i,y_\alpha;z^\alpha,v_\alpha))=y^\alpha v_\alpha + y_\alpha
v^\alpha.
\end{equation}
 Thus, the pairing $\widetilde{T{<\cdot,\cdot>}}$ is non-singular
 and, therefore it induces an isomorphism between the vector
 bundles $\rho^*(TE)\to E$ and $\rho^*(TE^*)^*\to E$ which we also
 we denote by $\widetilde{T{<\cdot,\cdot>}}$, that is,
 \[
 \widetilde{T{<\cdot,\cdot>}}:\rho^*(TE)\to \rho^*(TE^*)^*.
 \]
 From (\ref{Tilpai1}), it follows that
 \begin{equation}\label{Tilpai2}
 \widetilde{T{<\cdot,\cdot>}}(x^i,y^\gamma;z^\gamma,
 v^\gamma)=(x^i,v^\gamma;z^\gamma,y^\gamma).
 \end{equation}
 Next, we consider the isomorphism of vector bundles $A^*_E:{\cal
 L}^\tau E\to \rho^*(TE^*)^*$ given by
 \[
 A_E^*=\widetilde{T{<\cdot,\cdot>}}\circ \sigma_E
 \]
 $\sigma_E:{\cal L}^\tau E\to \rho^*(TE)$ being the canonical
 involution introduced in Section \ref{seccion4}. Using
 (\ref{Tilpai2}) and Theorem \ref{t4.4}, we deduce that the local
 expression of the map $A^*_E$ is
 \begin{equation}\label{AEs}
A_E^*(x^i,y^\gamma;z^\gamma,v^\gamma)= (x^i,v^\gamma+
C_{\alpha\beta}^{\gamma}y^\beta z^\alpha ;y^\gamma,z^\gamma).
\end{equation}

Finally, the isomorphism $A_E:\rho^*(TE^*)\to ({\cal L}^\tau E)^*$
between the vector bundles $\rho^*(TE^*)\to E$ and $({\cal L}^\tau
E)^*\to E$ is just the dual map to $A_E^*:{\cal L}^\tau E\to
\rho^*(TE^*)^*.$ From (\ref{AEs}) we conclude that
 \begin{equation}\label{AE}
A_E(x^i,y_\alpha;z^\alpha,v_\alpha)= (x^i,z^\alpha; v_\alpha+
C_{\alpha\beta}^{\gamma}y_\gamma z^\beta ,y_\alpha).
\end{equation}

Diagram (\ref{Diag}) will be called {\it Tulczyjew's triple for
the Lie algebroid} $E$.

\begin{remark}
{\rm If $E$ is the standard Lie algebroid $TM$ then the vector
bundle isomorphisms $A_{TM}:T(T^*M)\to T^*(TM)$ and
$\flat_{T^*M}:T(T^*M)\to T^*(T^*M)$ were considered by Tulczyjew
\cite{Tul1,Tul2} and Diagram (\ref{Diag}) is just Tulczyjew's
triple.}
\end{remark}

\setcounter{equation}{0}
\section{The prolongation of a symplectic Lie
algebroid}\label{seccion6} First of all, we will introduce the
definition of a symplectic Lie algebroid.
\begin{definition}
A Lie algebroid $(E,\lcf\cdot,\cdot\rcf, \rho)$ over a manifold
$M$ is said to be symplectic if it admits a symplectic section
$\Omega$, that is, $\Omega$ is a section of the vector bundle
$\wedge^2E^*\to M$ such that:
\begin{enumerate}
\item For all $x\in M$, the $2$-form $\Omega(x):E_x\times E_x\to
\R$ on the vector space $E_x$ is nondegenerate and
\item
$\Omega$ is a $2$-cocycle, i.e., $d^{E}\Omega=0.$
\end{enumerate}\end{definition}
\begin{examples} {\rm \begin{enumerate}
\item Let $(M,\Omega)$ be a symplectic manifold. Then the tangent
bundle $TM$ to $M$ is a symplectic Lie algebroid.

\item Let $(E,\lcf\cdot,\cdot\rcf,\rho)$ be an arbitrary Lie
algebroid and $(\lcf\cdot,\cdot\rcf^{\tau^*},\rho^{\tau^*})$ be
the Lie algebroid structure on the prolongation ${\cal
L}^{\tau^*}E$ of $E$ over the bundle projection $\tau^*:E^*\to M$
of the dual vector bundle $E^*$ to $E$. Then $({\cal
L}^{\tau^*}E,\lcf\cdot,\cdot\rcf^{\tau^*},\rho^{\tau^*})$ is a
symplectic Lie algebroid and $\Omega_E$ is a symplectic section of
${\cal L}^{\tau^*} E,$ $\Omega_E$ being the canonical symplectic
$2$-section of ${\cal L}^{\tau^*}E$ (see Section
\ref{seccion3.2}).

\end{enumerate}}
\end{examples}

In this section, we will prove that the prolongation of a
symplectic Lie algebroid over the vector bundle projection is a
symplectic Lie algebroid.

For this purpose, we will need some previous results.

Let $(E,\lcf\cdot,\cdot\rcf,\rho)$ be a Lie algebroid over a
manifold $M$ and $\tau:E\to M$ be the vector bundle projection.
Denote by $(\lcf\cdot,\cdot\rcf^\tau,\rho^\tau)$ the Lie algebroid
structure on the prolongation ${\cal L}^\tau E$ of $E$ over
$\tau$. If $pr_1:E\times TE\to E$ is the canonical projection on
the first factor then the pair $(pr_{1|{\cal L}^\tau E}, \tau)$ is
a morphism between the Lie algebroids $({\cal L}^\tau E,
\lcf\cdot,\cdot\rcf^\tau,\rho^\tau)$ and
$(E,\lcf\cdot,\cdot\rcf,\rho)$ (see Section \ref{seccion2.1.1}).
Thus, if $\alpha\in \Gamma(\wedge^k E^*)$ we may consider the
section $\alpha^{\bf v}$ of the vector bundle $\wedge^k({\cal
L}^\tau E)^*\to M$ defined by
\begin{equation}\label{veralpha}
\alpha^{\bf v}=(pr_{1|{\cal L}^\tau E}, \tau)^*( \alpha).
\end{equation}
$\alpha^{\bf v}$ is called {\it the vertical lift } to ${\cal
L}^\tau E$ of $\alpha$ and we have that
\begin{equation}\label{difver}
d^{{\cal L}^\tau E}\alpha^{\bf v}=(d^E \alpha)^{\bf v}.
\end{equation}

Note that if $\{e_\alpha\}$ is a local basis of $\Gamma(E)$ and
$\{e^\alpha\}$ is the dual basis to $\{e_\alpha\}$ then
\[
(e^\alpha)^{\bf v}((e_\beta)^{\bf
c})=\delta_{\alpha\beta},\;\;\;\; (e^\alpha)^{\bf
v}((e_\beta)^{\bf v})=0,
\]
for all $\alpha$ and $\beta$. Moreover, if $\gamma\in \Gamma
(\wedge^k E^*)$ and
\[
\gamma =\gamma_{\alpha_1,\dots ,\alpha_k}e^{\alpha_1}\wedge\dots
\wedge e^{\alpha_k}
\]
we have that
\[
\gamma^{\bf v}=(\gamma_{\alpha_1,\dots ,\alpha_k} \circ \tau )
(e^{\alpha_1})^{\bf v}\wedge\dots \wedge (e^{\alpha_k})^{\bf v}.
\]

Now, suppose that $f \in C^{\infty}(M)$. Then, one may consider
the complete lift $f^c$ of $f$ to $E$. $f^c$ is a real
$C^{\infty}$-function on $E$ (see (\ref{fcv})). This construction
may be generalized as follows.

\begin{proposition}
If $\alpha$ is a section of the vector bundle $\wedge^k E^*\to M$,
then there exists a unique section $\alpha^{\bf c}$ of the vector
bundle $\wedge^k({\cal L}^\tau E)^*\to E$ such that
\begin{equation}\label{complete}
\begin{array}{rcl}
\alpha^{\bf c}(X_1^{\bf c},\dots , X_k^{\bf c})&=&\alpha(X_1,\dots
,X_k)^{c},\\ \alpha^{\bf c}(X_1^{\bf v},X_2^{\bf c}\dots ,
X_k^{\bf c})&=&\alpha(X_1,\dots ,X_k)^v,\\\alpha^{\bf c}(X_1^{\bf
v},\dots ,X_s^{\bf v},X_{s+1}^{\bf c},\dots , X_k^{\bf c})&=&0,
\;\;\; \mbox{ if } 2\leq s\leq k,
\end{array}
\end{equation}
for $X_1,\dots ,X_k\in \Gamma(E).$ Furthermore,
\begin{equation}\label{difcom}
d^{{\cal L}^\tau E}\alpha^{\bf c}=(d^E \alpha)^{\bf c}.
\end{equation}
\end{proposition}
\begin{proof}
We recall that if $\{X_i\}$ is a local basis of $\Gamma(E)$ then
$\{X_i^{\bf c},X_i^{\bf v}\}$ is a local basis of $\Gamma({\cal
L}^{\tau} E).$

On the other hand, if $X\in \Gamma(E)$ and $f, g \in C^\infty(M)$
then
\[
\begin{array}{ccc}
(fX)^{\bf c} = f^cX^{\bf v} + f^vX^{\bf c}, & \quad &  (fX)^{\bf
v} =  f^vX^{\bf v},\\ (fg)^c  =  f^c g^v + f^v g^c, & & (fg)^v  =
f^v g^v.
\end{array}
\]
Using the above facts, we deduce the first part of the
proposition.

Now, if $\alpha\in \Gamma(E^*)$ then, from (\ref{cvstru1}) and
(\ref{complete}) it follows that
$$
\begin{array}{rcl}
(d^{{\cal L}^\tau E}\alpha^{\bf c})(X^{\bf c},Y^{\bf c})&=&((d^E
\alpha)(X,Y))^c=(d^E \alpha)^{\bf c}(X^{\bf c},Y^{\bf c})\\
(d^{{\cal L}^\tau E}\alpha^{\bf c})(X^{\bf v},Y^{\bf c})&=&((d^E
\alpha)(X,Y))^{ v}=(d^E \alpha)^{\bf c}(X^{\bf v},Y^{\bf c})\\
(d^{{\cal L}^\tau E}\alpha^{\bf c})(X^{\bf v},Y^{\bf v})&=&0=(d^E
\alpha)^{\bf
c}(X^{\bf v},Y^{\bf v})\\
\end{array}
$$
for $X,Y\in \Gamma(E)$. Thus,
\begin{equation}\label{difcom1}
d^{{\cal L}^\tau E}\alpha^{\bf c}=(d^E \alpha)^{\bf c},\;\;\;
\mbox{ for } \alpha\in \Gamma(E^*).
\end{equation}
In addition, using (\ref{cvstru1}) and (\ref{complete}), we have
that
\begin{equation}\label{6.68'}
d^{{\cal L}^\tau E}f^c = (d^E f)^c, \makebox[.4cm]{} \mbox{ for }
f \in C^{\infty}(M) = \Gamma(\Lambda^0 E^*).
\end{equation}
 Moreover, using (\ref{veralpha}) and (\ref{complete}) we
obtain that
\begin{equation}\label{extcom}
(f \alpha_1\wedge \dots \wedge \alpha_r)^{\bf c}=f^c \alpha_1^{\bf
v}\wedge \dots \wedge \alpha_r^{\bf v} + f^v \displaystyle
\sum_{i=1}^{r} \alpha_1^{\bf v}\wedge \dots \wedge \alpha_i^{\bf
c}\wedge \dots\wedge \alpha_r^{\bf v}
\end{equation}
for $\alpha_i\in \Gamma(\wedge^{k_i}E^*), i\in \{1,\dots ,r\}.$
Therefore, from (\ref{difver}), (\ref{difcom1}), (\ref{6.68'}) and
(\ref{extcom}) we conclude that (\ref{difcom}) holds.
\end{proof}

 The section $\alpha^{\bf c}$ of the vector bundle
$\wedge^k(({\cal L}^TE)^*)\to E$ is called {\it the complete lift
of $\alpha$.}

If $\{e_\alpha\}$ is a local basis of $\Gamma(E)$ and
$\{e^\alpha\}$ is the dual basis to $\{e_\alpha\}$ then
\[
(e^\alpha)^{\bf c}((e_\beta)^{\bf c})=0,\;\;\;\; (e^\alpha)^{\bf
c}((e_\beta)^{\bf v})=\delta_{\alpha\beta},
\]
for all $\alpha$ and $\beta$. Furthermore, if $\gamma\in
\Gamma(\wedge^kE^*)$ and $\gamma=\gamma_{\alpha_1\dots
\alpha_k}e^{\alpha_1}\wedge \dots \wedge e^{\alpha_k},$ we have
that
\[
\gamma^{\bf c}=\gamma_{\alpha_1\dots \alpha_k}^c
(e^{\alpha_1})^{\bf v}\wedge \dots \wedge (e^{\alpha_k})^{\bf v} +
\displaystyle \sum_{i=1}^{k}(\gamma_{\alpha_1\dots \alpha_k}\circ
\tau) (e^{\alpha_1})^{\bf v}\wedge \dots \wedge
(e^{\alpha_i})^{\bf c}\wedge \dots \wedge (e^{\alpha_k})^{\bf v}.
\]
Thus, if $(x^i)$ are local coordinates defined an open subset $U$
of $M$ and $\{e_\alpha\}$ is a local basis of $\Gamma(E)$ on $U$
then
\[
\begin{array}{rcl}
\gamma^{\bf
c}(x^i,y^\alpha)&=&\rho_\beta^i\displaystyle\frac{\partial
\gamma_{\alpha_1\dots \alpha_k}}{\partial x^i} y^\beta
(e^{\alpha_1})^{\bf v}\wedge \dots \wedge (e^{\alpha_k})^{\bf v}\\
&& + \displaystyle\sum_j(\gamma_{\alpha_1\dots \alpha_k}\circ \tau
)(e^{\alpha_1})^{\bf v}\wedge \dots \wedge (e^{\alpha_j})^{\bf
c}\wedge \dots \wedge (e^{\alpha_k})^{\bf v},
\end{array}
\]
where $(x^i,y^\alpha)$ are the corresponding local coordinates on
$E$ and $\rho_\beta^i$ are the components of the anchor map $\rho$
with respect to $(x^i)$ and to $\{e_\alpha\}$.

Note that
\[
\{ (e^\alpha)^{\bf c}, (e^\alpha)^{\bf v}\}
\]
is a local basis of $\Gamma(({\cal L}^\tau E)^*).$ In fact,
$\{(e^\alpha)^{\bf c}, (e^\alpha)^{\bf v}\}$ is the dual basis to
the local basis of $\Gamma({\cal L}^\tau E)$
\[
\{(e_\alpha)^{\bf v}, (e_\alpha)^{\bf c}\}.
\]
\begin{remark}\label{r6.5}
{\rm If $E$ is the standard Lie algebroid $\tau_M:TM\to M$ and
$\alpha$ is a $k$-form on $M$, that is, $\alpha \in
\Gamma(\wedge^k(T^*M))$ then ${\cal L}^{\tau_M} E=T(TM)$ and
$\alpha^{\bf c}$ is the usual complete lift of $\alpha$ to $TM$
(see \cite{YI}).}
\end{remark}
Next, we will prove the main result of this section.
\begin{theorem}\label{t6.4}
Let $(E,\lcf\cdot,\cdot\rcf,\rho)$ be a symplectic Lie algebroid
with symplectic section $\Omega$ and $\tau:E\to M$ be the vector
bundle projection of $E$. Then, the prolongation ${\cal L}^\tau E$
of $E$ over $\tau$ is a symplectic Lie algebroid and the complete
lift $\Omega^{\bf c}$ of $\Omega$ to ${\cal L}^\tau E$ is  a
symplectic section of ${\cal L}^TE$.
\end{theorem}
\begin{proof}
It is clear that
\[
d^{{\cal L}^\tau E}\Omega^{\bf c}=(d^E \Omega)^{\bf c}=0.
\]
Moreover, if $\{e_{\alpha}\}$ is a local basis of $\Gamma(E)$ on
an open subset $U$ of $E$ and $\{e^\alpha\}$ is the dual basis to
$\{e_\alpha\}$ then
\[
\Omega=\frac{1}{2}\Omega_{\alpha\beta} e^\alpha\wedge e^\beta
\]
with $\Omega_{\alpha\beta}=-\Omega_{\beta\alpha}$ real functions
on $U.$

Thus, using (\ref{extcom}), we have that
\[
\Omega^{\bf c}=\frac{1}{2}\Omega^c_{\alpha\beta}(e^\alpha)^{\bf
v}\wedge (e^\beta)^{\bf v} + \Omega_{\alpha\beta}^{
v}(e^\alpha)^{\bf c}\wedge (e^\beta)^{\bf v}.
\]
Therefore, the local matrix associated to $\Omega^{\bf c}$ with
respect to the basis $\{(e^\alpha)^{\bf c}, (e^\alpha)^{\bf v}\}$
is $$\left(\begin{array}{cc}
0&(\Omega_{\alpha\beta})^v\\-(\Omega_{\alpha\beta})^v&
(\Omega_{\alpha\beta})^c \end{array}\right)$$ Consequently, the
rank of $\Omega^{\bf c}$ is $2n$ and $\Omega^{\bf c}$ is
nondegenerate.
\end{proof}

\begin{remark}\label{r6.4'}
{\rm Let $(M, \Omega)$ be a symplectic manifold. Then, using
Theorem \ref{t6.4}, we deduce a well-known result (see
\cite{Tul2}): the tangent bundle to $M$ is a symplectic manifold
and the complete lift $\Omega^c$ of $\Omega$ to $TM$ is a
symplectic $2$-form on $TM$. }
\end{remark}

\begin{example}
{\rm Let $(E,\lcf\cdot,\cdot\rcf,\rho)$ be a Lie algebroid of rank
$n$ over a manifold $M$ of dimension $m$ and $\tau: E \to M$ be
the vector bundle projection. Denote by $\lambda_{E}$ and
$\Omega_{E}$ the Liouville section and the canonical symplectic
section of the prolongation ${\cal L}^{\tau^*}E$ of $E$ over the
vector bundle projection $\tau^*:E^*\to M$ of the dual vector
bundle $E^*$ to $E$.

On the other hand, if $(x^i)$ are coordinates on $M$ and
$\{e_\alpha\}$ is a local basis of $\Gamma(E)$ then we may
consider  the corresponding coordinates $(x^i,y_\alpha)$
(respectively, $(x^i,y_\alpha;z^\alpha,v_\alpha)$) of $E^*$
(respectively, ${\cal L}^{\tau^*}E$) and the corresponding local
basis $\{\tilde{e}_\alpha,\bar{e}_\alpha\}$ of $\Gamma({\cal
L}^{\tau^*}E)$ (see Section \ref{seccion3.1}). Thus, if ${\cal
L}^{\tau^{\tau^*}}({\cal L}^{\tau^*}E)$ is the prolongation of
${\cal L}^{\tau^*}E$ over the vector bundle projection
$\tau^{\tau^*}:{\cal L}^{\tau^*}E\to E^*$ then
\[
\{ \tilde{e}_\alpha^{\bf c},  \bar{e}_\alpha^{\bf c},
\tilde{e}_\alpha^{\bf v}, \bar{e}_\alpha^{\bf v}\}
\]
is a local basis of $\Gamma({\cal L}^{\tau^{\tau^*}}({\cal
L}^{\tau^*}E))$ and if $\{ \tilde{e}^\alpha, \bar{e}^\alpha\}$ is
the dual basis to $\{ \tilde{e}_\alpha, \bar{e}_\alpha\}$ then
\[
\{ (\tilde{e}^\alpha)^{\bf v}, (\bar{e}^\alpha)^{\bf v},
(\tilde{e}^\alpha)^{\bf c}, (\bar{e}^\alpha)^{\bf c}\}
\]
is the dual basis of $\{ \tilde{e}_\alpha^{\bf c},
\bar{e}_\alpha^{\bf c}, \tilde{e}_\alpha^{\bf v},
\bar{e}_\alpha^{\bf v}\}.$

Moreover, using (\ref{locsym}) and (\ref{extcom}), we deduce that
the local expressions of the complete lifts $\lambda_{E}^{\bf c}$
and $\Omega_E^{\bf c}$ of $\lambda_{E}$ and $\Omega_E$ are
\[
\begin{array}{lcl}
\lambda^{\bf c}_E (x^i, y_{\alpha}; z^{\alpha}, v_{\alpha}) & = &
y_{\alpha}^c (\tilde{e}^{\alpha})^{\bf v} + y_{\alpha}^v
(\tilde{e}^{\alpha})^{\bf c}, \\ \Omega_E^{\bf c}(x^i, y_{\alpha};
z^{\alpha}, v_{\alpha}) &=&(\tilde{e}^\alpha)^{\bf c}\wedge
(\bar{e}^\alpha)^{\bf v} + (\tilde{e}^\alpha)^{\bf v}\wedge
(\bar{e}^\alpha)^{\bf c}+
\displaystyle\frac{1}{2}(C_{\alpha\beta}^\gamma y_\gamma)^{
c}(\tilde{e}^\alpha)^{\bf v}\wedge (\tilde{e}^\beta)^{\bf v} \\ &
& + C_{\alpha\beta}^\gamma y_\gamma (\tilde{e}^\alpha)^{\bf
c}\wedge (\tilde{e}^\beta)^{\bf v}.
\end{array}
\]
Therefore, from (\ref{dualstru}), we conclude that
\begin{equation}\label{OmEc}
\begin{array}{lcl}
\lambda^{\bf c}_E (x^i, y_{\alpha}; z^{\alpha}, v_{\alpha}) & = &
y_{\alpha} (\tilde{e}^{\alpha})^{\bf c} +
v_{\alpha}(\tilde{e}^{\alpha})^{\bf v}, \\ \Omega_E^{\bf c} (x^i,
y_{\alpha}; z^{\alpha}, v_{\alpha}) &=&(\tilde{e}^\alpha)^{\bf
c}\wedge (\bar{e}^\alpha)^{\bf v} + (\tilde{e}^\alpha)^{\bf
v}\wedge (\bar{e}^\alpha)^{\bf c} \\ & &
+\displaystyle\frac{1}{2}(\rho_\mu^i\displaystyle\frac{\partial
C_{\alpha\beta}^\gamma}{\partial x^i}z^\mu y_\gamma +
C_{\alpha\beta}^\gamma v_\gamma )(\tilde{e}^\alpha)^{\bf v}\wedge
(\tilde{e}^\beta)^{\bf v}
\\ & & + C_{\alpha\beta}^\gamma y_\gamma
(\tilde{e}^\alpha)^{\bf c}\wedge (\tilde{e}^\beta)^{\bf v},
\end{array}
\end{equation}
$\rho^i_\alpha$ and $C_{\alpha\beta}^\gamma$ being the structure
functions of the Lie algebroid $(E,\lcf\cdot,\cdot\rcf,\rho)$ with
respect to the coordinates $(x^i)$ and to the basis
$\{e_\alpha\}.$ }
\end{example}

\setcounter{equation}{0}
\section{Lagrangian Lie subalgebroids in symplectic Lie
algebroids}

First of all, we will introduce the notion of a Lagrangian Lie
subalgebroid of a symplectic Lie algebroid.

\begin{definition}\label{d7.2}
Let $(E,\lcf\cdot,\cdot\rcf,\rho)$ be a symplectic Lie algebroid
with symplectic section $\Omega$ and $j:F\to E$, $i:N\to M$ be a
Lie subalgebroid (see Remark \ref{prolon-acc}). Then, the Lie
subalgebroid is said to be Lagrangian if $j(F_x)$ is a Lagrangian
subspace of the symplectic vector space
$(E_{i(x)},\Omega_{i(x)}),$ for all $x\in N.$
\end{definition}
Definition \ref{d7.2} implies that:
\begin{enumerate}
\item $rank F=\displaystyle\frac{1}{2}rank E$ and \item
$(\Omega (i(x)))_{|j(F_x)\times j(F_x)}=0,$ for all $x\in N.$
\end{enumerate}

\begin{remark}\label{r7.3}
{\rm Let $(M,\Omega)$ be a symplectic manifold, $S$ be a
submanifold of $M$ and $i: S \to M$ be the canonical inclusion.
Then, the standard Lie algebroid $\tau_M:TM\to M$ is symplectic
and $i: S \to M$, $j = Ti: TS \to TM$ is a Lie subalgebroid of
$\tau_{M}: TM \to M$. Moreover, one may prove that $S$ is a
Lagrangian submanifold of $M$ in the usual sense if and only if
the Lie subalgebroid $i:S\to M, j=Ti:TS\to TM$ of $\tau_M:TM\to M$
is Lagrangian.}
\end{remark}

Let $(E,\lcf\cdot,\cdot\rcf,\rho)$ be a Lie algebroid over $M.$
Then, the prolongation ${\cal L}^{\tau^*}E$ of $E$ over the vector
bundle projection $\tau^*:E^*\to M$ is a symplectic Lie algebroid.
Moreover, if $x$ is a point of $M$ and $E_x^*$ is the fiber of
$E^*$ over the point $x$, we will denote by
\[
j_x:TE_x^*\to {\cal L}^{\tau^*}E,\;\;\;\; i_x:E^*_x\to E^*
\]
the maps given by
\[
j_x(v)=(0(x),v),\;\;\; i_x(\alpha)=\alpha
\]
for $v\in T E^*_x$ and $\alpha\in E_x^*$, where $0:M\to E$ is the
zero section of $E$. Note that if $v\in T E_x^*$, $(T\tau^*)(v)=0$
and thus $(0(x),v)\in {\cal L}^{\tau^*}E$.

On the other hand, if $\gamma\in \Gamma(E^*)$ we will denote by
$F_\gamma$ the vector bundle over $\gamma(M)$ given by
\begin{equation}\label{fgamma}
F_\gamma=\{(b,(T\gamma)(\rho(b)))\in E\times TE^*/b\in E\}
\end{equation}
and by $j_\gamma:F_\gamma\to {\cal L}^{\tau^*}E$ and
$i_\gamma:\gamma(M)\to E^*$ the canonical inclusions. Note that
the pair $((Id,T\gamma\circ \rho),\gamma)$ is an isomorphism
between the vector bundles $E$ and $F_\gamma$, where the map
$(Id,T\gamma\circ \rho):E\to F_\gamma$ is given by
\[
(Id,T\gamma\circ \rho)(b)=(b,T\gamma(\rho(b))),\;\;\; \mbox{ for }
b\in E.
\]
Thus, $F_\gamma$ is a Lie algebroid over $\gamma(M)$. Moreover, we
may prove the following results.

\begin{proposition}\label{p7.4}
Let $(E,\lcf\cdot,\cdot\rcf,\rho)$ be a Lie algebroid of rank $n$
over a manifold $M$ of dimension $m$ and  ${\cal L}^{\tau^*}E$ be
the prolongation of $E$ over the vector bundle projection
$\tau^*:E^*\to M.$
\begin{enumerate}
\item If $x$ is a point of $M$ then $j_x:TE_x^*\to {\cal
L}^{\tau^*}E$ and $i_x:E^*_x\to E^*$ is a Lagrangian Lie
subalgebroid of the symplectic Lie algebroid ${\cal L}^{\tau^*}E.$
\item If $\gamma\in \Gamma(E^*)$ then $j_\gamma:F_\gamma\to  {\cal
L}^{\tau^*}E$ and $i_{\gamma}:\gamma(M)\to E^*$ is a Lagrangian
Lie subalgebroid of the symplectic Lie algebroid ${\cal
L}^{\tau^*}E$ if and only if $\gamma$ is a $1$-cocycle for the
cohomology complex of the Lie algebroid $E$, that is,
$d^{E}\gamma=0.$
\end{enumerate}
\end{proposition}
\begin{proof}
$(i)$ It is clear that the rank of the vector bundle
$\tau_{E_x^*}:TE_{x}^*\to E_x^{*}$ is $n=\frac{1}{2}rank ({\cal
L}^{\tau^*}E).$

Moreover, if $(x^i)$ are local coordinates on $M$, $\{e_\alpha\}$
is a local basis of $\Gamma(E),$ $(x^i,y_\alpha)$ are the
corresponding coordinates on $E^*$ and
$\{\tilde{e}_\alpha,\bar{e}_\alpha\}$ is the corresponding local
basis of $\Gamma({\cal L}^{\tau^*}E)$ then, from (\ref{tilbar}),
it follows that
\begin{equation}\label{7.64'}
\begin{array}{l}
(j_x, i_x)^*(\tilde{e}^{\alpha}) = 0, \makebox[.4cm]{} (j_x,
i_x)^*(\bar{e}^{\alpha}) = d^{TE^*_x}(y_{\alpha} \circ i_{x}), \\
j_x(\displaystyle \frac{\partial }{\partial y_{\alpha}}_{|\mu})
=\bar{e}_{\alpha}(\mu), \makebox[.4cm]{} \mbox{ for all } \mu \in
E^*_x.
\end{array}
\end{equation}

This, using (\ref{dif*}), implies that $j_x:TE_x^*\to {\cal
L}^{\tau^*}E$ and $i_x:E_x^*\to E^*$ is a morphism of Lie
algebroids. Thus, since $j_x$ is injective and $i_x$ is an
injective inmersion, we deduce that $j_x:TE_x^*\to {\cal
L}^{\tau^*}E$ and $i_x:E_x^*\to E^*$ is a Lie subalgebroid of
${\cal L}^{\tau^*}E.$ Finally, from (\ref{locsym}) and
(\ref{7.64'}), we conclude that
\[
(\Omega_E(i_x(\mu)))_{|j_x(T_\mu E_x^*)\times j_x(T_\mu E_x^*)}=0,
\]
for all $\mu \in E_x^*.$

\medskip

$(ii)$ If $\gamma$ is a section of $E^*$ then the Lie algebroids
$E\to M$ and $F_\gamma\to \gamma(M)$ are isomorphic and, under
this isomorphism, the inclusions $j_\gamma$ and $i_\gamma$ are the
maps $(Id,T\gamma\circ \rho):E\to {\cal L}^{\tau^*}E$ and
$\gamma:M\to E^*$, respectively. Furthermore, it clear that the
map $(Id,T\gamma\circ \rho)$ is injective and that $\gamma:M\to
E^*$ is an injective inmersion.

On the other hand, from Theorem \ref{t3.1''}, we have that
\[
((Id,T\gamma\circ \rho),\gamma)^*(\Omega_E)=-d\gamma.
\]
Therefore, the Lie subalgebroid $j_\gamma:F_\gamma\to {\cal
L}^{\tau^*}E$ and $i_\gamma:\gamma(M)\to E^*$  is Lagrangian if
and only if $\gamma$ is a $1$-cocycle.
\end{proof}

\begin{remark}\label{r7.4'}{\rm Using Remark \ref{r7.3} and
applying Proposition \ref{p7.4} to the particular case when $E$ is
the standard Lie algebroid $TM$, we deduce two well-known results
(see, for instance, \cite{AM}):
\begin{enumerate}
\item If $\Omega_{TM}$ is the canonical symplectic $2$-form
on $T^{*}M$ and $x$ is a point of $M$ then the cotangent space to
$M$ at $x$, $T_{x}^{*}M$, is a Lagrangian submanifold of the
symplectic manifold $(T^*M, \Omega_{TM})$.

\item If $\gamma: M \to T^*M$ is a $1$-form on $M$ then the
submanifold $\gamma(M)$ is Lagrangian in the symplectic manifold
$(T^*M, \Omega_{M})$ if and only if $\gamma$ is a closed $1$-form.

\end{enumerate}
}
\end{remark}

Let $(E,\lcf\cdot,\cdot\rcf,\rho)$ be a Lie algebroid of rank $n$
over a manifold $M$ of dimension $m$ and $S$ be a submanifold of
$E$. Denote by $\tau:E\to M$ the vector bundle projection, by
$i:S\to E$ the canonical inclusion and by $\tau^S:S\to M$ the map
given by
\[
\tau^S=\tau\circ i.
\]
If there exists a natural number $c$ such that
\begin{equation}\label{proS}
\dim (\rho(E_{\tau^S(x)}) + (T_x\tau^S)(T_xS))=c,\mbox{ for all
}x\in S,
\end{equation}
then, using the results of Section \ref{seccion2.1.1}, one can
consider the prolongation of the Lie algebroid $E$ over the map
$\tau^S$,
\[
{\cal L}^{\tau^S}E=\{(b,v)\in E\times TS/ \rho(b)=(T\tau^S)(v)\},
\]
which is a Lie algebroid over $S$ of rank $n+s-c,$ where $s=\dim
S$. Moreover, the maps $(Id,Ti):{\cal L}^{\tau^S}E\to{\cal
L}^{\tau}E$ and $i:S\to E$ define a Lie subalgebroid of the
prolongation of $E$ over the bundle projection $\tau$, where $(Id,
Ti)$ is the map given by
\[
(Id, Ti)(b,v)=(b, Ti(v)), \mbox{ for }(b,v)\in {\cal L}^{\tau^S}E.
\]
\begin{examples}\label{e7.4}
{\rm $(i)$ Let $(E,\lcf\cdot,\cdot\rcf,\rho)$ be a Lie algebroid
of rank $n$ over a manifold $M$ and $X$ be a section of $E.$
Suppose that the submanifold $S$ is $X(M)$. Then, it is easy prove
that the condition (\ref{proS}) holds and that $c$ is just $m=\dim
M=\dim S.$ Thus, one may consider the prolongation ${\cal
L}^{\tau^S}E$ of $E$ over the map $\tau^S:S=X(M)\to M$ and the
rank of the vector bundle ${\cal L}^{\tau^S} E\to S$ is $n$.

Note that if $x$ is a point of $M$, the fiber of ${\cal
L}^{\tau^S}E$ over $X(x)\in S$ is the vector space
\[
({\cal L}^{\tau^S}E)_{X(x)}=\{(b,(T_xX)(\rho(b)))/b\in E_x\}.
\]
Therefore, if $Y$ is a section of $E$ then, using (\ref{comverl}),
we deduce that
\[
TX(\rho(Y))=(Y^{c} - \lcf X,Y\rcf^{v})\circ X,
\]
where $Z^c$ (respectively, $Z^v$) denotes the complete
(respectively, vertical) lift of a section $Z$ of $E$ to a section
of the vector bundle $TE\to E$. Consequently (see (\ref{Defvc})),
\begin{equation}\label{vier}
TX(\rho(Y))=\{\rho^\tau(Y^{\bf c} - \lcf X,Y\rcf^{\bf v})\}\circ
X.
\end{equation}
Here, $Z^{\bf c}$ (respectively, $Z^{\bf v}$) is the complete
(respectively, vertical) lift of a section $Z$ of $E$ to a section
of the vector bundle ${\cal L}^\tau E\to E.$ Now, from
(\ref{vier}), we obtain that
\[
Y^{\bf c}_{|S}-\lcf X, Y\rcf^{\bf v}_{|S}
\]
is a section of the vector bundle $S\to {\cal L}^{\tau^s}E.$ Thus,
if $\{e_\alpha\}$ is a local basis of $\Gamma(E)$ it follows that
\[
\{ e_{\alpha|S}^{\bf c} - \lcf X,e_\alpha\rcf^{\bf v}_{|S}\}
\]
is a local basis of $\Gamma({\cal L}^{\tau^S} E).$

\medskip

$(ii)$ Let $(E,\lcf\cdot,\cdot\rcf,\rho)$ be a Lie algebroid of
rank $n$ over a manifold of $M$ of dimension $m$. Denote by ${\cal
L}^{\tau^*} E$ the prolongation of $E$ over the vector bundle
projection $\tau^*:E^*\to M$. ${\cal L}^{\tau^*}E$ is a Lie
algebroid over $E^*$. On the other hand, let $\rho^*(TE^*)$ be the
pull-back of the vector bundle $T\tau^*:TE^*\to TM$ over the
anchor map $\rho:E\to TM$. $\rho^*(TE^*)$ is a vector bundle over
$E$. Moreover, as we know, the total spaces of these vector
bundles coincide, that is,
\[
\rho^*(TE^*)={\cal L}^{\tau^*} E.
\]
Now, suppose that $\tilde{X}$ is a section of the vector bundle
$\rho^*(TE^*)\to E$. Then, $S=\tilde{X}(E)$ is a submanifold of
${\cal L}^{\tau^*}E$ of dimension $m+n$. Furthermore, if we
consider on ${\cal L}^{\tau^*}E$  the Lie algebroid structure
$(\lcf\cdot,\cdot \rcf^{\tau^*},\rho^{\tau^*})$ (see Section
\ref{seccion3.1}) then condition (\ref{proS}) holds for the
submanifold $S$ and the natural number $c$ is $n+m=\dim S$. Thus,
the prolongation ${\cal L}^{(\tau^{\tau^*})^S}({\cal
L}^{\tau^*}E)$ of  the Lie algebroid $({\cal
L}^{\tau^*}E,\lcf\cdot,\cdot\rcf^{\tau^*},\rho^{\tau^*})$ over the
map $\tau^{\tau^*}\circ i=(\tau^{\tau^*})^S:S\to {\cal
L}^{\tau^*}E\to E^*$ is a Lie algebroid over $S$ of rank $2n$. In
fact, if $(x^i)$ are local coordinates on an open subset $U$ of
$M$, $\{e_\alpha\}$ is a basis of $\tau^{-1}(U) \to U$,
$(x,y)\equiv (x^i,y^\alpha)$ (respectively $(x,y;z,v)\equiv
(x^i,y_\alpha;z^\alpha,v_\alpha)$) are the corresponding
coordinates on $E$ (respectively, ${\cal L}^{\tau^*}E)$ and the
local expression of $\tilde{X}$ in this coordinates is
\[
\tilde{X}(x,y)=\tilde{X}(x^i,y^\alpha)=(x^i,\tilde{X}^\alpha;
y^\alpha,\tilde{X}^{'\alpha})\] then
$\{\tilde{e}_\alpha^{\tilde{X}},\bar{e}_\alpha^{\tilde{X}}\}$ is a
local basis of sections of ${\cal L}^{(\tau^{\tau^*})^S}({\cal
L}^{\tau^*}E)\to S$, where
\[
\tilde{e}_\alpha^{\tilde{X}}:S\to {\cal L}^{\tau^*}E\times
TS,\;\;\;\; \bar{e}_\alpha^{\tilde{X}}:S\to {\cal
L}^{\tau^*}E\times TS
\]
are defined by
\begin{equation}\label{bXtilde}
\begin{array}{rcl}
\tilde{e}_\alpha^{\tilde{X}}(\tilde{X}(x,y))&=&
(\tilde{e}_\alpha((\tau^{\tau^*})^S(\tilde{X}(x,y))) +
\rho_\alpha^i(x)\displaystyle\frac{\partial
\tilde{X}^\beta}{\partial
x^i}_{|(x,y)}\bar{e}_\beta((\tau^{\tau^*})^S(\tilde{X}(x,y))),\\&&
\rho_\alpha^i(x)\displaystyle\frac{\partial }{\partial
x^i}_{|\tilde{X}(x,y)} +
\rho_\alpha^i(x)\displaystyle\frac{\partial
\tilde{X}^\beta}{\partial x^i}_{|(x,y)}\displaystyle\frac{\partial
}{\partial y_\beta}_{|\tilde{X}(x,y)}\\&& +
\rho^i_\alpha(x)\displaystyle\frac{\partial
\tilde{X'}^\beta}{\partial x^i}_{|(x,y)}\frac{\partial }{\partial
v_{\beta}}_{|\tilde{X}(x,y)}),\\
\bar{e}_\alpha^{\tilde{X}}(\tilde{X}(x,y))&=&
(\displaystyle\frac{\partial \tilde{X}^\beta}{\partial
y^\alpha}_{|(x,y)}\bar{e}_\beta((\tau^{\tau^*})^S(\tilde{X}(x,y))),
\displaystyle\frac{\partial }{\partial
z^\alpha}_{|\tilde{X}(x,y)}\\&& +\displaystyle\frac{\partial
\tilde{X}^\beta}{\partial
y^\alpha}_{|(x,y)}\displaystyle\frac{\partial }{\partial
y_\beta}_{|\tilde{X}(x,y)} + \displaystyle\frac{\partial
\tilde{X'}^\beta}{\partial
y^\alpha}_{|(x,y)}\displaystyle\frac{\partial }{\partial
v_{\beta}}_{|\tilde{X}(x,y)}).
\end{array}
\end{equation}
Here, $\rho_\alpha^i$ are the components of the anchor map with
respect to $(x^i)$ and $\{e_\alpha\}$ and
$\{\tilde{e}_\alpha,\bar{e}_\alpha\}$ is the corresponding local
basis of $\Gamma({\cal L}^{\tau^*}E)$. Using the local basis
$\{\tilde{e}_{\alpha}^{\tilde{X}}, \bar{e}_{\alpha}^{\tilde{X}}\}$
of $\Gamma({\cal L}^{(\tau^{\tau^*})^S}({\cal L}^{\tau^*}E))$ one
may introduce, in a natural way, local coordinates $(x^i,
y^{\alpha}; z^{\alpha}_{\tilde{X}}, v_{\tilde{X}}^{\alpha})$ on
${\cal L}^{(\tau^{\tau^*})^S}({\cal L}^{\tau^*}E)$ as follows. If
$\omega_{\tilde{X}} \in {\cal L}^{(\tau^{\tau^*})^S}({\cal
L}^{\tau^*}E)_{\tilde{X}(a)}$, with $a \in E$, then $(x^{i},
y^{\alpha})$ are the coordinates of $a$ and
\[
\omega_{\tilde{X}} =
z_{\tilde{X}}^{\alpha}\tilde{e}_{\alpha}^{\tilde{X}}(\tilde{X}(a))
+v_{\tilde{X}}^{\alpha}\bar{e}_{\alpha}^{\tilde{X}}(\tilde{X}(a)).
\]
Moreover, if $(\lcf\cdot,\cdot\rcf^S, \rho^S)$ is the Lie
algebroid structure on the vector bundle ${\cal
L}^{(\tau^{\tau^*})^S}({\cal L}^{\tau^*}E) \to S$ then, using
(\ref{2.10'}), (\ref{2.10''}), (\ref{dualstru}) and
(\ref{bXtilde}), we obtain that
\begin{equation}\label{e7.75'}
\begin{array}{lcl}
\lcf \tilde{e}_{\alpha}^{\tilde{X}}, \tilde{e}_{\beta}^{\tilde{X}}
\rcf^S = C_{\alpha \beta}^{\gamma} \tilde{e}_{\gamma}^{\tilde{X}},
\makebox[.4cm]{} \lcf \tilde{e}_{\alpha}^{\tilde{X}},
\bar{e}_{\beta}^{\tilde{X}} \rcf^S = \lcf
\bar{e}_{\alpha}^{\tilde{X}}, \bar{e}_{\beta}^{\tilde{X}} \rcf^S =
0, \\ \rho^{S}(\tilde{e}_{\alpha}^{\tilde{X}}) =
(T\tilde{X})(\rho^i_{\alpha} \displaystyle
\frac{\partial}{\partial x^i}), \makebox[.4cm]{}
\rho^{S}(\bar{e}_{\alpha}^{\tilde{X}}) =
(T\tilde{X})(\displaystyle \frac{\partial}{\partial y^\alpha}),
\end{array}
\end{equation}
for all $\alpha$ and $\beta$. Now, we consider the map
$\Theta^{\tilde{X}}: {\cal L}^{\tau}E \to {\cal
L}^{(\tau^{\tau^*})^S}({\cal L}^{\tau^*}E)$ defined by
\[
\Theta^{\tilde{X}}(a, v) = ((a,
(T_{\tilde{X}(b)}(\tau^{\tau^*})^S)((T_{b}\tilde{X})(v))),
(T_{b}\tilde{X})(v)),
\]
for $(a, v) \in ({\cal L}^{\tau}E)_{b} \subseteq E_{\tau(b)}
\times T_{b}E$, with $b \in E$. Note that $\tau^* \circ
(\tau^{\tau^*})^S \circ \tilde{X} = \tau $ and, thus, $(a,
(T_{\tilde{X}(b)}(\tau^{\tau^*})^S)((T_{b}\tilde{X})(v))) \in
({\cal L}^{\tau^*}E)_{(\tau^{\tau^*})^S(\tilde{X}(b))} \subseteq
E_{\tau(b)} \times T_{\tau^{\tau^*}(\tilde{X}(b))}E^*$ which
implies that
\[
\begin{array}{lcl}
((a, (T_{\tilde{X}(b)}(\tau^{\tau^*})^S)((T_{b}\tilde{X})(v))),
(T_{b}\tilde{X})(v)) \in {\cal L}^{(\tau^{\tau^*})^S}({\cal
L}^{\tau^*}E)_{\tilde{X}(b)} \\ \subseteq ({\cal
L}^{\tau^*}E)_{(\tau^{\tau^*})^S(\tilde{X}(b))} \times
T_{\tilde{X}(b)}S.
\end{array}
\]
Furthermore, if $\{\tilde{T}_{\alpha}, V_{\alpha}\}$ is the local
basis of $\Gamma({\cal L}^{\tau}E)$ considered in Remark
\ref{r2.3} then a direct computation, using (\ref{bXtilde}),
proves that
\begin{equation}\label{7.75''}
\Theta^{\tilde{X}}(\tilde{T}_{\alpha}(a)) =
\tilde{e}_{\alpha}^{\tilde{X}}(\tilde{X}(a)), \makebox[.4cm]{}
\Theta^{\tilde{X}}(V_{\alpha}(a)) =
\bar{e}_{\alpha}^{\tilde{X}}(\tilde{X}(a)),
\end{equation}
for all $a \in \tau^{-1}(U)$. Therefore, from (\ref{tilde}),
(\ref{e7.75'}) and (\ref{7.75''}), we conclude that the pair
$(\Theta^{\tilde{X}}, \tilde{X})$ is an isomorphism between the
Lie algebroids ${\cal L}^{\tau}E \to E$ and ${\cal
L}^{(\tau^{\tau^*})^S}({\cal L}^{\tau^*}E) \to S$. Note that the
local expression of $\Theta^{\tilde{X}}$ in the local coordinates
$(x^i, y^{\alpha}; z^{\alpha}, v^{\alpha})$ and $(x^i, y^{\alpha};
z^{\alpha}_{\tilde{X}}, v^{\alpha}_{\tilde{X}})$ on ${\cal
L}^{\tau}E$ and ${\cal L}^{(\tau^{\tau^*})^S}({\cal L}^{\tau^*}E)$
is just the identity, that is,
\[
\Theta^{\tilde{X}}(x^i, y^{\alpha}; z^{\alpha}, v^{\alpha}) =
(x^i, y^{\alpha}; z^{\alpha}, v^{\alpha}).
\]
}
\end{examples}

We recall that if $E$ is a symplectic Lie algebroid with
symplectic section $\Omega$ then the prolongation ${\cal L}^\tau
E$ of $E$ over the vector bundle projection $\tau:E\to M$ is a
symplectic Lie algebroid with symplectic section $\Omega^{\bf c}$,
the complete lift of $\Omega$ (see Section \ref{seccion6}).

\begin{proposition}\label{p7.5}
Let $(E,\lcf\cdot,\cdot\rcf,\rho)$ be a symplectic Lie algebroid
with symplectic section $\Omega$ and $X$ be a section of $E.$
Denote by $S$ the submanifold of $E$ defined by $S=X(E)$, by
$i:S\to E$ the canonical inclusion, by $\alpha_X$  the section of
$E^*$ given by $\alpha_X=i_{X}\Omega$ and by $\tau^S:S\to M$ the
map defined by $\tau^S=\tau\circ i$, $\tau:E\to M$ being the
vector bundle projection. Then, the Lie subalgebroid
$(Id,Ti):{\cal L}^{\tau^S}E\to {\cal L}^\tau E, i:S\to E$ of the
symplectic Lie algebroid $(({\cal
L}^{\tau}E,\lcf\cdot,\cdot\rcf^\tau,\rho^\tau),\Omega^{\bf c})$ is
Lagrangian if and only if $\alpha_X$ is a $1$-cocycle.
\end{proposition}
\begin{proof}
Suppose that $Y$ and $Z$ are sections of $E$. Since $\Omega$ is a
$2$-cocycle, it follows that (see (\ref{dif}))
\begin{equation}\label{difalpha}
d\alpha_X(Y,Z)=\rho(X)(\Omega(Y,Z))-\Omega(Y,\lcf
X,Z\rcf)-\Omega(\lcf X, Y\rcf,Z).
\end{equation}

On the other hand, using (\ref{complete}), we obtain that
\[
\Omega^{\bf c}(Y^{\bf c} - \lcf X,Y\rcf^{\bf v},Z^{\bf c}-\lcf
X,Z\rcf^{\bf v})=\Omega(Y,Z)^{c} - \Omega(Y,\lcf X,Z\rcf)^{v} -
\Omega(\lcf X,Y\rcf, Z)^{v}
\]
which implies that (see (\ref{fcv}))
\begin{equation}\label{Omegac}
\Omega^{\bf c}(Y^{\bf c} - \lcf X,Y\rcf^{\bf v},Z^{\bf c}- \lcf
X,Z\rcf^{\bf v})\circ X=\rho(X)(\Omega(Y,Z))-\Omega(Y,\lcf
X,Z\rcf)-\Omega(\lcf X,Y\rcf,Z).
\end{equation}

Therefore, from (\ref{difalpha}), (\ref{Omegac}) and taking into
account that the rank of ${\cal L}^{\tau^S}E$ is $n=\frac{1}{2}
rank ({\cal L}^\tau E)$, we conclude that the Lie subalgebroid
$(Id,Ti):{\cal L}^{\tau^S}E\to {\cal L}^\tau E,$ $i:S\to E$ is
Lagrangian if and only if $\alpha_X$ is a cocycle (see Example
\ref{e7.4}, $(i)$).
\end{proof}

\begin{remark}{\rm Let $(M,\Omega)$ be a symplectic manifold,
$\tau_M:TM\to M$ be the standard symplectic Lie algebroid and $X$
be a vector field on $M$. Then, the tangent bundle $TM$ of $M$ is
a symplectic manifold with symplectic form the complete lift
$\Omega^{\bf c}$ of $\Omega$ to $TM$ (see Remark \ref{r6.4'}).
Moreover, using Remark \ref{r7.3} and Proposition \ref{p7.5}, we
deduce a well-known result: the submanifold $X(M)$ of $TM$ is
Lagrangian if and only if $X$ is a locally Hamiltonian vector
field of $M.$ }
\end{remark}

Let $(E,\lcf\cdot,\cdot,\rcf,\rho)$ be a Lie algebroid over a
manifold $M$ and denote by $\Omega_E$ the canonical symplectic
section of the Lie algebroid $({\cal
L}^{\tau^*}E,\lcf\cdot,\cdot\rcf^{\tau^*},\rho^{\tau^*})$ (see
Sections \ref{seccion3.1} and \ref{seccion3.2}) and by
$\rho^*(TE^*)\to E$ the pull-back of the vector bundle
$T\tau^*:TE^*\to TM$ over the anchor map $\rho:E\to TM$. As we
know, $\rho^*(TE^*)={\cal L}^{\tau^*}E$. Now, suppose that
$\tilde{X}:E\to \rho^*(TE^*)={\cal L}^{\tau^*}E$ is a section of
$\rho^*(TE^*)\to E.$ Then, $S=\tilde{X}(E)$ is a submanifold of
${\cal L}^{\tau^*}E$ and one may consider the Lie subalgebroid
$(Id,Ti):{\cal L}^{(\tau^{\tau^*})^S}({\cal L}^{\tau^*}E)\to {\cal
L}^{\tau^{\tau^*}}({\cal L}^{\tau^*}E)$, $i:S\to {\cal
L}^{\tau^*}E$ of the symplectic Lie algebroid ${\cal
L}^{\tau^{\tau^*}}({\cal L}^{\tau^*}E)\to {\cal L}^{\tau^*}E$ (see
Example \ref{e7.4}, $(ii)$). We remark that the symplectic section
of ${\cal L}^{\tau^{\tau^*}}({\cal L}^{\tau^*}E)\to {\cal
L}^{\tau^*}E$ is the complete lift $\Omega_E^{\bf c}$ of
$\Omega_E$ to the prolongation of ${\cal L}^{\tau^*}E$ over the
bundle projection $\tau^{\tau^*}:{\cal L}^{\tau^*}E\to E^*$.

On the other hand, let $A_E:\rho^*(TE^*)\to ({\cal L}^\tau E)^*$
be the canonical isomorphism between the vector bundles
$\rho^*(TE^*)\to E$ and $({\cal L}^\tau E)^*\to E$ considered in
Section \ref{seccion5} (see (\ref{AE}))
 and $\alpha_{\tilde{X}}$
be the section of $({\cal L}^\tau E)^*\to E$ given by
\begin{equation}\label{7.77'}
\alpha_{\tilde{X}}=A_E\circ \tilde{X}.
\end{equation}
Then, we have the following result.

\begin{proposition}\label{p7.8}
The Lie subalgebroid $(Id,Ti):{\cal L}^{(\tau^{\tau^*})^S}({\cal
L}^{\tau^*}E)\to {\cal L}^{\tau^{\tau^*}}({\cal L}^{\tau^*}E)$,
$i:S\to ({\cal L}^\tau E)^*$ is Lagrangian if and only if the
section $\alpha_{\tilde{X}}$ is $1$-cocycle of the Lie algebroid
$({\cal L}^\tau E, \linebreak
\lcf\cdot,\cdot\rcf^\tau,\rho^\tau).$
\end{proposition}
\begin{proof}
Suppose that $(x^i)$ are local coordinates on an open subset $U$
of $M$ and that $\{e_\alpha\}$ is basis of the vector bundle
$\tau^{-1}(U) \to U$. Then, we will use the following notation:
\begin{itemize}
\item $\{ \tilde{T}_{\alpha}, \tilde{V}_{\alpha}\}$ (respectively,
$\{\tilde{e}_{\alpha}, \bar{e}_{\alpha}\}$) is the local basis of
$\Gamma({\cal L}^{\tau}E)$ (respectively, $\Gamma({\cal
L}^{\tau^*}E)$) considered in Remark \ref{r2.3} (respectively,
Section \ref{seccion3.1}). $\{ \tilde{T}^{\alpha},
\tilde{V}^{\alpha}\}$ (respectively, $\{\tilde{e}^{\alpha},
\bar{e}^{\alpha}\}$) is the dual basis of $\{ \tilde{T}_{\alpha},
\tilde{V}_{\alpha}\}$ (respectively, $\{\tilde{e}_{\alpha},
\bar{e}_{\alpha}\}$).

\item $(x^{i}, y^{\alpha})$ (respectively, $(x^{i}, y_{\alpha};
z^{\alpha}, v_{\alpha})$) are the corresponding local coordinates
on $E$ (respectively, ${\cal L}^{\tau^*}E$).

\item $\rho^{i}_{\alpha}$ and $C_{\alpha \beta}^{\gamma}$ are the
structure functions of $E$ with respect to $(x^i)$ and
$\{e_{\alpha}\}$.

\item $\Theta^{\tilde{X}}: {\cal L}^{\tau}E \to {\cal
L}^{(\tau^{\tau^*})^S}({\cal L}^{\tau^*}E)$ is the isomorphism
between the Lie algebroids ${\cal L}^{\tau}E \to E$ and ${\cal
L}^{(\tau^{\tau^*})^S}({\cal L}^{\tau^*}E) \to S$ (over the map
$\tilde{X}: E \to S = \tilde{X}(E)$) considered in Example
\ref{e7.4}, $(ii)$.

\end{itemize}
Then,
\[
\{\tilde{e}_\alpha^{\bf c},\bar{e}_\alpha^{\bf
c},\tilde{e}_\alpha^{\bf v},\bar{e}_\alpha^{\bf v}\}
\]
is a local basis of $\Gamma({\cal L}^{\tau^{\tau^*}}({\cal
L}^{\tau^*}E)),$ where $\tilde{e}_\alpha^{\bf c}$ and
$\bar{e}_\alpha^{\bf c}$ (respectively, $\tilde{e}_\alpha^{\bf v}$
and $\bar{e}_\alpha^{\bf v}$) are the complete (respectively,
vertical) lift of $\tilde{e}_\alpha$ and $\bar{e}_\alpha$ to
${\cal L}^{\tau^{\tau^*}}({\cal L}^{\tau^*} E).$ Moreover,
\[
\{(\tilde{e}^\alpha)^{\bf v},(\bar{e}^\alpha)^{\bf
v},(\tilde{e}^\alpha)^{\bf c},(\bar{e}^\alpha)^{\bf c}\}
\]
is the dual basis of $\{\tilde{e}_\alpha^{\bf
c},\bar{e}_\alpha^{\bf c},\tilde{e}_\alpha^{\bf
v},\bar{e}_\alpha^{\bf v}\}$.

Now, assume that the local expression of the section $\tilde{X}$
is
\begin{equation}\label{7.77''}
\tilde{X}(x^i,y^\alpha)=(x^i,\tilde{X}^\alpha; y^\alpha,
\tilde{X}'^\alpha).
\end{equation}

Next, we consider the local basis
$\{\tilde{e}_\alpha^{\tilde{X}},\bar{e}_\alpha^{\tilde{X}}\}$ of
$\Gamma({\cal L}^{(\tau^{\tau^*})^S}({\cal L}^{\tau^*}E))$ given
by (\ref{bXtilde}). A direct computation, using (\ref{Defvc}),
(\ref{comverl}), (\ref{dualstru}) and (\ref{bXtilde}), proves that
\begin{equation}\label{comput}
\begin{array}{rcl}
\tilde{e}_\alpha^{\tilde{X}}(\tilde{X}(x^i,y^\gamma))&=&\tilde{e}^{\bf
c}_\alpha(\tilde{X}(x^i,y^\gamma)) +
\rho_\alpha^i(x^i)\displaystyle\frac{\partial
\tilde{X}^\beta}{\partial x^i}_{|(x^i,y^\gamma)}\bar{e}_\beta^{\bf
c}(\tilde{X}(x^i,y^\gamma))\\&&+
C_{\alpha\beta}^\gamma(x^i)y^\beta \tilde{e}_\gamma^{\bf
v}(\tilde{X}(x^i,y^\gamma)) +
\rho_\alpha^i(x^i)\displaystyle\frac{\partial
\tilde{X}'^\beta}{\partial
x^i}_{|(x^i,y^\gamma)}\bar{e}_\beta^{\bf
v}(\tilde{X}(x^i,y^\gamma)),\\
\bar{e}_\alpha^{\tilde{X}}(\tilde{X}(x^i,y^\gamma))&=&\displaystyle\frac{\partial
\tilde{X}^\beta}{\partial
y^\alpha}_{|(x^i,y^\gamma)}\bar{e}_\beta^{\bf
c}(\tilde{X}(x^i,y^\gamma))+ \tilde{e}_\alpha^{\bf
v}(\tilde{X}(x^i,y^\gamma))
\\&& + \displaystyle\frac{\partial \tilde{X}'^\beta}{\partial
y^\alpha}_{|(x^i,y^\gamma)}\bar{e}_\beta^{\bf
v}(\tilde{X}(x^i,y^\gamma)).
\end{array}
\end{equation}

Thus, if $\lambda_{E}$ is the Liouville section of ${\cal
L}^{\tau^*}E \to E$ then, from (\ref{OmEc}), (\ref{7.75''}) and
(\ref{comput}), we obtain that
\[
((Id, Ti) \circ \Theta^{\tilde{X}}, i \circ
\tilde{X})^*(\lambda_{E}^{\bf c}) = (\tilde{X}'^{\alpha} +
\tilde{X}^{\beta} C_{\alpha \gamma}^{\beta}
y^{\gamma})\tilde{T}^{\alpha} + \tilde{X}^{\alpha}
\tilde{V}^{\alpha}.
\]
Therefore, using (\ref{AE}), (\ref{7.77'}) and (\ref{7.77''}), it
follows that
\[
((Id, Ti) \circ \Theta^{\tilde{X}}, i \circ
\tilde{X})^*(\lambda_{E}^{\bf c}) = \alpha_{\tilde{X}}.
\]
Now, since $\Omega_{E}^{\bf c} = (-d^{{\cal L}^{\tau^*}E}
\lambda_{E})^{\bf c} = -d^{{\cal L}^{\tau^{\tau^*}}({\cal
L}^{\tau^*} E)} \lambda_{E}^{\bf c}$ (see Proposition
\ref{complete}), we have that
\begin{equation}\label{7.79}
((Id, Ti) \circ \Theta^{\tilde{X}}, i \circ
\tilde{X})^*(\Omega_{E}^{\bf c}) = -d^{{\cal L}^{\tau}E}
\alpha_{\tilde{X}}.
\end{equation}
Note that the pair $((Id, Ti) \circ \Theta^{\tilde{X}}, i \circ
\tilde{X})$ is a morphism between the Lie algebroids ${\cal
L}^{\tau}E \to E$ and ${\cal L}^{\tau^{\tau^*}}({\cal
L}^{\tau^*}E) \to {\cal L}^{\tau^*}E$.

Consequently, using (\ref{7.79}) and since the rank of the vector
bundle ${\cal L}^{(\tau^{\tau^*})^S}({\cal L}^{\tau^*} E)\to S$ is
$2n$, we deduce the result.
\end{proof}

\setcounter{equation}{0}
\section[Lagrangian submanifolds, Tulczyjew's triple and dynamics]{Lagrangian submanifolds,
Tulczyjew's triple and Eu\-ler-Lagrange (Hamilton) equations}

Let $(E,\lcf\cdot,\cdot\rcf,\rho)$ be a symplectic Lie algebroid
over a manifold $M$ with symplectic section $\Omega$. Then, the
prolongation ${\cal L}^\tau E$ of $E$ over the vector bundle
projection $\tau:E\to M$ is a symplectic Lie algebroid with
symplectic section $\Omega^{\bf c}$, the complete lift of $\Omega$
to ${\cal L}^\tau E$ (see Theorem \ref{t6.4}).

\begin{definition}\label{d8.1}
Let $S$ a submanifold of the symplectic Lie algebroid $E$ and
$i:S\to E$ be the canonical inclusion. $S$ is said to be
Lagrangian if condition (\ref{proS}) holds and the corresponding
Lie subalgebroid $(Id,Ti):{\cal L}^{\tau^S}E\to {\cal L}^\tau E,$
$i:S\to E$ of the symplectic Lie algebroid $({\cal L}^\tau
E,\lcf\cdot,\cdot\rcf^\tau,\rho^\tau)$ is Lagrangian.
\end{definition}

\begin{remark}{\rm Let $M$ be a symplectic manifold with symplectic
$2$-form $\Omega$ and $S$ be a submanifold of $M$. Denote by $i: S
\to M$ the canonical inclusion and by $Ti: TS \to TM$ the tangent
map to $i$. Then, $Ti$ is an injective inmersion and $TS$ is a
submanifold of $TM$, the standard Lie algebroid $\tau_{M}: TM \to
M$ is symplectic, the prolongation ${\cal L}^{\tau_{M}}(TM)$ of
$\tau_{M}: TM \to M$ over $\tau_{M}$ is the standard Lie algebroid
$\tau_{TM}: T(T(M)) \to TM$ and $\Omega^{\bf c}$ is the usual
complete lift of $\Omega$ to $TM$ (see Remark \ref{r6.5}).
Moreover, the Lie subalgebroid $(Id, T(Ti)): {\cal
L}^{\tau_{M}^{TS}}(TM) \to {\cal L}^{\tau_{M}}(TM) = T(TM)$, $Ti:
TS \to TM$ is just the standard Lie algebroid $\tau_{TS}: T(TS)
\to TS$. Thus, $TS$ is a Lagrangian submanifold of $TM$ in the
sense of Definition \ref{d8.1} if and only if $TS$ is a Lagrangian
submanifold of the symplectic manifold $(TM, \Omega^{\bf c})$ in
the usual sense (see Remark \ref{r7.3}). On the other hand, we
have that
\begin{equation}\label{8.78'}
(Ti)^*(\Omega^{\bf c}) = (i^*(\Omega))^{\bf c},
\end{equation}
where $(i^*(\Omega))^{\bf c}$ is the usual complete lift of the
$2$-form $i^*\Omega$ to $TS$. From (\ref{8.78'}), it follows that
$TS$ is a Lagrangian submanifold the symplectic manifold $(TM,
\Omega^{\bf c})$ in the usual sense if and only if $S$ is a
Lagrangian submanifold of the symplectic manifold $(M, \Omega)$ in
the usual sense. Therefore, we conclude that the following
conditions are equivalent:
\begin{enumerate}
\item $S$ is a Lagrangian submanifold of the symplectic manifold
$(M, \Omega)$ in the usual sense.

\item $TS$ is a Lagrangian submanifold of the symplectic manifold
$(TM, \Omega^{\bf c})$ in the usual sense.

\item $TS$ is a Lagrangian submanifold of $TM$ in the sense of
Definition \ref{d8.1}.
\end{enumerate}
}
\end{remark}

From Proposition \ref{p7.5}, we deduce

\begin{corollary}\label{c8.3} Let $(E,\lcf\cdot,\cdot\rcf,\rho)$ be a
symplectic Lie algebroid over a manifold $M$ with symplectic
section $\Omega$ and $X$ be a section of $E$. If $\alpha_X$ is the
section of $E^*$ given by
\[
\alpha_X=i_X\Omega
\]
and $\alpha_X$ is a $1$-cocycle of $E$ then $S=X(M)$ is a
Lagrangian submanifold of $E$.
\end{corollary}

If $(E,\lcf\cdot,\cdot\rcf,\rho)$ is a Lie algebroid over a
manifold $M$, we will denote by $\rho^*(TE^*)\to E$ the pull-back
of the vector bundle $T\tau^*:TE^*\to TM$ over the anchor map
$\rho:E\to TM,$ by ${\cal L}^{\tau^*}E$ the prolongation of $E$
over the vector bundle projection $\tau^*:E^*\to M$ and by
$A_E:\rho^*(TE^*)\to ({\cal L}^{\tau} E)^*$ the isomorphism of
vector bundles considered in Section \ref{seccion5}.

Using Proposition \ref{p7.8}, we have the following result
\begin{corollary}\label{c8.2} Let $(E,\lcf\cdot,\cdot\rcf,\rho)$ be a Lie algebroid
over a manifold $M$ and $\tilde{X}$ be a section of the vector
bundle $\rho^*(TE^*)\to E.$ If $\alpha_{\tilde{X}}$ is the section
of $({\cal L}^\tau E)^*\to E$ given by
\[
\alpha_{\tilde{X}}=A_E\circ \tilde{X}
\]
and $\alpha_{\tilde X}$ is a $1$-cocycle of Lie algebroid $({\cal
L}^\tau E,\lcf\cdot,\cdot\rcf^{\tau},\rho^\tau)$ then
$S=\tilde{X}(E)$ is a Lagrangian submanifold of the symplectic Lie
algebroid ${\cal L}^{\tau^*}E.$
\end{corollary}

Now, let $(E,\lcf\cdot,\cdot\rcf,\rho)$ be a Lie algebroid over a
manifold $M$ and $H:E^*\to \R$ be a Hamiltonian function. If
$\Omega_E$ is the canonical symplectic section of ${\cal
L}^{\tau^*}E,$ then there exists a unique section $\xi_H$ of
${\cal L}^{\tau^*}E\to E^*$ such that
\[
i_{\xi_H}\Omega_E = dH.
\]
Moreover, from Corollary \ref{c8.3}, we deduce that
$S_H=\xi_H(E^*)$ is a Lagrangian submanifold of ${\cal
L}^{\tau^*}E.$

On the other hand, it is clear that there exists a bijective
correspondence $\Psi_H$ between the set of curves in $S_H$ and the
set of curves in $E^*$. In fact, if $c: I \to E^*$ is a curve in
$E^*$ then the corresponding curve in $S_{H}$ is $\xi_{H} \circ c:
I \to S_{H}$.

A curve $\gamma$ in $S_H$,
\[
\gamma:I\to S_H\subseteq {\cal L}^{\tau^*}E\subseteq E\times
TE^*,\;\; t\to (\gamma_1(t),\gamma_2(t))
\]
is said to be {\it admissible } if the curve $\gamma_2:I\to TE^*,
t\to \gamma_2(t),$ is a tangent lift, that is,
\[
\gamma_2(t)=\dot{c}(t),
\]
where $c:I\to E^*$ is the curve in $E^*$ given by $\tau_{E^*}\circ
\gamma_2,$ $\tau_{E^*}:TE^*\to E^*$ being the canonical
projection.

\begin{theorem}
Under the bijection $\Psi_H$, the admissible curves in the
Lagrangian submanifold $S_H$ correspond with the solutions of the
Hamilton equations for $H$.
\end{theorem}
\begin{proof}
Let $\gamma:I\to S_H\subseteq {\cal L}^{\tau^*}E\subseteq E\times
TE^*$ be an admissible curve in $S_H,$
\[
\gamma(t)=(\gamma_1(t),\gamma_2(t)),
\]
for all $t$. Then, $\gamma_2(t)=\dot{c}(t)$, for all $t$, where
$c:I\to E^*$ is the curve in $E^*$ given by $c=\tau_{E^*}\circ
\gamma_2.$

Now, since $\xi_H$ is a section of the vector bundle
$\tau^{\tau^*}:{\cal L}^{\tau^*}E\to E^*$ and $\gamma(I)\subseteq
S_H=\xi_H(E^*),$ it follows that
\begin{equation}\label{cgamma0}
\xi_H(c(t))=\gamma(t),\;\;\;\mbox{ for all } t
\end{equation}
that is, $c=\Psi_H(\gamma).$ Thus, from (\ref{cgamma0}), we obtain
that
\[
\rho^{\tau^*}(\xi_H)\circ c=\gamma_2=\dot{c}, \]
 that is, $c$ is
an integral curve of the vector field $\rho^{\tau^*}(\xi_H)$ and,
therefore, $c$ is a solution of the Hamilton equations associated
with $H$ (see Section \ref{seccion3.3}).

Conversely, assume that $c:I\to E^*$ is a solution of the Hamilton
equations associated with $H$, that is, $c$ is an integral curve
of the vector field $\rho^{\tau^*}(\xi_H)$ or, equivalently,
\begin{equation}\label{ccdot0}
\rho^{\tau^*}(\xi_H)\circ c=\dot{c}.
\end{equation}
Then, $\gamma=\xi_H\circ c$ is a curve in $S_H$ and, from
(\ref{ccdot0}), we deduce that $\gamma$ is admissible.
\end{proof}

 Next, suppose that $L:E\to \R$ is a Lagrangian
function. Then, from Corollary \ref{c8.2}, we obtain that
$S_L=(A_E^{-1}\circ d^{{\cal L}^{\tau}E}L)(E)$ is a Lagrangian
submanifold of the symplectic Lie algebroid ${\cal L}^{\tau^*}E.$

On the other hand, we have a bijective correspondence $\Psi_L$
between the set of curves in $S_L$ and the set of curves in $E$.
In fact, if $\gamma: I \to S_{L}$ is a curve in $S_{L}$ then there
exists a unique curve $c: I \to E$ in $E$ such that
\[
A_{E}(\gamma(t))) = (d^{{\cal L}^{\tau}E}L)(c(t)), \mbox{ for all
} t.
\]
Note that
\[
pr_{1}(\gamma(t)) = (\tau^{\tau})^*(A_E(\gamma(t))) =
(\tau^{\tau})^*((d^{{\cal L}^{\tau}E}L)(c(t))) = c(t),
\makebox[.25cm]{} \mbox{ for all } t,
\]
where $pr_{1}: {\cal L}^{\tau^*}E \subseteq E \times TE^* \to E$
is the canonical projection on the first factor  and
$(\tau^{\tau})^*: ({\cal L}^{\tau}E)^* \to E$ is the vector bundle
projection. Thus,
\[
\gamma(t) = (c(t), \gamma_{2}(t)) \in S_{L} \subseteq {\cal
L}^{\tau^*}E \subseteq E \times TE^*, \makebox[.3cm]{} \mbox{ for
all } t.
\]
A curve $\gamma$ in $S_L$
\[
\gamma:I\to S_L\subseteq {\cal L}^{\tau^*}E\subseteq E\times
TE^*,\;\;\; t\to (c(t),\gamma_2(t))
\]
is said to be {\it admissible} if the curve
\[
\gamma_2:I\to TE^*,\;\;\; t\to \gamma_2(t)
\]
is a tangent lift, that is, $\gamma_2(t)=\dot{c}^*(t),$ where
$c^*:I\to E^*$ is the curve in $E^*$ given by $c^* =
\tau_{E^*}\circ \gamma_2.$

\begin{theorem}
Under the bijection $\Psi_L$, the admissible curves in the
Lagrangian submanifold $S_L$ correspond with the solutions of the
Euler-Lagrange equations for $L.$
\end{theorem}
\begin{proof}
Suppose that $(x^i)$ are local coordinates on $M$ and that
$\{e_\alpha\}$  is a local basis of $\Gamma(E)$. Denote by
$(x^i,y^\alpha)$ (respectively, $(x^i,y_\alpha)$ and
$(x^i,y_\alpha;z^\alpha,v_\alpha))$ the corresponding coordinates
on $E$ (respectively, $E^*$ and ${\cal L}^{\tau^*}E).$ Then, using
(\ref{difpro}), (\ref{TtilT}) and (\ref{AE}), it follows that the
submanifold $S_L$ is characterized by the following equations
\begin{equation}\label{e8.80^0}
y_\alpha=\frac{\partial L}{\partial y^\alpha},\;\;\;
z^\alpha=y^\alpha,\;\;\; v_\alpha=\rho^i_\alpha\frac{\partial
L}{\partial x^i}-C_{\alpha\beta}^\gamma\frac{\partial L}{\partial
y^\gamma}y^\beta,
\end{equation}
for all $\alpha\in  \{1,\dots ,n\}$.

Now, let $\gamma: I \to S_{L}$ be an admissible curve in $S_{L}$
\[
\gamma(t) = (c(t), \gamma_{2}(t)) \in S_{L} \subseteq {\cal
L}^{\tau^*}E \subseteq E \times TE^*, \mbox{ for all } t
\]
and denote by $c^*: I \to E^*$ the curve in $E^*$ satisfying
\begin{equation}\label{e8.80^2}
\gamma_{2}(t) = \dot{c}^*(t), \mbox{ for all } t,
\end{equation}
i.e.,
\begin{equation}\label{e8.80^3}
c^*(t) = \tau_{E^*}(\gamma_{2}(t)), \mbox{ for all } t.
\end{equation}
If the local expressions of $\gamma$ and $c$ are
\[
\gamma(t) = (x^{i}(t), y_{\alpha}(t); z^{\alpha}(t),
v_{\alpha}(t)), \makebox[.3cm]{} c(t) = (x^{i}(t), y^{\alpha}(t)),
\]
then we have that
\begin{equation}\label{e8.80^5}
y^{\alpha}(t) = z^{\alpha}(t), \mbox{ for all } \alpha.
\end{equation}
Moreover, from (\ref{e8.80^2}) and (\ref{e8.80^3}), we deduce that
\begin{equation}\label{e8.80^6}
c^*(t) = (x^i(t), y_{\alpha}(t)), \makebox[.3cm]{} \gamma_{2}(t) =
\displaystyle \frac{dx^i}{dt} \frac{\partial}{\partial
x^i}_{|c^*(t)} + \frac{dy_{\alpha}}{dt} \frac{\partial}{\partial
y_{\alpha}}_{|c^*(t)}.
\end{equation}
Thus,
\begin{equation}\label{e8.80^7}
v_{\alpha}(t) = \displaystyle \frac{dy_{\alpha}}{dt}, \mbox{ for
all } \alpha.
\end{equation}
Therefore, using (\ref{e8.80^0}), (\ref{e8.80^5}),
(\ref{e8.80^6}), (\ref{e8.80^7}) and the fact that $\rho(c(t)) =
(T\tau^*)(\gamma_{2}(t))$, it follows that
\[
\frac{d x^i}{dt}=\rho_\alpha^i y^\alpha,\;\;\;
\frac{d}{dt}(\frac{\partial L}{\partial
y^\alpha})=\rho_\alpha^i\frac{\partial L}{\partial
x^i}-C_{\alpha\beta}^\gamma y^\beta\frac{\partial L}{\partial
y^\gamma},
\]
for all $i$ and $\alpha$, that is, $c$ is a solution of the
Euler-Lagrange equations for $L$.

Conversely, let $c: I \to E$ be a solution of the Euler-Lagrange
equations for $L$ and $\gamma: I \to S_{L}$ be the corresponding
curve in $S_{L}$,
\[
c = \Psi_{L}(\gamma).
\]
Suppose that
\[
\gamma(t) = (c(t), \gamma_{2}(t)) \in {\cal L}^{\tau^*} E
\subseteq E \times TE^*, \mbox{ for all } t
\]
and denote by $c^*: I \to E^*$ the curve in $E^*$ given by
\[
c^* = \tau_{E^*} \circ \gamma_{2}.
\]
If the local expressions of $\gamma$ and $c$ are
\[
\gamma(t) = (x^{i}(t), y_{\alpha}(t); z^{\alpha}(t),
v_{\alpha}(t)), \makebox[.3cm]{} c(t) = (x^{i}(t), y^{\alpha}(t)),
\]
then
\begin{equation}\label{8.80^8}
y^{\alpha}(t) = z^{\alpha}(t), \mbox{ for all }  \alpha
\end{equation}
and the local expressions of $c^*$ and $\gamma_{2}$ are
\[
c^*(t) = (x^{i}(t), y_{\alpha}(t)), \makebox[.3cm]{} \gamma_{2}(t)
= \displaystyle z^{\alpha}(t)\rho_{\alpha}^{i}(x^j(t))
\frac{\partial}{\partial x^i}_{|c^*(t)} + v_{\alpha}(t)
\frac{\partial}{\partial y_{\alpha}}_{|c^*(t)}.
\]
Thus, using (\ref{e8.80^0}) and the fact that $c$ is a solution of
the Euler-Lagrange equations for $L$, we deduce that
\[
\gamma_{2}(t) = \dot{c}^*(t), \mbox{ for all } t,
\]
which implies that $\gamma$ is admissible.
\end{proof}

Now, assume that the Lagrangian function $L:E\to \R$ is
hyperregular and denote by $\omega_L,E_L$ and $\xi_L$ the
Poincar\'{e}-Cartan $2$-section, the energy function and the
Euler-Lagrange section associated with $L$, respectively. Then,
$\omega_L$ is a symplectic  section of the Lie algebroid $({\cal
L}^\tau E, \lcf\cdot,\cdot\rcf^\tau, \rho^\tau)$ and
\[
i_{\xi_L}\omega_L=d^{{\cal L}^{\tau}E}E_L,
\]
(see Section \ref{seccion2.2.2}). Moreover, from Corollary
\ref{c8.3}, we deduce that
\[
S_{\xi_L}=\xi_L(E)
\]
is a Lagrangian submanifold of the symplectic Lie algebroid ${\cal
L}^\tau E.$

On the other hand, it is clear that there exists a bijective
correspondence $\Psi_{S_{\xi_L}}$ between the set of curves in
$S_{\xi_L}$ and the set of curves in $E$.

A curve $\gamma$ in $S_{\xi_L}$
\[
\gamma:I\to S_{\xi_L}\subseteq {\cal L}^\tau E\subseteq E\times
TE,\;\;\; t\to (\gamma_1(t),\gamma_2(t)),
\]
is said to be {\it admissible} if the curve
\[
\gamma_2:I\to TE,\;\;\; t\to \gamma_2(t)
\]
is a tangent lift, that is,
\[
\gamma_2(t)=\dot{c}(t),\mbox{ for all } t,
\]
where $c:I\to E$ the curve in $E$ defined by $c=\tau_E\circ
\gamma_2,$ $\tau_E:TE\to E$ being the canonical projection.

\begin{theorem}
If the Lagrangian $L$ is hyperregular then under the bijection
$\psi_{S_{\xi_L}}$ the admissible curves in the Lagrangian
submanifold $S_{\xi_L}$ correspond with the solutions of the
Euler-Lagrange equations for $L.$
\end{theorem}
\begin{proof}
Let $\gamma:I\to S_{\xi_L}\subseteq {\cal L}^\tau E\subseteq
E\times TE$ be an admissible curve in $S_{\xi_L}$
\[
\gamma(t)=(\gamma_1(t),\gamma_{2}(t)),\mbox{ for all } t.
\]
Then,
\[
\gamma_2(t)=\dot{c}(t),\mbox{ for all } t,
\]
where $c:I\to E$ is the curve in $E$ given by $c=\tau_E\circ
\gamma_2.$

Now, since $\xi_L$ is a section of the vector bundle
$\tau^{\tau}:{\cal L}^\tau E\to E$ and $\gamma(I)\subseteq
S_{\xi_L}=\xi_L(E),$ it follows that
\begin{equation}\label{cgamma1}
\xi_L(c(t))=\gamma(t),\mbox{ for all }t
\end{equation}
that is, $c=\Psi_{S_{\xi_L}}(\gamma).$

Thus, from (\ref{cgamma1}), we obtain that
\[
\rho^\tau(\xi_L)\circ c=\gamma_2=\dot{c},
\]
that is, $c$ is an integral curve of the vector field
$\rho^\tau(\xi_L)$ and, therefore, $c$ is a solution of the
Euler-Lagrange equations associated with $L$ (see Section
\ref{seccion2.2.2}).

Conversely, assume that $c:I\to E$ is a solution of the
Euler-Lagrange equations associated with $L$, that is, $c:I\to E$
is an integral curve of the vector field $\rho^{\tau}(\xi_L)$ or,
equivalently,
\begin{equation}\label{ccdot1}
\rho^\tau(\xi_L)\circ c=\dot{c}.
\end{equation}
Then, $\gamma=\xi_L\circ c$ is a curve in $S_{\xi_L}$ and, from
(\ref{ccdot1}), we deduce that $\gamma$ is admissible.
\end{proof}

 If $L:E\to \R$ is hyperregular then the Legendre
transformation $Leg_L:E\to E^*$ associated with $L$ is a global
diffeomorphism and we may consider the Lie algebroid isomorphism
${\cal L} Leg_L:{\cal L}^\tau E\to {\cal L}^{\tau^*}E$ given by
(\ref{LLegL}) and the Hamiltonian function $H:E^*\to \R$ defined
by $H=E_L\circ Leg_L^{-1}$ (see Section \ref{seccion3.4}).

Thus, we have:
\begin{itemize}
\item The Lagrangian submanifolds $S_L$ and $S_H$ of the symplectic
Lie algebroid ${\cal L}^{\tau^*} E$.
\item The Lagrangian submanifold $S_{\xi_L}$ of the symplectic Lie
algebroid ${\cal L}^\tau E.$
\end{itemize}
\begin{theorem}
If the Lagrangian function $L:E\to \R$ is hyperregular and
$H:E^*\to \R$ is the corresponding Hamiltonian function then the
Lagrangian submanifolds $S_L$ and $S_H$ are equal and
\begin{equation}\label{RelLag}
{\cal L}Leg_L(S_{\xi_L})=S_L=S_H.
\end{equation}
\end{theorem}
\begin{proof}
Using (\ref{Related}), we obtain that
\begin{equation}\label{conrel}
A_E\circ \xi_H\circ Leg_L=A_E\circ {\cal L}Leg_L\circ \xi_L.
\end{equation}
Now, suppose that $(x^i)$ are local coordinates in $M$ and that
$\{e_\alpha\}$ is a local basis of $\Gamma(E)$. Denote by
$(x^i,y^\alpha)$ the corresponding coordinates on $E$ and by
$(x^i,y^\alpha;z^\alpha,v^\alpha)$ (respectively,
$(x^i,y_\alpha;z^\alpha,v_\alpha)$) the corresponding coordinates
on ${\cal L}^\tau E$ (respectively, ${\cal L}^{\tau^*}E)$. Then,
from (\ref{locxiL}), (\ref{LLegloc}) and (\ref{AE}), we deduce
that
\[
(A_E\circ {\cal L}Leg_L\circ
\xi_L)(x^i,y^\alpha)=(x^i,y^\alpha;\rho_\alpha^i\frac{\partial
L}{\partial x^i},\frac{\partial L}{\partial y^\alpha}).\] Thus,
from (\ref{Difftil}), it follows that
\[
(A_E\circ {\cal L}{Leg_L}\circ \xi_L)(x^i,y^\alpha)=d^{{\cal
L}^{\tau}E}L(x^i,y^\alpha),
\]
that is (see (\ref{conrel})),
\[
A_E\circ \xi_H\circ Leg_L=d^{{\cal L}^{\tau}E}L.
\]
Therefore, $S_L=S_H.$

On the other hand, using (\ref{Related}), we obtain that
(\ref{RelLag}) holds.
\end{proof}

 \setcounter{equation}{0}
\section[An application]{An application: Lagrangian submanifolds in prolongations of
Atiyah algebroids and Lagrange (Hamil\-ton)-Poincar\'e equations}
\subsection{Prolongations of the Atiyah algebroid associated with
a principal bundle} Let $\pi:Q\to M$ be a principal bundle with
structural group $G$, $\phi:G\times Q\to Q$ be the free action of
$G$ on $Q$ and $\tau_{Q}|G: TQ/G\to M$ be the Atiyah algebroid
associated with $\pi:Q\to M$ (see Section \ref{seccion2.1.3}).

The tangent action $\phi^T$ of $G$ on $TQ$ is free and thus, $TQ$
is the total space of a principal bundle over $TQ/G$ with
structural group $G$. The  canonical projection $\pi_T:TQ\to TQ/G$
is just the bundle projection.

Now, let ${\cal L}^{(\tau_Q|G)} (TQ/G)$ be the prolongation of the
Atiyah algebroid $\tau_{Q}|G:TQ/G\to M$ by the vector bundle
projection $\tau_Q|G:TQ/G\to M$, and denote by
$(\phi^T)^{T^*}:G\times T^*(TQ)\to T^*(TQ)$ the cotangent lift of
the tangent action $\phi^T:G\times TQ\to TQ.$
\begin{theorem}\label{t9.1} Let $\pi:Q\to M$ be a principal bundle with
structural group $G$ and $\tau_Q|G:TQ/G\to M$ be the Atiyah
algebroid associated with the principal bundle. Then:
\begin{enumerate}
\item The Lie algebroid ${\cal L}^{(\tau_Q|G)}(TQ/G)$ and the
Atiyah algebroid associated with the principal bundle $\pi_T:TQ\to
TQ/G$ are isomorphic.
\item The dual vector bundle to ${\cal L}^{(\tau_Q|G)}(TQ/G)$ is
isomorphic to the quotient vector bundle of $\pi_{TQ}:T^*(TQ)\to
TQ$ by the action $(\phi^T)^{T^*}$ of $G$ on $T^*(TQ)$.
\end{enumerate}
\end{theorem}
\begin{proof}
$(i)$ The Atiyah algebroid associated with the principal bundle
$\pi_T:TQ\to TQ/G$ is the quotient vector bundle
$\tau_{TQ}|G:T(TQ)/G\to TQ/G$ of $\tau_{TQ}:T(TQ)\to TQ$ by the
action $(\phi^T)^T$ of $G$ on $T(TQ).$

On the other hand, we have that the fiber of ${\cal
L}^{(\tau_Q|G)}(TQ/G)$ over $[u_q]\in TQ/G$ is the subspace of
$(TQ/G)_{\pi(q)}\times T_{[u_q]}(TQ/G)$ defined by
\[
\begin{array}{rcl}
{\cal L}^{(\tau_Q|G)} (TQ/G)_{[u_q]}&=& \{([v_q],X_{[u_q]})\in
(TQ/G)_{\pi(q)}\times T_{[u_q]}(TQ/G)/\\
&& (T_q\pi)(v_q)=(T_{[u_q]}(\tau_Q|G))(X_{[u_q]})\}.
\end{array}
\]
Now, we define the morphism $(\pi_T\circ T\tau_{Q}, T\pi_{T})$
between the vector bundles $\tau_{TQ}:T(TQ)\to TQ$ and
$(\tau_Q|G)^{(\tau_Q|G)}:{\cal L}^{(\tau_Q|G)}(TQ/G)\to TQ/G$ over
the map $\pi_T:TQ\to TQ/G$ as follows,
\begin{equation}\label{Ttaupi}
(\pi_T\circ
T\tau_Q,T\pi_T)(X_{u_q})=(\pi_T((T_{u_q}\tau_Q)(X_{u_q})),
(T_{u_q}\pi_T)(X_{u_q}))
\end{equation}
for $X_{u_q}\in T_{u_q}(TQ)$, with $u_q\in T_qQ$.

Since the following diagram

\vspace{-20pt}

\begin{picture}(375,90)(40,20)
\put(180,20){\makebox(0,0){$TQ/G$}}
\put(230,25){$\tau_Q|G$}\put(200,20){\vector(1,0){80}}
\put(310,20){\makebox(0,0){$M=Q/G$}} \put(160,50){$\pi_T$}
\put(180,70){\vector(0,-1){40}} \put(320,50){$\pi$}
\put(310,70){\vector(0,-1){40}} \put(180,80){\makebox(0,0){$TQ$}}
\put(225,85){$\tau_Q$}\put(200,80){\vector(1,0){80}}
\put(310,80){\makebox(0,0){$Q$}}
\end{picture}

is commutative, one deduces that
\[
(\pi_T\circ T\tau_Q,T\pi_T)(X_{u_q})\in {\cal
L}^{(\tau_Q|G)}(TQ/G)_{[u_q]}
\]
and, thus, the map $(\pi_T\circ T\tau_Q,T\pi_T)$ is well-defined.

Next, we will proceed in two steps.

{\it First step: We will prove that the map $(\pi_T\circ T\tau_Q,
T\pi_T)$ induces an isomorphism $\widetilde{(\pi_T\circ T\tau_Q,
T\pi_T)}$   between the vector bundles $\tau_{TQ}|G:T(TQ)/G\to
TQ/G$ and $(\tau_Q|G)^{(\tau_Q|G)}:{\cal L}^{(\tau_Q|G)}(TQ/G)\to
TQ/G.$ }

It is clear that
\[
(\pi_T\circ T\tau_Q,T\pi_T)_{|T_{u_q}(TQ)}:T_{u_q}(TQ)\to {\cal
L}^{(\tau_Q|G)}(TQ/G)_{[u_q]} \] is linear. In addition, this map
is injective. In fact, if $(\pi_T\circ T\tau_Q,
T\pi_T)(X_{u_q})=0$ then $(T_{u_q}\pi_T)(X_{u_q})=0$ and there
exists $\xi\in {\frak g}\cong T_eG$, $e$ being the identity
element of $G$, such that
\begin{equation}\label{Xuq}
X_{u_q}=(T_e(\phi^T)_{u_q})(\xi),
\end{equation}
where $(\phi^T)_{u_q}:G\to TQ$ is defined by
\[
(\phi^T)_{u_q}(g)=(\phi^T)_g(u_q)=(T_q\phi_g)(u_q), \;\;\; \mbox{
for } g\in G.
\]
Therefore, using that $\pi_T((T_{u_q}\tau_Q)(X_{u_q}))=0$, and
hence $(T_{u_q}\tau_Q)(X_{u_q})=0,$ we have that
\[
0=T_e(\tau_Q\circ (\phi^T)_{u_q})(\xi)=(T_e\phi_q)(\xi),
\]
$\phi_q:G\to Q$ being the injective immersion given by
\[
\phi_q(g)=\phi_g(q) = \phi(g, q),\;\;\; \mbox{ for }g\in G.
\]
Consequently, $\xi=0$ and $X_{u_q}=0$ (see (\ref{Xuq})).

We have proved that the linear map $(\pi_T\circ
T\tau_Q,T\pi_T)_{|T_{u_q}(TQ)}$ is injective which implies that
$(\pi_T\circ T\tau_Q,T\pi_T)_{|T_{u_q}(TQ)}$ is a linear
isomorphism (note that the dimensions of the spaces $T_{u_q}(TQ)$
and ${\cal L}^{(\tau_Q|G)}(TQ/G)_{[u_q]})$ are equal).

Furthermore, since the following diagram

\vspace{-20pt}

\begin{picture}(375,90)(40,20)
\put(180,20){\makebox(0,0){$TQ$}}
\put(240,25){$\tau_Q$}\put(200,20){\vector(1,0){80}}
\put(300,20){\makebox(0,0){$Q$}} \put(190,50){$\phi_g^T=T\phi_g$}
\put(180,70){\vector(0,-1){40}} \put(310,50){$\phi_g$}
\put(300,70){\vector(0,-1){40}} \put(180,80){\makebox(0,0){$TQ$}}
\put(240,85){$\tau_Q$}\put(200,80){\vector(1,0){80}}
\put(300,80){\makebox(0,0){$Q$}}

\put(160,80){\vector(-2,-1){40}} \put(160,25){\vector(-2,1){40}}
\put(100,52){\makebox(0,0){$TQ/G$}}\put(130,20){$\pi_T$}\put(130,80){$\pi_T$}
\end{picture}

is commutative, we deduce that $(\pi_T\circ T\tau_Q,T\pi_T)$
induces a morphism $(\widetilde{\pi_T\circ T\tau_Q,T\pi_T})$
between the vector bundles $\tau_{TQ}|G:T(TQ)/G\to TQ/G$ and
$(\tau_Q|G)^{(\tau_Q|G)}: {\cal L}^{(\tau_Q|G)} (TQ/G)\to TQ/G$ in
such a way that the following diagram

\vspace{-20pt}

\begin{picture}(375,90)(40,20)
\put(180,20){\makebox(0,0){$T(TQ)/G$}}
\put(150,50){$(\pi_T)_T$}\put(180,70){\vector(0,-1){40}}
\put(180,80){\makebox(0,0){$T(TQ)$}} \put(210,85){$(\pi_T\circ
T\tau_Q,T\pi_T)$}\put(205,80){\vector(1,0){90}}
\put(340,80){\makebox(0,0){${\cal L}^{(\tau_Q|G)}(TQ/G)$}}

\put(210,20){\vector(3,1){150}}

\put(260,25){$(\widetilde{\pi_T\circ T\tau_Q,T\pi_T})$}
\end{picture}

is commutative, where $(\pi_T)_T:T(TQ)\to T(TQ)/G$ is the
canonical projection. In addition, if $u_q\in T_qQ$ then, since
the map
\[
((\pi_{T})_{T})_{|T_{u_q}(TQ)}:T_{u_q}(TQ)\to (T(TQ)/G)_{[u_q]}
\]

is a linear isomorphism, we conclude that $\widetilde{(\pi_T\circ
T\tau_Q,T\pi_T)}$ is a isomorphism between the vector bundles
$\tau_{TQ}|G:T(TQ)/G\to TQ/G$ and $(\tau_Q|G)^{(\tau_Q|G)}:{\cal
L}^{\tau_Q|G}(TQ/G)\to TQ/G.$

\medskip

{\it Second step: We will prove that the map
\[
\widetilde{(\pi_T\circ T\tau_Q,T\pi_T)}:T(TQ)/G\to {\cal
L}^{(\tau_Q|G)}(TQ/G)
\]
is an isomorphism between the Atiyah algebroid associated with the
principal bundle $\pi_T:TQ\to TQ/G$ and ${\cal
L}^{(\tau_Q|G)}(TQ/G).$}

Let $A:TQ\to {\frak g}$ be a (principal) connection on the
principal bundle $\pi:Q\to M=Q/G.$ We choose a local
trivialization of $\pi:Q\to M$ to be $U\times G$, where $U$ is an
open subset of $M$. Then, $G$ acts on $U\times G$ as follows,
\begin{equation}\label{Acphi}
\phi(g,(x,g'))=(x,gg'),\makebox[1cm]{for } g\in G \mbox{ and }
(x,g')\in U\times G.
\end{equation}

Assume that there are local coordinates $(x^i)$ on $U$ and that
$\{\xi_a\}$ is a basis of ${\frak g}$. Denote by $\{\xi_a^L\}$ the
corresponding left-invariant vector fields on $G$ and suppose that
\[
A(\frac{\partial }{\partial x^i}_{|(x,e)})=A_i^a(x)\xi_a,
\]
for $i\in \{1,\dots ,m\}$ and $x\in U.$ Then, as we know (see
Remark \ref{r2.1''}), the vector fields on $U\times G$
\[
\{e_i=\frac{\partial }{\partial x^i}-A_i^a\xi_a^L,e_b=\xi^L_b\}
\]
define a local basis $\{e'_i, e'_b\}$ of $\Gamma(TQ/G)$, such that
$\pi_T\circ e_i=e'_i\circ\pi$, and $\pi_T\circ e_b=e'_b\circ\pi$.
Thus, one may consider the local coordinates $(x^i,
\dot{x}^i,\bar{v}^b)$ on $TQ/G$ induced by the local basis
$\{e'_i, e'_b\}$ and the corresponding local basis
$\{\tilde{T}_i,\tilde{T}_b,\tilde{V}_i,\tilde{V}_b\}$ of
$\Gamma({\cal L}^{(\tau_Q|G)}(TQ/G))$. From (\ref{Atstfu}) and
(\ref{2.21'}), we have that
\begin{equation} \label{9.3'}
\begin{array}{rcl}\tilde{T}_i[u_q]=(e'_i(\pi(q)),\frac{\partial}{\partial
x^i}_{|[u_q]}),&& \tilde{T}_b[u_q]=(e'_b(\pi(q)), 0),\\
\tilde{V}_i[u_q]=(0,\frac{\partial }{\partial
\dot{x}^i}_{|[u_q]}), && \tilde{V}_b[u_q]=(0,\frac{\partial
}{\partial \bar{v}^b}_{|[u_q]}),
\end{array}
\end{equation}
for $u_ q\in T_qQ$, with $q\in Q.$

On the other hand, using the left-translations by elements of $G,$
one may identify the tangent bundle to $G$, $TG,$ with the product
manifold $G\times {\frak g}$ and, under this identification, the
tangent action of $G$ on $T(U\times G)\cong TU\times TG\cong
TU\times (G\times {\frak g})$ is given by (see (\ref{Acphi}))
\[
\phi^T(g,(v_x,(g',\xi)))=(v_x, (gg',\xi)),
\]
for $g\in G,$ $v_x\in T_xU$ and $(g',\xi)\in G\times {\frak g}.$
Therefore, $T(U\times G)/G\cong TU\times {\frak g}$ and the vector
fields on $T(U\times G)\cong TU\times (G\times {\frak g})$ defined
by

\begin{equation}\label{Xtilbar}
\begin{array}{rcl}\tilde{X}_i=\displaystyle\frac{\partial }{\partial
x^i}-A_i^a\xi_a^L,&& \tilde{X}_b=\xi_b^L,\\
\bar{X}_i=\displaystyle\frac{\partial }{\partial
\dot{x}^i},&&\bar{X}_b=\displaystyle\frac{\partial }{\partial
\bar{v}^b},
\end{array}
\end{equation}

are $\phi^T$-invariant and they define a (local) basis of
$\Gamma(T(TQ)/G).$ Moreover, it follows that
\[
\begin{array}{rcl}
T\tau_Q\circ \tilde{X}_i=e_i\circ\tau_Q,&& T\tau_Q\circ \tilde{X}_b=e_b\circ\tau_Q,\\
T\tau_Q\circ \bar{X}_i =0,&&T\tau_Q\circ \bar{X}_b =0,\\[15pt]
T\pi_T\circ \tilde{X}_i =\displaystyle\frac{\partial }{\partial
x^i}\circ\pi_T, &&T\pi_T\circ \tilde{X}_b
=0,\\[5pt]
T\pi_T\circ \bar{X}_i =\displaystyle\frac{\partial }{\partial
\dot{x}^i}\circ\pi_T, &&T\pi_T\circ \bar{X}_b
=\displaystyle\frac{\partial }{\partial \bar{v}^b}\circ\pi_T.
\end{array}
\]
This implies that
\begin{equation}\label{Relbas}
\begin{array}{rcl}
(\pi_T\circ T\tau_Q,T\pi_T)\circ
\tilde{X}_i=\tilde{T}_i\circ\pi_T,&&
(\pi_T\circ T\tau_Q,T\pi_T)\circ \tilde{X}_b=\tilde{T}_b\circ\pi_T,\\
(\pi_T\circ T\tau_Q,T\pi_T)\circ
\bar{X}_i=\tilde{V}_i\circ\pi_T,&& (\pi_T\circ
T\tau_Q,T\pi_T)\circ \bar{X}_b=\tilde{V}_b\circ\pi_T.
\end{array}
\end{equation}
Furthermore, if $c_{ab}^c$ are the structure constants of ${\frak
g}$ with respect to the basis $\{\xi_a\}$, $B:TQ\oplus TQ\to
{\frak g}$ is the curvature of $A$ and
\[
B(\frac{\partial }{\partial x^i}_{|(x,e)}, \frac{\partial
}{\partial x^j}_{|(x,e)})=B_{ij}^a(x)\xi_a, \makebox[.4cm]{}
\mbox{for} \; \; x \in U,
\]
then, a direct computation proves that,
\[
[\tilde{X}_i,\tilde{X}_j]=-B_{ij}^a\tilde{X}_a,\;\;\;
[\tilde{X}_i, \tilde{X}_a]=c_{ab}^cA_i^b\tilde{X}_c,\;\;\;
[\tilde{X}_a, \tilde{X}_b]=c_{ab}^c\tilde{X}_c,
\]
and the rest of the Lie brackets of the vector fields
$\{\tilde{X}_i, \tilde{X}_a,\bar{X}_i,\bar{X}_a\}$ are zero. Thus,
from (\ref{Atstfu}) and (\ref{tilde}), we conclude that
$(\widetilde{\pi_T\circ T\tau_Q,T\pi_T})$ is a Lie algebroid
isomorphism.

\medskip

$(ii)$ It follows using $(i)$ and the results of Section
\ref{seccion2.1.3} (see Example \ref{ex2.1'}, $(b)$).
\end{proof}

\begin{remark}\label{r9.1'}{\rm As we know ${\cal L}^{\tau_Q}(TQ)\cong T(TQ)$
(see Remark \ref{r2.2}). In addition, if $\{X_i, X_b\}$ is the
local basis of $\Gamma(TQ)={\frak X}(Q)$ given by
\[
X_i=\frac{\partial }{\partial x^i}-A_i^a\xi_a^L,\;\;\; X_b=\xi_b^L
\]
then the corresponding (local) basis of $\Gamma({\cal
L}^{\tau_Q}(TQ))\cong \Gamma(T(TQ))={\frak X}(TQ)$ is
$\{\tilde{X}_i,\tilde{X}_b,\bar{X_i},\bar{X_b}\}$, where
$\tilde{X}_i,\tilde{X}_b,\bar{X}_i$ and $\bar{X}_b$ are the local
vector fields on $T(U\times G)\cong TU\times (G\times {\frak g})$
defined by $(\ref{Xtilbar})$. One may deduce this result using
(\ref{2.21'}), the fact that the anchor map of $\tau_Q:TQ\to Q$ is
the identity and the following equalities
\[
X_i^v=\frac{\partial }{\partial \dot{x_i}},\;\;\;\;
X_b^v=\frac{\partial }{\partial \bar{v}^b},
\]
where $X_i^v$ (respectively $X_b^v$) is the vertical lift of $X_i$
(respectively, $X_b$).}
\end{remark}

If $\pi:Q\to M$ is a principal bundle with structural group $G$
and $\phi:G\times Q\to Q$ is the free action of $G$ on $Q$ then,
as we know (see Section \ref{seccion2.1.3}), the dual vector
bundle to the Atiyah algebroid $\tau_Q|G:TQ/G\to M$ is isomorphic
to the quotient vector bundle of the cotangent bundle
$\pi_Q:T^*Q\to Q$ by the cotangent action $\phi^{T^*}$ of $G$ on
$T^*Q,$ that is, the vector bundles $(\tau_Q|G)^*:(TQ/G)^*\to M$
and $\pi_Q|G:T^*Q/G\to M$ are isomorphic. Since $\phi^{T^*}$ is a
free action, $T^*Q$ is the total space of a principal bundle over
$T^*Q/G$ with structural group $G$. The canonical projection
$\pi_{T^*}:T^*Q\to T^*Q/G$ is just the bundle projection.

Now, denote by ${\cal L}^{(\tau_Q|G)^*}(TQ/G)$ the prolongation of
the Atiyah algebroid  $\tau_Q|G:TQ/G\to M$ by the vector bundle
projection $\pi_Q|G:T^*Q/G\to M$ and by
$(\phi^{T^*})^{T^*}:G\times T^*(T^*Q)\to T^*(T^*Q)$ the cotangent
lift of the cotangent action $\phi^{T^*}:G\times T^*Q\to T^*Q.$

\begin{theorem}\label{t9.2}
Let $\pi:Q\to M$ be a principal bundle with structural group $G$
and $\tau_Q|G:TQ/G\to M$ be the Atiyah algebroid associated with
the principal bundle. Then:
\begin{enumerate}
\item The Lie algebroid ${\cal L}^{(\tau_Q|G)^*}(TQ/G)$ and the
Atiyah algebroid associated with the principal bundle
$\pi_{T^*}:T^*Q\to T^*Q/G$ are isomorphic.
\item
The dual vector bundle to ${\cal L}^{(\tau_Q|G)^*}(TQ/G)$ is
isomorphic to the quotient vector bundle of
$\pi_{T^*Q}:T^*(T^*Q)\to T^*Q$ by the action $(\phi^{T^*})^{T^*}$
of $G$ on $T^*(T^*Q).$
\end{enumerate}
\end{theorem}
\begin{proof}
$(i)$ The Atiyah algebroid associated with the principal bundle
$\pi_{T^*}:T^*Q\to T^*Q/G$ is the quotient vector bundle
$\tau_{T^*Q}|G:T(T^*Q)/G\to T^*Q/G$ of $\tau_{T^*Q}:T(T^*Q)\to
T^*Q $ by the action $(\phi^{T^*})^{T}$ of $G$ on $T(T^*Q)$.

On the other hand, we have that the fiber of ${\cal
L}^{(\tau_Q|G)^*}(TQ/G)$ over $[\alpha_q]\in T^*Q/G,$ with
$\alpha_q\in T^*_qQ$, is the subspace of $(TQ/G)_{\pi(q)}\times
T_{[\alpha_q]}(T^*Q/G)$ defined by
\[
\begin{array}{rcl}{\cal
L}^{(\tau_Q|G)^*}(TQ/G)_{[\alpha_q]}&=&\{([v_q],X_{[\alpha_q]})\in
(TQ/G)_{\pi(q)}\times
T_{[\alpha_q]}(T^*Q/G)/\\&&(T_q\pi)(v_q)=(T_{[\alpha_q]}(\pi_Q|G))(X_{[\alpha_q]})\}.
\end{array}
\]
Now, we define the morphism $(\pi_T\circ T\pi_Q,T\pi_{T^*})$
between the vector bundles $\tau_{T^*Q}:T(T^*Q)\to T^*Q$ and
$(\tau_Q|G)^{(\tau_Q|G)^*}:{\cal L}^{(\tau_Q|G)^*}(TQ/G)\to
T^*Q/G$ over the map $\pi_{T^*}:T^*Q\to T^*Q/G$ as follows,
\begin{equation}\label{Tpiipii}
(\pi_T\circ T\pi_Q,T\pi_{T^*})(X_{\alpha_q}) =(\pi_T(
(T_{\alpha_q}\pi_Q)(X_{\alpha_q})),
(T_{\alpha_q}\pi_{T^*})(X_{\alpha_q}))
\end{equation}
 for
$X_{\alpha_q}\in T_{\alpha_q}(T^*Q),$ with $\alpha_q\in T_q^*Q.$

Since the following diagram

\vspace{-20pt}

\begin{picture}(375,90)(40,20)
\put(180,20){\makebox(0,0){$T^*Q/G$}}
\put(240,25){$\pi_Q|G$}\put(210,20){\vector(1,0){80}}
\put(330,20){\makebox(0,0){$M=Q/G$}} \put(160,50){$\pi_{T^*}$}
\put(180,70){\vector(0,-1){40}} \put(320,50){$\pi$}
\put(310,70){\vector(0,-1){40}}
\put(180,80){\makebox(0,0){$T^*Q$}}
\put(250,85){$\pi_Q$}\put(210,80){\vector(1,0){80}}
\put(310,80){\makebox(0,0){$Q$}}
\end{picture}

is commutative, one deduces that
\[
(\pi_T\circ T\pi_Q,T\pi_{T^*})(X_{\alpha_q})\in {\cal
L}^{(\tau_Q|G)^*}(TQ/G)_{[\alpha_q]}
\]
and, thus, the map $(\pi_T\circ T\pi_Q,T\pi_{T^*})$ is
well-defined.

If $\alpha_q\in T_q^*Q,$ we will denote by
$(\phi^{T^*})_{\alpha_q}:G\to T^*Q$ and by $\phi_q:G\to Q$ the
maps given by
\[
(\Phi^{T^*})_{\alpha_q}(g)=(\Phi^{T^*})_g(\alpha_q)=(T^*\phi_{g^{-1}})(\alpha_q),
\]
\[\phi_q(g)=\phi_g(q)=\phi(g,q),
\]
for $g\in G$. Then, proceeding as in the first step of the proof
of Theorem \ref{t9.1} and using that $\pi_Q\circ
(\phi^{T^*})_{\alpha_q}=\phi_q$ and the fact that the following
diagram

\begin{picture}(375,90)(40,20)
\put(180,20){\makebox(0,0){$T^*Q$}}
\put(250,25){$\pi_Q$}\put(210,20){\vector(1,0){80}}
\put(310,20){\makebox(0,0){$Q$}}
\put(190,50){$\phi_g^{T^*}=T^*\phi_{g^{-1}}$}
\put(180,70){\vector(0,-1){40}} \put(320,50){$\phi_g$}
\put(310,70){\vector(0,-1){40}}
\put(180,80){\makebox(0,0){$T^*Q$}}
\put(250,85){$\pi_Q$}\put(210,80){\vector(1,0){80}}
\put(310,80){\makebox(0,0){$Q$}}

\put(160,80){\vector(-2,-1){40}} \put(160,25){\vector(-2,1){40}}
\put(100,52){\makebox(0,0){$T^*Q/G$}}\put(130,20){$\pi_{T^*}$}\put(130,80){$\pi_{T^*}$}
\end{picture}

is commutative, we deduce that $(\pi_T\circ T\pi_Q,T\pi_{T^*})$
induces an isomorphism $(\widetilde{\pi_T\circ
T\pi_Q,T\pi_{T^*}})$ between the vector bundles
$\tau_{T^*Q}|G:T(T^*Q)/G\to T^*Q/G$ and
$(\tau_Q|G)^{(\tau_Q|G)^*}: {\cal L}^{(\tau_Q|G)^*}(TQ/G)
\linebreak \to T^*Q/G.$

On the other hand, proceeding as in the second step of the proof
of Theorem \ref{t9.1} and using $(\ref{Atstfu})$, (\ref{tilbar})
and (\ref{dualstru}), we conclude that $\widetilde{(\pi_T\circ
T\pi_Q,T\pi_{T^*})}$ is a Lie algebroid isomorphism.

\medskip

$(ii)$ It follows using $(i)$ and the results of Section
\ref{seccion2.1.3} (see Example \ref{ex2.1'}, $(b)$).
\end{proof}

\subsection[Lagrangian submanifolds and Hamilton-Poincar\'{e} equations]{Lagrangian
submanifolds in prolongations of Atiyah algebroids and
Hamilton-Poincar\'{e} equations}\label{seccion9.2}

Let $\pi:Q\to M$ be a principal bundle with structural group $G$,
$\phi:G\times Q\to Q$ be the free action of $G$ on $Q$ and
$\tau_Q|G: TQ/G\to M$ be the Atiyah algebroid associated with the
principal bundle $\pi:Q\to M$. Then, the dual bundle to
$\tau_Q|G:TQ/G\to M$ may be identified with the quotient vector
bundle $\pi_{Q}|G:T^*Q/G\to M$ of the cotangent bundle
$\pi_Q:T^*Q\to Q$ by the cotangent action $\phi^{T^*}$ of $G$ on
$T^*Q.$

Now, denote by $(\pi_T\circ T\pi_Q,T\pi_{T^*}):T(T^*Q)\to {\cal
L}^{(\tau_Q|G)^*}(TQ/G)$ the map given by (\ref{Tpiipii}). Then,
the pair $((\pi_T\circ T\pi_Q,T\pi_{T^*}),\pi_{T^*})$ is a
morphism between the vector bundles $\tau_{T^*Q}:T(T^*Q)\to T^*Q$
and $(\tau_Q|G)^{(\tau_Q|G)^*}: {\cal L}^{(\tau_Q|G)^*}(TQ/G)\to
T^*Q/G.$ We remark that ${\cal L}^{\pi_Q}(TQ)\cong T(T^*Q)$ and,
thus, the Lie algebroids $\tau_{T^*Q}:T(T^*Q)\to T^*Q$ and
$(\tau_Q|G)^{(\tau_Q|G)^*}: {\cal L}^{(\tau_Q|G)^*}(TQ/G)\to
T^*Q/G$ are symplectic (see Section \ref{seccion3.2}).

\begin{theorem}\label{t9.3}
$(i)$ The pair $((\pi_T\circ T\pi_Q,T\pi_{T^*}),\pi_{T^*})$ is a
symplectomorphism between the symplectic Lie algebroids
$\tau_{T^*Q}:T(T^*Q)\to T^*Q$ and $(\tau_Q|G)^{(\tau_Q|G)^*}:
{\cal L}^{(\tau_Q|G)^*}(TQ/G)\to T^*Q/G$. In other words, we have:
\medskip

$(i_a)$ The pair $((\pi_T\circ T\pi_Q,T\pi_{T^*}),\pi_{T^*})$ is a
morphism between the Lie algebroids $\tau_{T^*Q}:T(T^*Q)\to T^*Q$
and $(\tau_Q|G)^{(\tau_Q|G)^*}:{\cal  L}^{(\tau_Q|G)^*}(TQ/G)\to
T^*Q/G.$

$(i_b)$ If $\Omega_{TQ}$ (respectively, $\Omega_{TQ/G}$) is the
canonical symplectic section of $\tau_{T^*Q}:T(T^*Q)\to T^*Q$
(respectively, $(\tau_Q|G)^{(\tau_Q|G)^*}:{\cal
L}^{(\tau_Q|G)^*}(TQ/G)\to T^*Q/G)$ then
\begin{equation}\label{9.6'}
((\pi_T\circ
T\pi_Q,T\pi_{T^*}),\pi_{T^*})^*(\Omega_{TQ/G})=\Omega_{TQ}.
\end{equation}

$(ii)$ Let $h:T^*Q/G\to \R$ be a Hamiltonian function and
$H:T^*Q\to \R$ be the corresponding $G$-invariant Hamiltonian on
$T^*Q$
\begin{equation}\label{9.6''}
H=h\circ \pi_{T^*}.
\end{equation}
If $\xi_H\in \Gamma(T(T^*Q)) \cong {\frak X}(T^*Q)$ (respectively,
$\xi_h\in \Gamma({\cal L}^{(\tau_Q|G)^*}(TQ/G)))$ is the
Hamiltonian section associated with $H$ (respectively, $h$) then
\[
(\pi_T\circ T\pi_Q,T\pi_{T^*})\circ \xi_H=\xi_h\circ \pi_{T^*}.
\]
\end{theorem}
\begin{proof}
$(i)$ We consider the Atiyah algebroid $\tau_{T^*Q}|G:T(T^*Q)/G\to
T^*Q/G$ associated with the principal bundle $\pi_{T^*}:T^*Q\to
T^*Q/G.$ If $\pi_{TT^*}:T(T^*Q)\to T(T^*Q)/G$ is the canonical
projection, we have that the pair $(\pi_{TT^*},\pi_{T^*})$ is a
Lie algebroid morphism (see Section \ref{seccion2.1.3}).

Now, denote by $\widetilde{(\pi_T\circ
T\pi_Q,T\pi_{T^*})}:T(T^*Q)/G\to {\cal L}^{(\tau_Q|G)^*}(TQ/G)$
the isomorphism between the Lie algebroids
$\tau_{T^*Q}|G:T(T^*Q)/G\to T^*Q/G$ and
$(\tau_Q|G)^{(\tau_Q|G)^*}: {\cal L}^{(\tau_Q|G)^*}(TQ/G)\to
T^*Q/G$ considered in the proof of the Theorem \ref{t9.2}. It
follows that
\[
(\pi_T\circ T\pi_Q,T\pi_{T^*})= \widetilde{(\pi_T\circ
T\pi_Q,T\pi_{T^*})} \circ\pi_{TT^*}.
\]
This proves $(i_a).$

Next, we will prove $(i_b).$ If $\lambda_{TQ}$ (respectively,
$\lambda_{TQ/G}$) is the Liouville section of
$\tau_{T^*Q}:T(T^*Q)\to T^*Q$ (respectively,
$(\tau_Q|G)^{(\tau_Q|G)^*}:{\cal L}^{(\tau_Q|G)^*}(TQ/G)\to
T^*Q/G$) then, using (\ref{Lio}) and (\ref{Tpiipii}), we deduce
that
\[
((\pi_T\circ
T\pi_Q,T\pi_{T^*}),\pi_{T^*})^*(\lambda_{TQ/G})=\lambda_{TQ}.
\]
Thus, from $(i_a)$ and since $\Omega_{TQ/G}=-d^{{\cal
L}^{(\tau_Q|G)^*}(TQ/G)}\lambda_{TQ/G}$ and
$\Omega_{TQ}=-d^{T(T^*Q)}\lambda_{TQ},$ we conclude that
\[
((\pi_T\circ
T\pi_Q,T\pi_{T^*}),\pi_{T^*})^*(\Omega_{TQ/G})=\Omega_{TQ}.
\]
$(ii)$ Using $(i_a)$ and $(\ref{9.6''})$, we obtain that
\begin{equation}\label{dHh}
((\pi_T\circ T\pi_Q,T\pi_{T^*}),\pi_{T^*})^*(d^{{\cal
L}^{(\tau_Q|G)^*}(TQ/G)}h)=d^{T(T^*Q)}H.
\end{equation}
Therefore, from (\ref{defxiH}), (\ref{9.6'}) and (\ref{dHh}), it
follows that
\[
\begin{array}{l}
(i_{(\pi_T\circ
T\pi_Q,T\pi_{T^*})(\xi_H(\alpha_q))}\Omega_{TQ/G}(\pi_{T^*}(\alpha_q)))
((\pi_T\circ T\pi_Q,T\pi_{T^*})(X_{\alpha_q}))\\=
(i_{\xi_h(\pi_{T^*}(\alpha_q))}\Omega_{TQ/G}(\pi_{T^*}(\alpha_q)))((\pi_T\circ
T\pi_Q, T\pi_{T^*})(X_{\alpha_q})),
\end{array}
\]
for $\alpha_q\in T^*_qQ$ and $X_{\alpha_q}\in T_{\alpha_q}(T^*Q).$
This implies that
\[
(\pi_T\circ T\pi_Q,
T\pi_{T^*})(\xi_H(\alpha_q))=\xi_h(\pi_{T^*}(\alpha_q)).
\]
\end{proof}

 Now, we prove the following result.
\begin{corollary}\label{c9.4}
Let $h:T^*Q/G\to \R$ be a Hamiltonian function and $H:T^*Q\to \R$
be  the corresponding $G$-invariant Hamiltonian on $T^*Q$. Then,
the solutions of the Hamilton equations for $h$ are just the
solutions of the Hamilton-Poincar\'{e} equations for $H$.
\end{corollary}
\begin{proof}
We will give two proofs.

\medskip
{\it First proof (as a consequence of Theorem \ref{t9.3})}: Let
$\rho$ be the anchor map of the Atiyah algebroid $\tau_Q|G:
TQ/G\to M$ and $\rho^{(\tau_Q|G)^*}:{\cal
L}^{(\tau_Q|G)^*}(TQ/G)\to T(T^*Q/G)$ be the anchor map of the Lie
algebroid ${\cal L}^{(\tau_Q|G)^*}(TQ/G).$ Using Theorem
\ref{t9.3}, we deduce that the vector field $\xi_H\in {\frak
X}(T^*Q)$ is $\pi_{T^*}$-projectable on the vector field
$\rho^{(\tau_Q|G)^*}(\xi_h)\in {\frak X}(T^*Q/G).$ Thus, the
projections, via $\pi_{T^*}$, of the integral curves of $\xi_H$
are the integral curves of the vector field
$\rho^{(\tau_Q|G)^*}(\xi_h).$ But, the integral curves of $\xi_H$
and $\rho^{(\tau_Q|G)^*}(\xi_h)$ are the solutions of the Hamilton
equations for $H$ and $h$, respectively (see Section
\ref{seccion3.3}). Finally, since the projections (via
$\pi_{T^*}$) of the solutions of the Hamilton equations for $H$
are the solutions of the Hamilton-Poincar\'{e} equations for $H$ (see
\cite{CMPR}) the result follows.

\medskip

{\it Second proof (a direct local proof)}: Let $A:TQ\to {\frak g}$
be a (principal) connection on the principal bundle $\pi:Q\to M$
and $B:TQ\oplus TQ\to {\frak g}$ be the curvature of $A$. We
choose a local trivialization of $\pi:Q\to M$ to be $U\times G$,
where $U$ is an open subset of $M$ such that there are local
coordinates $(x^i)$ on $U$. Suppose that $\{\xi_a\}$ is a basis of
${\frak g}$, that $c_{ab}^c$ are the structure constants of
${\frak g}$ with respect to the basis $\{\xi_a\}$ and that $A_i^a$
(respectively, $B_{ij}^a$) are the components of $A$
(respectively, $B$) with respect to the local coordinates $(x^i)$
and the basis $\{\xi_a\}$ (see (\ref{ABcom})).

Denote by $\{e_i, e_a\}$ the local basis of $G$-invariant vector
fields on $Q$ given by (\ref{basis}), by
$(x^i,\dot{x}^i,\bar{v}^a)$ the corresponding local fibred
coordinates on $TQ/G$ and by $(x^i,p_i,\bar{p}_a)$ the (dual)
local fibred coordinates on $T^*Q/G$. Then, using (\ref{Atstfu})
and (\ref{Hameq}), we derive the Hamilton equations for $h$
\[
\frac{dx^i}{dt}=\frac{\partial h}{\partial
p_i},\;\;\;\frac{dp_i}{dt}=-\frac{\partial h}{\partial x^i} +
B_{ij}^a\bar{p}_a\frac{\partial h}{\partial
p_j}-c_{ab}^cA_i^b\bar{p}_c\frac{\partial h}{\partial
\bar{p}_a},\]
\[
\frac{d\bar{p}_a}{dt}=c_{ab}^cA_i^b\bar{p}_c\frac{\partial
h}{\partial p_i}-c_{ab}^c\bar{p}_c\frac{\partial h}{\partial
\bar{p}_b},
\]
which are just the Hamilton-Poincar\'{e} equations associated with the
$G$-invariant Hamiltonian $H$ (see \cite{CMPR}).
\end{proof}

As we know (see Section \ref{seccion3.1}), the local basis $\{e_i,
e_a\}$  of $\Gamma(TQ)$ induces a local basis
$\{\tilde{e}_i,\tilde{e}_a,\bar{e}_i, \linebreak \bar{e}_a\}$ of
$\Gamma({\cal L}^{(\tau_Q|G)^*}(TQ/G))$ and we may consider the
corresponding local coordinates $(x^i,y_i,y_a; \linebreak z^i,
z^a, v_i,v_a)$ on ${\cal L}^{(\tau_Q|G)^*}(TQ/G)$ (see again
Section \ref{seccion3.1}).

Since the Lie algebroids $(\tau_Q|G)^{(\tau_Q|G)^*}:{\cal
L}^{(\tau_Q|G)^*}(TQ/G)\to T^*Q/G$ and $\tau_{T^*Q}|G:T(T^*Q)/G
\linebreak \to T^*Q/G$ are isomorphic and ${\cal
L}^{(\tau_Q|G)^*}(TQ/G)\subseteq TQ/G\times T(T^*Q/G)$, we will
adopt the following notation for the above coordinates
$$(x^i,p_i,\bar{p}_a;\dot{x}^i,\bar{v}^a,\dot{p}_i,\dot{\bar{p}}_a).$$
We recall that $(x^i,\dot{x}_i,\bar{v}^a)$ and
$(x^i,p_i,\bar{p}_a)$ are the local coordinates on $TQ/G$ and
$T^*Q/G$, respectively (see the second proof of Corollary
\ref{c9.4}).

Next, using the coordinates
$(x^i,p_i,\bar{p}_a;\dot{x}^i,\bar{v}^a,\dot{p}_i,\dot{\bar{p}}_a),$
we will obtain the local equations defining the Lagrangian
submanifold $S_h=\xi_h(T^*Q/G)$ of the sympletic Lie algebroid
$({\cal L}^{(\tau_Q|G)^*}(TQ/G),$ $ \Omega_{TQ/G})$, $h:T^*Q/G\to
\R$ being a Hamiltonian function.

Using  (\ref{Atstfu}) and (\ref{xiH0}), we deduce that the local
expression of $\xi_h$ is

\begin{equation}\label{99'}
\begin{array}{rcl}
\xi_h(x^i,p_i,\bar{p}_a)&=&\displaystyle\frac{\partial h}{\partial
p_i}\tilde{e}_i + \displaystyle\frac{\partial h}{\partial
\bar{p}_a}\tilde{e}_a-(\displaystyle\frac{\partial h}{\partial
x^i}-B_{ij}^a\bar{p}_a\displaystyle\frac{\partial h}{\partial
p_j}-c_{bd}^aA_i^b\bar{p}_a\displaystyle\frac{\partial h}{\partial
\bar{p}_d})\bar{e}_i \\&&
+(c_{ab}^cA_i^b\bar{p}_c\displaystyle\frac{\partial h}{\partial
p_i}-c_{ab}^c\bar{p}_c\displaystyle\frac{\partial h}{\partial
\bar{p}_b})\bar{e}_a.
\end{array}
\end{equation}

Thus, the local equations defining the submanifold $S_h$ are
$$\begin{array}{rcl} \bar{v}^a&=&\displaystyle\frac{\partial
h}{\partial \bar{p}_a},\\[8pt]
\dot{x}^i&=&\displaystyle\frac{\partial h}{\partial p_i},\;\;\;
\dot{p_i}=-(\displaystyle\frac{\partial h}{\partial
x^i}-B_{ij}^a\bar{p}_a\displaystyle\frac{\partial h}{\partial
p_j}-c_{bd}^aA_i^b\bar{p}_a\displaystyle\frac{\partial h}{\partial
\bar{p}_d}),\\
\dot{\bar{p}}_a&=&c_{ab}^cA_i^b\bar{p}_c\displaystyle\frac{\partial
h}{\partial {p}_i}- c_{ab}^c\bar{p}_c\displaystyle\frac{\partial
h}{\partial \bar{p}_b} ,\end{array}$$

or, in other words,
\begin{equation}\label{vbara}
\bar{v}^a=\frac{\partial h}{\partial \bar{p}_a},
\end{equation}
\begin{equation}\label{Hpeq}
\begin{array}{l}
\displaystyle\frac{dx^i}{dt}=\frac{\partial h}{\partial
p_i},\;\;\;\;\;\; \displaystyle\frac{dp_i}{dt}=-(\frac{\partial
h}{\partial x^i}-B_{ij}^a\bar{p}_a\displaystyle\frac{\partial
h}{\partial p_j}-c_{bd}^aA_i^b\bar{p}_a\displaystyle\frac{\partial
h}{\partial
\bar{p}_d}),\\\displaystyle\frac{d\bar{p}_a}{dt}=c_{ab}^cA_i^b\bar{p}_c\displaystyle\frac{\partial
h}{\partial {p}_i}-c_{ab}^c\bar{p}_c\frac{\partial h}{\partial
\bar{p}_b} .\end{array}\end{equation}

Eqs. (\ref{vbara}) give the definition of the components of the
{\it locked body angular velocity} (in the terminology of
\cite{BKMM}) and Eqs. (\ref{Hpeq}) are just the {\it
Hamilton-Poincar\'{e} equations} for the $G$-invariant Hamiltonian
$H=h\circ \pi_{T^*}.$

Finally, we will discuss the relation between the solutions of the
Hamilton-Jacobi equation for the Hamiltonians $h$ and $H$.

Suppose that $\alpha \in \Gamma (T^*Q/G)$ is a $1$-cocycle of the
Atiyah algebroid $\tau_Q|G: TQ/G \to M= Q/G$. Then, since the pair
$(\pi_T,\pi)$ is a morphism between the Lie algebroids
$\tau_Q:TQ\to Q$ and $\tau_Q|G:TQ/G\to M=Q/G$ (see Section
\ref{seccion2.1.3}), we deduce that the section
$\tilde{\alpha}\in\Gamma(T^*Q)$ given by
\begin{equation}\label{piddpi}
\tilde{\alpha} = (\pi_T,\pi)^*\alpha,
\end{equation}
is also a 1-cocycle or, in other words, $\tilde{\alpha}$ is a
closed $1$-form on $Q$. It is clear that $\tilde{\alpha}$ is
$G$-invariant. Conversely, if $\tilde{\alpha}$ is a $G$-invariant
closed $1$-form on $Q$ then, using that $(\pi_{T})_{|T_q Q}: T_q Q
\to (TQ/G)_{\pi(q)}$ is a linear isomorphism, for all $q \in Q$,
we deduce that there exists a unique $1$-cocycle $\alpha \in
\Gamma(T^*Q/G)$ of the Atiyah algebroid $\tau_Q|G: TQ/G \to M$
such that (\ref{piddpi}) holds.


\begin{proposition}
There exists a one-to-one correspondence between the solutions of
the Hamilton-Jacobi equation for $h$ and the $G$-invariant
solutions of the Hamilton-Jacobi equation for $H$.
\end{proposition}
\begin{proof}
We recall that a 1-cocycle $\alpha \in \Gamma(T^*Q/G)$
(respectively, $\tilde{\alpha} \in \Gamma (T^*Q)$) is a solution
of the Hamilton-Jacobi equation for $h$ (respectively, $H$) if
$d^{TQ/G}(h\circ \alpha)=0$ (respectively, $d^{TQ}(H\circ
\tilde{\alpha})=0).$

Now, assume that $\alpha$ is a $1$-cocycle of the Atiyah algebroid
$\tau_Q|G: TQ/G \to M$ and denote by $\tilde{\alpha}$ the cocycle
defined by (\ref{piddpi}).  We obtain that
\[
\pi_{T^*}\circ d^{TQ}(H \circ \tilde{\alpha}) = d^{TQ/G}(h \circ
\alpha)\circ \pi.
\]
Thus, using that
\[
\pi_{T^*|T_q^*Q}:T_q^*Q\to (T^*Q/G)_{\pi(q)}
\]
is a linear isomorphism, for all $q\in Q$, we conclude that
\[
d^{TQ/G}(h\circ\alpha)=0\Leftrightarrow d^{TQ}(H\circ
\tilde{\alpha})=0
\]
which proves the result.
\end{proof}


\subsection[Lagrangian submanifolds and Lagrange-Poincar\'{e} equations]{Lagrangian
submanifolds in prolongations of Atiyah algebroids and
Lagrange-Poincar\'{e} equations}\label{seccion9.3}

Let $\pi:Q\to M$ be a principal bundle with structural group $G$,
$\tau_Q|G:TQ/G\to M$ be the Atiyah algebroid associated with the
principal bundle $\pi:Q\to M$ and $\pi_T:TQ\to TQ/G$ be the
canonical projection.

\begin{theorem}
The solutions of the Euler-Lagrange equations for a Lagrangian
$l:TQ/G\to \R$ are the solutions of the Lagrange-Poincar\'{e}
equations for the corresponding $G$-invariant Lagrangian $L$ given
by $L=l\circ \pi_T$.
\end{theorem}
\begin{proof}
Let $A:TQ\to {\frak g}$ be a (principal) connection on the
principal bundle $\pi:Q\to M$ and $B:TQ\oplus TQ\to {\frak g}$ be
the curvature of $A$. We choose a local trivialization of
$\pi:Q\to M$ to be $U\times G,$ where $U$ is an open subset on $M$
such that there are local coordinates $(x^i)$ on $U$. Suppose that
$\{\xi_a\}$ is a basis of ${\frak g}$, that $c_{ab}^c$ are the
structure constants of ${\frak g}$ with respect to the basis
$\{\xi_a\}$ and that $A_i^a$ (respectively, $B_{ij}^a)$ are the
components of $A$ (respectively, $B$) with respect to the local
coordinates $(x^i)$ and the basis $\{\xi_a\}$ (see (\ref{ABcom})).

Denote by $\{e_i, e_a\}$ the local basis of $G$-invariant vector
fields on $Q$ given by (\ref{basis}) and by
$(x^i,\dot{x}^i,\bar{v}^a)$ the corresponding local fibred
coordinates on $TQ/G$. Then, using (\ref{Atstfu}) and
(\ref{2.38'}), we derive the Euler-Lagrange equations for $l$
\[
\frac{\partial l}{\partial x^j}-\frac{d}{dt}(\frac{\partial
l}{\partial \dot{x}^j})=\frac{\partial l}{\partial
\bar{v}^a}(B_{ij}^a\dot{x}^i + c_{db}^aA_j^b\bar{v}^d),\mbox{ for
all } j,
\]
\[
\frac{d}{dt}(\frac{\partial l}{\partial \bar{v}^b})=\frac{\partial
l}{\partial
\bar{v}^a}(c_{db}^a\bar{v}^d-c_{db}^aA_i^d\dot{x}^i),\mbox{ for
all } b,\]

which are just the {\it Lagrange-Poincar\'{e} equations} associated
with the $G$-invariant Lagrangian $L$ (see \cite{CMR}).
\end{proof}

Now, let $A_{TQ/G}:{\cal L}^{(\tau_Q|G)^*}(TQ/G)\equiv
\rho^*(T(T^*Q/G))\to {\cal L}^{(\tau_Q|G)}(TQ/G)^*$ be the
isomorphism between the vector bundles
$pr_1:\rho^*(T(T^*Q/G))\equiv {\cal L}^{(\tau_Q|G)^*}(TQ/G)\to
TQ/G$ and $((\tau_Q|G)^{(\tau_Q|G)})^*:{\cal
L}^{(\tau_Q|G)}(TQ/G)^*\to TQ/G$ defined in Section \ref{seccion5}
(see (\ref{AE})) and $\Omega_{TQ/G}$ be the canonical symplectic
section associated with the Atiyah algebroid $\tau_Q|G:TQ/G\to M.$

Next, we will obtain the local equations defining the Lagrangian
submanifold $S_l=(A_{TQ/G}^{-1}\circ d^{{\cal
L}^{(\tau_Q|G)}(TQ/G)}l)(TQ/G)$ of the symplectic Lie algebroid
$({\cal L }^{{(\tau_Q|G)}^*}(TQ/G),\Omega_{TQ/G})$.

The local basis $\{e_i, e_a\}$ induces a local basis
$\{\tilde{T}_i,\tilde{T}_a,\tilde{V}_i,\tilde{V}_a\}$ of
$\Gamma({\cal L}^{(\tau_Q|G)}(TQ/G))$ (see Remark \ref{r2.3}) and
we may consider the corresponding local coordinates
$(x^i,\dot{x}^i,\bar{v}^a; z^i,z^a,v^i,$ $v^a)$ on ${\cal
L}^{(\tau_Q|G)}(TQ/G)$ (see again Remark \ref{r2.3}). We will
denote by $(x^i,\dot{x}^i,\bar{v}^a; z_i,z_a,v_i,v_a)$ the dual
coordinates on the dual bundle ${\cal L}^{(\tau_Q|G)}(TQ/G)^*$ to
${\cal L}^{(\tau_Q|G)}(TQ/G).$

On  the other hand, we will use the notation of Section
\ref{seccion9.2} for the local coordinates on ${\cal
L}^{(\tau_Q|G)^*}(TQ/G),$ that is, $(x^i,p_i,\bar{p}_a;
\dot{x}^i,\bar{v}^a,\dot{p}_i,\dot{\bar{p}}_a)$.

Then, from (\ref{Atstfu}), (\ref{Difftil}) and (\ref{AE}), we
deduce that
\[
(d^{{\cal
L}^{(\tau_Q|G)}(TQ/G)}l)(x^i,\dot{x}^i,\bar{v}^a)=\frac{\partial
l}{\partial x^i}\tilde{T}^i + \frac{\partial l}{\partial
\dot{x}^i}\tilde{V}^i + \frac{\partial l}{\partial
\bar{v}^a}\tilde{V}^a,\]
\[
\begin{array}{rcl}
A_{TQ/G}^{-1}(x^i,\dot{x}^i,\bar{v}^a;
z_i,z_a,v_i,v_a)&=&(x^i,v_i,v_a;\dot{x}^i,\bar{v}^a,
z_i+B_{ij}^c\dot{x}^jv_c-c_{ab}^cA_i^bv_c\bar{v}^a,\\&&
z_a-c_{ab}^c\bar{v}^bv_c + c_{ab}^cA_j^b\dot{x}^jv_c),
\end{array}
\]
where $\{\tilde{T}^i,\tilde{T}^a,\tilde{V}^i,\tilde{V}^a\}$ is the
dual basis of
$\{\tilde{T}_i,\tilde{T}_a,\tilde{V}_i,\tilde{V}_a\}.$

Thus, the local equations defining the Lagrangian submanifold
$S_l$ of the symplectic Lie algebroid $({\cal
L}^{(\tau_Q|G)^*}(TQ/G), \Omega_{TQ/G})$ are
\[
\begin{array}{l}
p_i=\displaystyle\frac{\partial l}{\partial \dot{x}^i},\;\;\;
\bar{p}_a=\displaystyle\frac{\partial l}{\displaystyle\partial
\bar{v}^a},\\[5pt] \dot{p_i}=\displaystyle\frac{\partial
l}{\partial x^i} + B_{ij}^a \dot{x}^j\displaystyle\frac{\partial
l}{\partial \bar{v}^a} - c_{db}^a
A_i^b\bar{v}^d\displaystyle\frac{\partial l}{\partial
\bar{v}^a},\\[5pt] \dot{\bar{p}}_b=\displaystyle\frac{\partial
l}{\partial \bar{v}^a} (c_{db}^a\bar{v}^d-c_{db}^aA_i^d\dot{x}^i),
\end{array}
\]
or, in other words,
\begin{equation}\label{Defmom}
p_i=\displaystyle\frac{\partial l}{\partial \dot{x}^i},\;\;\;\;\;
\bar{p}_a=\displaystyle\frac{\partial l}{\partial \bar{v}^a},
\end{equation}

\begin{equation}\label{ELequsub}
\begin{array}{l}
\displaystyle\frac{\partial l}{\partial
x^j}-\displaystyle\frac{dp_j}{dt}=\bar{p}_a(B_{ij}^a\dot{x}^i +
c_{db}^aA_j^b\bar{v}^d),\\[8pt] \displaystyle\frac{d
p_b}{dt}=\bar{p}_a(c_{db}^a \bar{v}^d-A_i^dc_{db}^a\dot{x}^i).
\end{array}
\end{equation}

Eqs. (\ref{Defmom}) give the definition of the momenta and Eqs.
(\ref{ELequsub}) are just the Lagrange-Poincar\'{e} equations for the
$G$-invariant Lagrangian $L.$

Now, let $(\pi_T\circ T\tau_Q, T\pi_T):T(TQ)\to {\cal
L}^{(\tau_Q|G)}(TQ/G)$ be the map given by (\ref{Ttaupi}). Then,
the pair $((\pi_T\circ T\tau_Q,T\pi_T),\pi_T)$ is a morphism
between the vector bundles $\tau_{TQ}:T(TQ)\to TQ$ and
$(\tau_Q|G)^{(\tau_Q|G)}:{\cal L}^{(\tau_Q|G)}(TQ/G)\to TQ/G.$ We
recall that the Lie algebroids $\tau_{TQ}:T(TQ)\to TQ$ and
$\tau_Q^{\tau_Q}:{\cal L}^{\tau_Q}(TQ)\to TQ$ are isomorphic.
\begin{theorem}
$(i)$ The pair $((\pi_T\circ T\tau_Q,T\pi_T),\pi_T)$ is morphism
between the Lie algebroids $\tau_{TQ}\cong
\tau_Q^{\tau_Q}:T(TQ)\cong {\cal L}^{\tau_Q}(TQ)\to TQ$ and
$(\tau_Q|G)^{(\tau_Q|G)}:{\cal L}^{(\tau_Q|G)}(TQ/G)\to TQ/G.$

\medskip

$(ii)$ Let $l:TQ/G\to \R$ be a Lagrangian function and $L:TQ\to
\R$ be the corresponding $G$-invariant Lagrangian on $TQ$,
$L=l\circ \pi_T$. If $\omega_l\in \Gamma(\wedge^2({\cal
L}^{(\tau_Q|G)}(TQ/G)^*))$ and $E_l\in C^\infty(TQ/G)$
(respectively, $\omega_L\in \Gamma(\wedge^2({\cal
L}^{\tau_Q}(TQ)^*))\cong \Gamma(\wedge^2(T^*(TQ)))$ and $E_L\in
C^\infty(TQ))$ are the Poincar\'{e}-Cartan $2$-section and the
Lagrangian energy associated with $l$ (respectively, $L$) then
\begin{equation}\label{OmegalL}
((\pi_T\circ T\tau_Q,T\pi_T),\pi_T)^*(\omega_l)=\omega_L,
\end{equation}
\begin{equation}\label{ElL}
E_l\circ \pi_T=E_L.
\end{equation}

\medskip

$(iii)$ If $Leg_l:TQ/G\to T^*Q/G$ (respectively, $Leg_L:TQ\to
T^*Q$) is the Legendre transformation associated with $l$
(respectively, $L$) then
\[
Leg_l\circ \pi_T=\pi_{T^*}\circ Leg_L,
\]
that is, the following diagram is commutative

\begin{picture}(375,85)(40,20)
\put(180,20){\makebox(0,0){$TQ/G$}}
\put(240,25){$Leg_l$}\put(210,20){\vector(1,0){80}}
\put(310,20){\makebox(0,0){$T^*Q/G$}} \put(155,50){$\pi_T$}
\put(180,70){\vector(0,-1){40}} \put(320,50){$\pi_{T^*}$}
\put(310,70){\vector(0,-1){40}} \put(180,80){\makebox(0,0){$TQ$}}
\put(240,85){$Leg_L$}\put(210,80){\vector(1,0){80}}
\put(310,80){\makebox(0,0){$T^*Q$}}
\end{picture}

\end{theorem}
\begin{proof}
$(i)$ We consider the Atiyah algebroid $\tau_{TQ}|G:T(TQ)/G\to
TQ/G$ associated with the principal bundle $\pi_T:TQ\to TQ/G$. If
$\pi_{TT}:T(TQ)\to T(TQ)/G$ is the canonical projection, we have
that the pair $(\pi_{TT},\pi_T)$ is a Lie algebroid morphism (see
Section \ref{seccion2.1.3}).

Now, denote by $\widetilde{(\pi_T\circ T\tau_Q,T\pi_T)}:T(TQ)/G\to
{\cal L}^{(\tau_Q|G)}(TQ/G)$ the isomorphism between the Lie
algebroids $\tau_{TQ}|G:T(TQ)/G\to TQ/G$ and
$(\tau_Q|G)^{(\tau_Q|G)}:{\cal L}^{(\tau_Q|G)}(TQ/G)\to TQ/G$
considered in the proof of Theorem \ref{t9.1}. It follows that
\[
(\pi_T\circ T\tau_Q,T\pi_T)=\widetilde{(\pi_T\circ
T\tau_Q,T\pi_T)}\circ \pi_{TT}.
\]
This proves $(i)$.

\medskip

$(ii)$ From (\ref{2.24'}), (\ref{Relbas}) and Remark \ref{r9.1'},
we deduce that the following diagram is commutative

\begin{picture}(375,90)(40,10)
\put(160,20){\makebox(0,0){${\cal L}^{(\tau_Q|G)}(TQ/G)$}}
\put(250,25){$S^{TQ/G}$}\put(210,20){\vector(1,0){80}}
\put(340,20){\makebox(0,0){${\cal L}^{(\tau_Q|G)}(TQ/G)$}}
\put(90,50){$(\pi_T\circ T\tau_Q,T\pi_T)$}
\put(180,70){\vector(0,-1){40}} \put(320,50){$(\pi_T\circ
T\tau_Q,T\pi_T)$} \put(310,70){\vector(0,-1){40}}
\put(180,80){\makebox(0,0){$T(TQ)$}}
\put(250,85){$S^{TQ}$}\put(210,80){\vector(1,0){80}}
\put(310,80){\makebox(0,0){$T(TQ)$}}
\end{picture}

where $S^{TQ}$ (respectively, $S^{TQ/G}$) is the vertical
endomorphism associated with the Lie algebroid $\tau_Q:TQ\to Q$
(respectively, the Atiyah algebroid $\tau_Q|G:TQ/G\to M$). Thus,
if $\theta_L$ (respectively, $\theta_l$) is the Poincar\'{e}-Cartan
$1$-section associated with $L$ (respectively, $l$) then, using
the first part of the theorem, (\ref{cartan1}) and the fact that
$L=l\circ \pi_T,$ we obtain that
\begin{equation}\label{therel}
((\pi_T\circ T\tau_Q,T\pi_T),\pi_T)^*(\theta_l)=\theta_L.
\end{equation}
Therefore, using again the first part of the Theorem and
(\ref{cartan2}), it follows that
\[
((\pi_T\circ T\tau_Q,T\pi_T),\pi_T)^*(\omega_l)=\omega_L.
\]
On the other hand, from (\ref{2.20'}), (\ref{Relbas}) and Remark
\ref{r9.1'}, we have the following diagram is commutative

\begin{picture}(375,90)(40,10)
\put(180,20){\makebox(0,0){$TQ/G$}}
\put(240,25){$\Delta^{TQ/G}$}\put(210,20){\vector(1,0){80}}
\put(340,20){\makebox(0,0){${\cal L}^{(\tau_Q|G)}(TQ/G)$}}
\put(160,50){$\pi_T$} \put(180,70){\vector(0,-1){40}}
\put(320,50){$(\pi_T\circ T\tau_Q,T\pi_T)$}
\put(310,70){\vector(0,-1){40}} \put(180,80){\makebox(0,0){$TQ$}}
\put(250,85){$\Delta^{TQ}$}\put(210,80){\vector(1,0){80}}
\put(310,80){\makebox(0,0){$T(TQ)$}}
\end{picture}

where $\Delta^{TQ}$ (respectively, $\Delta^{TQ/G}$) is the Euler
section associated with the Lie algebroid $\tau_Q:TQ\to Q$
(respectively, the Atiyah algebroid $\tau_Q|G:TQ/G\to M$).
Consequently, using the first part of the theorem and the fact
that $L=l\circ \pi_T$, we conclude that
\[
E_l\circ \pi_T=E_L.
\]
$(iii)$ From (\ref{LegL}) and (\ref{therel}), we deduce the
result.
\end{proof}

Now we prove the following
\begin{corollary} Let $l:TQ/G\to \R$ be a Lagrangian function and
$L:TQ\to \R$ be the corresponding $G$-invariant Lagrangian on
$TQ$, $ L=l\circ \pi_T.$ Then, $L$ is regular if and only if $l$
is regular.
\end{corollary}
\begin{proof}
The map
\[
(\pi_T\circ T\tau_Q,T\pi_T)_{|T_{v_q}(TQ)}:T_{v_q}(TQ)\to {\cal
L}^{(\tau_Q|G)}(TQ/G)_{[v_q]}
\]
is a linear isomorphism, for all $v_q\in T_qQ$ (see the proof of
Theorem \ref{t9.1}).

On the other hand, $L$ (respectively, $l$) is regular if and only
if $\omega_L$ (respectively, $\omega_l$) is a symplectic section
of $\tau_{TQ}:T(TQ)\to TQ$ (respectively, $(\tau_Q|G)^{(\tau_Q|G)}
: {\cal L}^{(\tau_{Q}|G)}(TQ/G) \to TQ/G$).

Thus, using (\ref{OmegalL}), the result follows.
\end{proof}

Assume that the Lagrangian function $l:TQ/G\to \R$ is regular and
denote by $\xi_l\in \Gamma( {\cal L}^{(\tau_Q|G)}(TQ/G))$ the
Euler-Lagrange section associated with $l$. We recall that $\xi_l$
is characterized by the equation
\[
i_{\xi_l}\omega_l=d^{{\cal L}^{(\tau_Q|G)}(TQ/G)}E_l.
\]

Next, we will obtain the local equations defining the Lagrangian
submanifold  $S_{\xi_l}=\xi_l(TQ/G)$ of the symplectic Lie
algebroid $({\cal L}^{(\tau_Q|G)}(TQ/G),\omega_l^{\bf c})$,
$\omega_l^{\bf c}$ being the complete lift of $\omega_l$.

Let $A:TQ\to {\frak g}$ be a (principal) connection on the
principal bundle $\pi:Q\to M$ and $B:TQ\oplus TQ\to {\frak g}$ be
the curvature of $A$. We choose a local trivialization of
$\pi:Q\to M$ to be $U\times G$, where $U$ is an open subset of $M$
such that there are local coordinates $(x^i)$ on $U$. Suppose that
$\{\xi_a\}$ is a basis of ${\frak g}$, that $c_{ab}^c$ are the
structure constants of ${\frak g}$ with respect to the basis
$\{\xi_a\}$ and that $A_i^a$ (respectively, $B_{ij}^a$) are the
components of $A$ (respectively, $B$) with respect to the local
coordinates $(x^i)$ and the basis $\{\xi_a\}$ (see (\ref{ABcom})).

Denote by $\{e_i, e_a\}$ the local basis of $G$-invariant vector
fields on $Q$ given by (\ref{basis}) and by
$(x^i,\dot{x}^i,\bar{v}^a)$ the corresponding local fibred
coordinates on $TQ/G$. $\{ e_i, e_a\}$ induces a local basis
$\{\tilde{T}_i,\tilde{T}_a,\tilde{V}_i,\tilde{V}_a\}$ of
$\Gamma({\cal L}^{(\tau_Q|G)}(TQ/G))$ (see (\ref{9.3'})) and we
have the corresponding local coordinates
$(x^i,\dot{x}^i,\bar{v}^a; z^i,z^a,$ $v^i,v^a)$ on ${\cal
L}^{(\tau_Q|G)}(TQ/G)$. Since the vector bundles
$\tau_{TQ}|G:T(TQ)/G\to TQ/G$ and $(\tau_Q|G)^{(\tau_Q|G)}:{\cal
L}^{(\tau_Q|G)}(TQ/G)\to TQ/G$ are isomorphic (see Theorem
\ref{t9.1}), we will adopt the following notation for the above
coordinates $$ (x^i,\dot{x}^i,\bar{v}^a;
z^i,z^a,\ddot{x}^i,\dot{\bar{v}}^a).$$

Now, we consider the regular matrix
\[\left(\begin{array}{ll}
W_{ij}&W_{ia}\\W_{ai}&W_{ab}\end{array}\right)=\left(\begin{array}{ll}
\displaystyle\frac{\partial^2 l}{\partial
\dot{x}^i\partial\dot{x}^j} &\displaystyle\frac{\partial^2
l}{\partial
\dot{x}^i\partial\bar{v}^a}\\[5pt]\displaystyle\frac{\partial^2
l}{\partial\bar{v}^a\partial
\dot{x}^i}&\displaystyle\frac{\partial^2 l}{\partial
\bar{v}^a\partial\bar{v}^b}\end{array}\right)
\]
 and denote by
 \[\left(\begin{array}{ll}
W^{ij}&W^{ia}\\W^{ai}&W^{ab}\end{array}\right)\] the inverse
matrix. Then, from  (\ref{Atstfu}) and (\ref{locxiL}), we deduce
that
\begin{equation}\label{xili}
\begin{array}{rcl}\xi_l&=&\dot{x}^i\tilde{T}_i + \bar{v}^a\tilde{T}_a +
(W^{ij}(\xi_l)_j + W^{ia}(\xi_l)_a)\tilde{V}_i\\&& +
(W^{ai}(\xi_l)_i + W^{ab}(\xi_l)_b)\tilde{V}_a,
\end{array}
\end{equation}

where

\begin{equation}\label{xilib}
\begin{array}{rcl}(\xi_l)_i&=&\displaystyle\frac{\partial l}{\partial x^i} -
\displaystyle\frac{\partial^2 l}{\partial x^j\partial
\dot{x}^i}\dot{x}^j + \displaystyle\frac{\partial l}{\partial
\bar{v}^a}(B_{ij}^a\dot{x}^j + c_{bd}^a\bar{v}^bA_i^d), \\[5pt]
(\xi_l)_b&=&-\displaystyle\frac{\partial^2l}{\partial x^i\partial
\bar{v}^b} \dot{x}^i + \displaystyle\frac{\partial l}{\partial
\bar{v}^c} (c_{bd}^cA_i^d\dot{x}^i + c_{ab}^c \bar{v}^a).
\end{array}
\end{equation}

Thus, using the coordinates $(x^i,\dot{x}^i,\bar{v}^a;
z^i,z^a,\ddot{x}^i,\dot{\bar{v}}^a)$, we obtain that the local
equations defining the Lagrangian submanifold $S_{\xi_l}$ are

\begin{equation}\label{zia}
z^i=\dot{x}^i,\;\;\; z^a=\bar{v}^a,\end{equation}
\begin{equation}\label{LPn}
\ddot{x}^i=W^{ij}(\xi_l)_j + W^{ib}(\xi_l)_b,\;\;\;
\dot{\bar{v}}^a=W^{aj}(\xi_l)_j + W^{ab}(\xi_l)_b.
\end{equation}

From (\ref{xilib}) and (\ref{LPn}), we conclude that
\[
\begin{array}{rcl}\displaystyle\frac{\partial l}{\partial x^j} -
\displaystyle\frac{d}{dt}(\displaystyle\frac{\partial l}{\partial
\dot{x}^j}) &=& \displaystyle\frac{\partial l}{\partial
\bar{v}^a}(B_{ij}^a\dot{x}^i + c_{db}^a\bar{v}^dA_j^b), \\[5pt]
\displaystyle\frac{d}{dt}(\displaystyle\frac{\partial l}{\partial
\bar{v}^b})&=& \displaystyle\frac{\partial l}{\partial \bar{v}^a}
(c_{db}^a\bar{v}^d-c_{db}^aA_i^d\dot{x}^i),
\end{array}
\]
which are just the Lagrange-Poincar\'{e} equations associated with the
$G$-invariant Lagrangian $L=l\circ \pi_T$.

\subsection{A particular example: Wong's equations} To illustrate
the theory that we have developed in this section, we will
consider an interesting example, that of {\it Wong's equations}.
Wong's equations arise in at least two different interesting
contexts. The first of these concerns the dynamics of a colored
particle in a Yang-Mills field and the second one is that of the
falling cat theorem (see \cite{M1,M2,M3}; see also \cite{CMR} and
references quoted therein).

Let $(M, g_{M})$ be a given Riemannian manifold, $G$ be a compact
Lie group with a bi-invariant Riemannian metric $\kappa$ and $\pi
: Q \to M$ be a principal bundle with structure group $G$. Suppose
that ${\frak g}$ is the Lie algebra of $G$, that $A: TQ \to {\frak
g}$ is a principal connection on $Q$ and that $B: TQ \oplus TQ \to
{\frak g}$ is the curvature of $A$.

If $q \in Q$ then, using the connection $A$, one may prove that
the tangent space to $Q$ at $q$, $T_{q}Q$, is isomorphic to the
vector space ${\frak g} \oplus T_{\pi(q)}M$. Thus, $\kappa$ and
$g_{M}$ induce a Riemannian metric $g_{Q}$ on $Q$ and we can
consider the kinetic energy $L: TQ \to \R$ associated with
$g_{Q}$. The Lagrangian $L$ is given by
\[
L(v_{q}) = \displaystyle \frac{1}{2}( \kappa_{e}(A(v_{q}),
A(v_{q})) + g_{\pi(q)}((T_{q}\pi)(v_{q}), (T_{q}\pi)(v_{q}))),
\]
for $v_{q} \in T_{q}Q$, $e$ being the identity element in $G$. It
is clear that $L$ is hyperregular and $G$-invariant.

On the other hand, since the Riemannian metric $g_{Q}$ is also
$G$-invariant, it induces a fiber metric $g_{TQ/G}$ on the
quotient vector bundle $\tau_{Q}|G: TQ/G \to M = Q/G$. The reduced
Lagrangian $l: TQ/G \to \R$ is just the kinetic energy of the
fiber metric $g_{TQ/G}$, that is,
\[
l[v_{q}] = \displaystyle \frac{1}{2} (\kappa_{e}(A(v_{q}),
A(v_{q})) + g_{\pi(q)}((T_{q}\pi)(v_{q}), (T_{q}\pi)(v_{q}))),
\]
for $v_{q} \in T_{q}Q$.

We have that $l$ is hyperregular. In fact, the Legendre
transformation associated with $l$ is just the vector bundle
isomorphism $\flat_{g_{TQ/G}}$ between $TQ/G$ and $T^{*}Q/G$
induced by the fiber metric $g_{TQ/G}$. Thus, the reduced
Hamiltonian $h: T^{*}Q/G \to \R$ is given by
\[
h[\alpha_{q}] = l(\flat_{g_{TQ/G}}^{-1}[\alpha_{q}])
\]
for $\alpha_{q} \in T^{*}_{q}Q$.

Now, we choose a local trivialization of $\pi: Q \to M$ to be $U
\times G$, where $U$ is an open subset of $M$ such that there are
local coordinates $(x^{i})$ on $U$. Suppose that $\{\xi_{a}\}$ is
a basis of ${\frak g}$, that $c_{ab}^{c}$ are the structure
constants of ${\frak g}$ with respect to the basis $\{\xi_{a}\}$,
that $A_{i}^{a}$ (respectively, $B_{ij}^{a}$) are the components
of $A$ (respectively, $B$) with respect to the local coordinates
$(x^{i})$ and the basis $\{\xi_{a}\}$ (see (\ref{ABcom})) and that
\[
\kappa_{e} = \kappa_{ab} \xi^{a} \otimes \xi^{b}, \makebox[.4cm]{}
g = g_{ij} dx^{i} \otimes dx^{j},
\]
where $\{\xi^{a}\}$ is the dual basis to $\{\xi_{a}\}$. Note that
since $\kappa$ is a bi-invariant metric on $G$, it follows that
\begin{equation}\label{bi-in}
c_{ab}^{c}\kappa_{cd} = c_{ad}^{c}\kappa_{cb}.
\end{equation}
Denote by $\{e_i, e_a\}$ the local basis of $G$-invariant vector
fields on $Q$ given by (\ref{basis}), by
$(x^i,\dot{x}^i,\bar{v}^a)$ the corresponding local fibred
coordinates on $TQ/G$ and by $(x^i,p_i,\bar{p}_a)$ the (dual)
coordinates on $T^*Q/G$. We have that
\begin{equation}\label{l}
l(x^i,\dot{x}^i,\bar{v}^a) = \displaystyle \frac{1}{2}
(\kappa_{ab} \bar{v}^a \bar{v}^b + g_{ij} \dot{x}^i \dot{x}^j),
\end{equation}
\begin{equation}\label{h}
h(x^i,p_i,\bar{p}_a) = \displaystyle \frac{1}{2} (\kappa^{ab}
\bar{p}_a \bar{p}_b + g^{ij} p_i p_j),
\end{equation}
where $(\kappa^{ab})$ (respectively, $(g^{ij})$) is the inverse
matrix of $(\kappa_{ab})$ (respectively, $(g_{ij})$). Thus, the
Hessian matrix of $l$, $W_{l}$, is
 \[\left(\begin{array}{ll}
g_{ij}&0\\0&\kappa_{ab}\end{array}\right)\] and the inverse matrix
of $W_{l}$ is
\[\left(\begin{array}{ll}
g^{ij}&0\\0&\kappa^{ab}\end{array}\right)\] The local basis
$\{e_i,e_a\}$ induces a local basis
$\{\tilde{e}_i,\tilde{e}_a,\bar{e}_i,\bar{e}_a\}$ of $\Gamma({\cal
L}^{(\tau_Q|G)^*}(TQ/G))$ and we may consider the corresponding
local coordinates
\[
(x^i,p_i,\bar{p}_a; \dot{x}^i,\bar{v}^a,
\dot{p}_i,\dot{\bar{p}}_a)
\]
on ${\cal L}^{(\tau_Q|G)^*}(TQ/G)$ (see Section \ref{seccion9.2}).

From (\ref{bi-in}) and (\ref{h}), we deduce that
\[
\displaystyle c_{ab}^{c} \bar{p}_c \frac{\partial h}{\partial
\bar{p}_{b}} = c_{ab}^{c} \kappa^{db} \bar{p}_c \bar{p}_{d} = 0.
\]
Thus, the local expression of the Hamiltonian section $\xi_{h}$ of
${\cal L}^{(\tau_Q|G)^*}(TQ/G)$ is (see (\ref{99'}) and (\ref{h}))
\[
\xi_h(x^i,p_i,\bar{p}_a) = (g^{ij}p_j) \tilde{e}_i + (\kappa^{bc}
\bar{p}_{c}) \tilde{e}_b-(\displaystyle \frac{1}{2} \frac{\partial
g^{jk}}{\partial x^i} p_j p_k - B_{ij}^a \bar{p}_a g^{jk}p_{k})
\bar{e}_i - (c_{ab}^c A_i^a \bar{p}_c g^{ij}p_{j}) \bar{e}_b.
\]
Therefore, the local equations defining the Lagrangian submanifold
$S_{l} = S_{h} = \xi_{h} (T^*Q/G)$ are

$$\begin{array}{rcl} \dot{x}^i &=& g^{ij}p_{j}, \;\;\; \bar{v}^b =
\kappa^{bc} \bar{p}_{c}, \;\;\; \mbox{for all $i$ and $b$,}
\\[8pt]
\dot{p_i}&=&-\displaystyle \frac{1}{2} \displaystyle
\frac{\partial g^{jk}}{\partial x^{i}} p_j p_k - \bar{p}_{a}
B_{ji}^a g^{jk}p_{k}, \;\;\; \mbox{for all $i$,} \\
\dot{\bar{p}}_b&=& -c_{db}^a A_i^d\bar{p}_a g^{ij}p_j, \;\;\;
\mbox{for all $b$}
\end{array}$$

or, in other words,
\[
\dot{x}^i = g^{ij}p_{j}, \;\;\; \bar{v}^b = \kappa^{bc}
\bar{p}_{c}, \;\;\; \mbox{for all $i$ and $b$,}
\]
\begin{equation}\label{2Wong}
\displaystyle \frac{dp_i}{dt}
 = -\displaystyle \frac{1}{2} \frac{\partial
g^{jk}}{\partial x^{i}} p_j p_k - \bar{p}_{a} B_{ji}^a
g^{jk}p_{k}, \;\;\; \mbox{for all $i$,}
\end{equation}
\begin{equation}\label{1Wong}
\displaystyle \frac{d\bar{p}_b}{dt} = -c_{db}^a A_i^d \bar{p}_a
\dot{x}^{i}, \;\;\; \mbox{for all $b$.}
\end{equation}
Eqs. (\ref{2Wong}) (respectively, Eqs. (\ref{1Wong})) is the {\it
second} (respectively, {\it first}) {\it Wong equation} (see
\cite{CMR}).

On the other hand, the local basis $\{e_i, e_a\}$ induces a local
basis $\{\tilde{T}_i,\tilde{T}_a,\tilde{V}_i,\tilde{V}_a\}$ of
\linebreak $\Gamma({\cal L}^{(\tau_Q|G)}(TQ/G))$ and we may
consider the corresponding local coordinates
\[
(x^i,\dot{x}^i,\bar{v}^a; z^i,z^a,\ddot{x}^i,\dot{\bar{v}}^a)
\]
on ${\cal L}^{(\tau_Q|G)}(TQ/G)$ (see Section \ref{seccion9.3}).

From (\ref{xili}), (\ref{xilib}), (\ref{bi-in}) and (\ref{l}), we
obtain that the Euler-Lagrange section $\xi_{l}$ associated with
$l$ is given by
\[
\begin{array}{rcl}\xi_l (x^k, \dot{x}^k, \bar{v}^c) & = &\dot{x}^i\tilde{T}_i + \bar{v}^b\tilde{T}_b +
g^{ij}(\displaystyle \frac{1}{2} \frac{\partial g_{kl}}{\partial
x^j} \dot{x}^k \dot{x}^l - \frac{\partial g_{jk}}{\partial x^l}
\dot{x}^k \dot{x}^l \\&& + \kappa_{ab} \bar{v}^b B_{jk}^{a}
\dot{x}^k) \tilde{V}_{i} + (c_{ac}^{b} \bar{v}^a A^c_i
\dot{x}^i)\tilde{V}_{b}.
\end{array}\]

Thus, the local equations defining the Lagrangian submanifold
$S_{\xi_{l}} = \xi_{l}(TQ/G)$ of \linebreak ${\cal
L}^{(\tau_{Q}|G)}(TQ/G)$ are
\[
 z^i = \dot{x}^i, \;\;\; z^b =
\bar{v}^b, \;\;\; \mbox{for all $i$ and $b$,}
\]
\begin{equation}\label{x2pti}
\ddot{x}^i = g^{ij}(\displaystyle \frac{1}{2} \frac{\partial
g_{kl}}{\partial x^j} \dot{x}^k \dot{x}^l - \frac{\partial
g_{jk}}{\partial x^l} \dot{x}^k \dot{x}^l + \kappa_{ab} \bar{v}^b
B_{jk}^{a} \dot{x}^k), \;\;\; \mbox{for all $i$},
\end{equation}
\begin{equation}\label{vbarptb}
\dot{\bar{v}}^b = c_{ac}^{b} \bar{v}^a A^c_i \dot{x}^i, \;\;\;
\mbox{for all $b$}.
\end{equation}

Now, put
\[
\bar{p}_{b} = \kappa_{bc}\bar{v}^c, \makebox[.4cm]{} p_i =
g_{ij}\dot{x}^j,
\]
for all $b$ and $i$.

Then, from (\ref{bi-in}) and (\ref{vbarptb}), it follows that
\[
\displaystyle \frac{d\bar{p}_{b}}{dt} = -\bar{p}_{a} c_{db}^{a}
A^d_i \dot{x}^i, \;\;\; \mbox{for all $b$},
\]
which is the first Wong equation. In addition, since $(g^{ij})$ is
the inverse matrix of $(g_{ij})$, we deduce that
\[
\displaystyle \frac{\partial g_{kl}}{\partial x^j} \dot{x}^k
\dot{x}^l = - \frac{\partial g^{kl}}{\partial x^j} p_k p_l.
\]
Therefore, using (\ref{bi-in}) and (\ref{x2pti}), we conclude that
\[
\displaystyle \frac{dp_i}{dt} = -\displaystyle \frac{1}{2}
\frac{\partial g^{jk}}{\partial x^{i}} p_j p_k - \bar{p}_{a}
B_{ji}^a g^{jk}p_{k}, \;\;\; \mbox{for all $i$,}
\]
which is the second Wong equation.

\vspace{.5cm}

{\bf Acknowledgments}. This work has been partially supported by
MICYT (Spain) Grants BFM 2001-2272, BFM 2003-01319 and BFM
2003-02532.

\end{document}